\documentclass{amsart}
\usepackage[margin=3cm]{geometry}
\usepackage{amsmath,amsfonts,amssymb,amsthm}
\usepackage{graphics}
\usepackage{epsfig, esint}
\usepackage{graphics,xcolor}
\usepackage{paralist}
\usepackage{dsfont}
\usepackage[hidelinks]{hyperref}

\numberwithin{equation}{section}

   \definecolor{labelkey}{gray}{.8}
   \definecolor{refkey}{gray}{.8}

\providecommand{\dx}{\, \mathrm{d} x}
\providecommand{\dy}{\, \mathrm{d} y}

\providecommand{\bfa}{{\bf a}}

\newcommand{\supp}{\mathop{\mathrm{supp}}}

\providecommand{\R}{\mathbb{R}}
\providecommand{\per}{{\rm {per}}}

\providecommand{\bfa}{{\bf a}}

\newcommand{\e}{\varepsilon}
\newcommand{\Z}{\mathbb Z}

\newcommand{\step}[1]{\medskip\noindent\textbf{Step #1. }}
\newcommand{\substep}[1]{\medskip\noindent\textit{Substep #1. }}

\newcommand{\ignore}[1]{}

\newtheorem{theorem}{Theorem}

\newtheorem{proposition}[theorem]{Proposition}
\newtheorem{remark}[theorem]{Remark}
\newtheorem{lemma}[theorem]{Lemma}
\newtheorem{corollary}[theorem]{Corollary}
\newtheorem{assumption}{Assumption}
\numberwithin{theorem}{section}

\author[L.~Koch]{Lukas Koch}
\address{University of Sussex\\Mathematics Department\\Falmer Campus\\ BN1 9QH Brighton, United Kingdom}
\email{lukas.koch@sussex.ac.uk}
\author[M. Sch\"affner]{Mathias Sch\"affner}
\address{MLU Halle, Institut f\"ur Mathematik\\
 Theodor Lieser Strasse 5,44227 Halle (Saale), Germany.}
\email{mathias.schaeffner@mathematik.uni-halle.de}

\title[]{Regularity for monotone Operators and applications to homogenization of $p$-Laplace type equations}

\begin{document}
\begin{abstract}In this manuscript, we provide local $L^q$-estimates for the gradient of solutions of a class of quasilinear equations whose principal part lacks strong monotonicity. 

These estimates are used to establish uniform large-scale $L^q$-estimates for the gradient of solutions of degenerate/singular quasilinear equations with oscillating coefficients and large-scale Lipschitz estimates for solutions of non-degenerate equations. 
\end{abstract}

\maketitle

\section{Introduction and main results}

We are interested in regularity properties of solutions of quasilinear elliptic equations in divergence form with rapidly oscillating periodic coefficients. We consider solutions of equations of the form
\begin{align}\label{eq:equationBasic}
\nabla \cdot \bfa\left(\tfrac x \e,\nabla u\right)=\nabla\cdot F \quad\mbox{in $B_1$},
\end{align}
where $\e\in(0,1]$ is a small parameter and $\bfa\colon \R^n \times \R^n \to \R^n$ is periodic and sufficiently smooth in the first variable and a monotone operator satisfying suitable growth and coercivity conditions modeled upon the $p$-Laplace operator in the second variable. More precisely, we assume
\begin{assumption}\label{ass:standard}
 Let $\mu\in[0,1]$, $1<p<\infty$ and $\Lambda\in[1,\infty)$ be given. Suppose that $\bfa:\R^n\times \R^{n}\to\R^{n}$ is a Carath\'eodory function satisfying for all $x\in\R^n$, $\xi_1,\xi_2,\xi\in \R^{n}$ that $\bfa(x,0)=0$ and 
\begin{align}
|\bfa(x,\xi_1)-\bfa(x,\xi_2)|\leq&\Lambda(\mu+|\xi_1|+|\xi_2|)^{p-2}|\xi_1-\xi_2|\label{eq:cont:intro}\\
\Lambda \langle \bfa(x,\xi_1)-\bfa(x,\xi_2),\xi_1-\xi_2\rangle\geq&(\mu+|\xi_1|+|\xi_2|)^{p-2}|\xi_1-\xi_2|^2\label{eq:coer:intro}.
\end{align}
\end{assumption}
Under Assumption~\ref{ass:standard} and supposing that $x\to \bfa(x,z)$ is $\Z^n$-periodic in $x$ and satisfies a suitable continuity property, we show that for every $q\in[p,\infty)$, there exists $C\geq1$ -- which is independent of $\e$ -- such that the following Calderon-Zygmund type estimate holds
\begin{align}\label{thm:Lqx:uniformeps:claimest:intro}
\|\mu+|\nabla u|\|_{\underline L^q(\frac12B)}\leq  C\left(\|\mu+\lvert \nabla u\rvert\|_{\underline L^p(B)}+ \| F\|_{\underline L^\frac{q}{p-1}(B)}^\frac1{p-1}\right),
\end{align}
see Theorem~\ref{thm:Lqx:uniformeps} below for the precise result. In the case of linear elliptic equations, i.e. ${\bfa(x,\xi)=A(x)\xi}$, where $A$ is a sufficiently smooth uniformly elliptic coefficient field and thus $\bfa$ satisfies \eqref{eq:cont:intro} and \eqref{eq:coer:intro} with $\mu=0$ and $p=2$, estimate \eqref{thm:Lqx:uniformeps:claimest:intro} is well-known and goes back to the seminal contributions of Avellaneda and Lin \cite{AL87,AL91}. Since the works \cite{AL87,AL91}, uniform regularity estimates of the form \eqref{thm:Lqx:uniformeps:claimest:intro}, often referred to as \textit{large-scale regularity}, play an important role in the development of homogenization theory. By now there are several far reaching extensions of \cite{AL87,AL91} covering possibly nonlinear equations with random coefficients \cite{AD16,AM16,GNO20}, see also textbooks \cite{Armstrong2019,Shen}. In particular, the work \cite{AD16} contains a suitable version of \eqref{thm:Lqx:uniformeps:claimest:intro} in the case of a nonlinear operator with random coefficients in the case $p=2$. However, to the best of our knowledge there are no general unconditional uniform/large-scale regularity results of the form \eqref{thm:Lqx:uniformeps:claimest:intro} prior to this work in the case $p\neq2$. 

\smallskip

Next we briefly discuss a key difficulty in the case $p\neq2$ and how we overcome this. The classical papers \cite{AL87,AL91}, and all following works on large-scale regularity,  use homogenization theory in some form: Solutions to the equation \eqref{eq:equationBasic} converge, as $\e>0$ tends to zero, to solutions $\overline u$ of     
\begin{align}\label{eq:equationBasichom}
\nabla \cdot \overline \bfa\left(\nabla\overline u\right)=\nabla \cdot F\qquad\mbox{in $B_1$},
\end{align}
where the autonomous homogenized operator $\overline \bfa:\R^n\to\R^n$ can be computed from $\bfa$ with help of certain cell-formulas. This is a classical fact in homogenization theory see e.g.\ the classical paper \cite{Fusco1986} for the nonlinear case of the present work or the textbook references \cite{Braides1998,Zhikov1994} (which are focused on the variational setting). In the linear case, that is $\bfa(x,\xi):=A(x)\xi$ with $A$ uniformly elliptic, the homogenized operator $\overline \bfa$ is given by $\overline \bfa(\xi)=\overline A\xi$ where $\overline A$ is a constant and uniformly elliptic. Hence, classical regularity theory can be applied to solutions of \eqref{eq:equationBasichom}. The regularity of $\overline u$ can be lifted to the solutions of \eqref{eq:equationBasic} by a nontrivial perturbation argument, see \cite{AL87,AL91}. In the nonlinear setting with $p=2$, a similar strategy can also be implemented, see \cite{AD16,AM16} (in the much more difficult case of random homogenization) or \cite{NS19} in the case of elliptic systems. The case $p=2$ is somewhat special as it is easy to show that $\overline \bfa$ satisfies the same growth and monotonicity conditions as $\bfa$ upon enlarging $\Lambda$. Hence, one has classical regularity theory for \eqref{eq:equationBasichom} at the disposal. In contrast, in the case $p\neq2$ it is not clear if the conditions \eqref{eq:cont:intro} and \eqref{eq:coer:intro} are closed under homogenization, that is if the homogenized operator $\overline \bfa$ satisfies the same conditions (upon possibly changing $\Lambda$), see e.g.\ \cite{Fusco1986,DalMaso1990}. As discussed for instance in \cite{Clozeau2023} it is clear in the case $p>2$ that $\overline \bfa$ satisfies
\begin{align}
|\overline \bfa(\xi_1)-\overline \bfa(\xi_2)|\leq&\overline \Lambda(\mu+|\xi_1|+|\xi_2|)^{p-2}|\xi_1-\xi_2|\label{eq:cont:intro:1}\\
\overline \Lambda \langle \overline \bfa(\xi_1)-\overline \bfa(\xi_2),\xi_1-\xi_2\rangle\geq& (\mu+|\xi_1-\xi_2|)^{p-2}|\xi_1-\xi_2|^2\label{eq:coer:intro:1},
\end{align}
for some $\overline \Lambda\in[1,\infty)$, see Section~\ref{sec:hom} below. Clearly, the monotonicity condition \eqref{eq:coer:intro:1} is weaker then \eqref{eq:coer:intro} and it is highly questionable if in general \eqref{eq:coer:intro} holds for $\overline \bfa$, see e.g.\ the discussions in \cite{Clozeau2023}, \cite[p.~398, 399]{CC04} or \cite[Section 7]{DalMaso1990}. In particular, $\overline \bfa$ satisfying \eqref{eq:cont:intro:1} and \eqref{eq:coer:intro:1} is not uniformly elliptic and classical regularity theory (e.g.\ \cite{Giusti2003}) cannot be applied. More can be said appealing to regularity results for non-uniformly elliptic equations. Indeed, equation \eqref{eq:equationBasichom} under \eqref{eq:cont:intro:1} and \eqref{eq:coer:intro:1} with $\mu=1$ can to some extent be covered by the regularity theory with nonstandard growth conditions pioneered by Marcellini~\cite{Mar91}, see \cite{MR21} for a recent overview. In a variational setting, that is $\bfa$ has a convex potential, the currently best results \cite{BS20,BS24} for this problem would yield Lipschitz estimates for solutions of \eqref{eq:equationBasichom} under \eqref{eq:cont:intro:1} and \eqref{eq:coer:intro:1} with $\mu=1$ provided $2\leq p$ and $p<\frac{2(n-1)}{n-3}$ if $n\geq4$ and $F$ is sufficiently regular, see \cite{Clozeau2023}. 

\smallskip

A key result of this contribution is to obtain a Calderon-Zygmund theory for solutions to \eqref{eq:equationBasichom} under the condition \eqref{eq:cont:intro:1} and \eqref{eq:coer:intro:1} \textit{without any restriction} on $p>2$ and $\mu\in[0,1]$. This is the content of Theorem~\ref{thm:Lqx}. Moreover, in the nondegenerate case $\mu=1$, we obtain Lipschitz-regularity for solutions of \eqref{eq:equationBasichom} with sufficiently smooth right-hand side, see Proposition~\ref{lem:HighRegularity} and Remark~\ref{lem:HighRegularity:rem}. This justifies regularity assumptions for the homogenized solution in \cite[Theorem~2.2]{Clozeau2023}. Moreover, we use this to obtain uniform Lipschitz estimates for solutions of \eqref{eq:equationBasic} with $F\equiv0$ under Assumption~\ref{ass:standard} with $\mu=1$, see Theorem~\ref{thm:Linfty:uniformeps}.

Before we discuss the precise results obtained in this manuscript, we mention some further related recent results. In \cite{Wang2019}, the authors provide quantitative homogenization and large-scale H\"older regularity results for $p$-Laplace type equations. The results in \cite{Wang2019} are \textit{conditional},  in the sense that the authors \textit{assume} that the homogenized operator $\overline \bfa$  satisfies the strong monotonicity condition \eqref{eq:coer:intro} which is, as mentioned above, not clear and justified in \cite{Wang2019} only in a specific case. Here, we recover some of the large-scale regularity results of \cite{Wang2019} without any additional assumptions, see Corollary~\ref{cor:largeScaleHolder} and Remark~\ref{rem:wang2019}. In the setting of periodic homogenization with a local defect \cite{Wolf2023} considered the $p$-Laplace operator and proved quantitative convergence results under the assumption that the periodic correctors are non-degenerate. Finally, the work \cite{Clozeau2023} provides, for a restricted range of $p\geq 2$, optimal quantitative stochastic homogenization results assuming \eqref{eq:cont:intro} and \eqref{eq:coer:intro} with $\mu=1$ (extending \cite{FN21} for the case $p=2$).

\subsection{Main results}

The first main result of this manuscript are local Calderon-Zygmund type estimates for solutions to
$$
\nabla \cdot \bfa(x,\nabla u)=\nabla \cdot F,
$$
under rather mild monotonicity and growth conditions on $\xi\mapsto \bfa(x,\xi)$. More precisely, we assume 
\begin{assumption}\label{ass:1} Let $\mu\in[0,1]$, $1<p<\infty$ and $\Lambda\in[1,\infty)$ be given. Suppose that $\bfa:\R^n\times \R^{n}\to\R^{n}$ is a Carath\'eodory function and satisfies 
\begin{align}\label{ass:1:monotone}
\Lambda\langle \bfa(x,\xi_1)-\bfa(x,\xi_2),\xi_1-\xi_2\rangle \geq 
\begin{cases}
|\xi_1-\xi_2|^p &\text{ if } p\geq 2\\
(\mu+\lvert \xi_1\rvert +\lvert \xi_2\rvert)^{p-2}\lvert \xi_1-\xi_2\rvert^2 &\text{ if } p<2
\end{cases}\quad\forall x,\xi_1,\xi_2\in\R^n,
\end{align}
and 
\begin{equation}\label{ass:1:growth}
|\bfa(x,\xi)|\leq \Lambda(\mu+|\xi|)^{p-1}\quad\forall x,\xi\in\R^n.
\end{equation}
\end{assumption}
Under the above assumptions we prove the following
\begin{theorem}\label{thm:Lqx}
Suppose that, for given $1<p<\infty$, Assumption~\ref{ass:1}  is satisfied. Moreover, assume that there is a nondecreasing, continuous modulus of continuity $\omega\colon [0,\infty)\to [0,\infty)$ with $\omega(0)=0$ such that
\begin{align}\label{ass:x}
\lvert \bfa(x,z)-\bfa(y,z)\rvert \leq \omega(\lvert x-y\rvert) \lvert (\mu+\lvert z\rvert)^{p-2}|z|.
\end{align}
Let $F\in L^{p\prime}(\Omega,\R^n)$ and let $u\in W^{1,p}(\Omega)$ be a weak solution to 
\begin{equation}\label{eq:thm:Lqx}
\nabla \cdot \bfa(x,\nabla u)=\nabla \cdot F\qquad\mbox{in $\Omega$.}
\end{equation}
Then, for all $q\in[p,\infty)$ and every $B=B_R(x)\Subset\Omega$ there exists $C=C(\Lambda,n,p,q,\omega,R)$ such that it holds
\begin{align}\label{est:thm:Lqx}
\|\mu+|\nabla u|\|_{\underline L^q(\frac12B)}\leq  C\left(\|\mu+\lvert \nabla u\rvert\|_{\underline L^p(B)}+ \| F\|_{\underline L^\frac{q}{p-1}(B)}^\frac1{p-1}\right).
\end{align}
\end{theorem}

The main novelty of Theorem~\ref{thm:Lqx} concerns the degenerate case $p>2$. Nonlinear Calder\'on-Zygmund estimates of the form \eqref{est:thm:Lqx} are typically derived under the strong monotonicity condition \eqref{eq:coer:intro}, that is
\begin{align}\label{eq:standard}
(\mu+|\xi_1|+|\xi_2|)^{p-2}|\xi_1-\xi_2|^2\lesssim \langle \bfa(x,\xi_1)-\bfa(x,\xi_2),\xi_1-\xi_2\rangle
\end{align}
see the classical references by Iwaniec \cite{Iwaniec83} and DiBenedetto \& Manfredi \cite{DiBM93}. We also refer the reader to \cite{Acerbi2007,DKS12,BWZ07,KZ99} for further results and \cite{Mingione17} for an overview . Here, we prove that \eqref{est:thm:Lqx} remains valid assuming \eqref{ass:1:monotone} instead of \eqref{eq:coer:intro}. Clearly, for $p>2$ condition \eqref{ass:1:monotone} is weaker than \eqref{eq:coer:intro}. To emphasize the difference between \eqref{eq:coer:intro} and \eqref{ass:1:monotone}, assume for the moment that $\bfa$ is smooth and autonomous. In this case \eqref{eq:coer:intro} ensures a quantitative degenerate ellipticity condition for the linearization, namely
$$
|z|^{p-2}|\xi|^2\lesssim \langle D\bfa(z)\xi,\xi\rangle\qquad\forall z,\xi\in\R^{n},
$$
which is a key ingredient in classical linearization based approaches. In contrast, \eqref{ass:1:monotone} does not provide such a quantitative lower bound and is therefore not directly compatible with linearization techniques in nonlinear Calderon-Zygmund theory.

 Finally, we comment that while we did not find a precise reference in the case $p\leq 2$, where Assumption \eqref{ass:1:monotone} is the standard assumption for an operator with $p$-growth, our result is certainly not new in this setting. For example, considering minimisers and $F=0$, a Lipschitz regularity result may be found in \cite[Theorem 2.7]{Fonseca2002}.

Theorem \ref{thm:Lqx} is relevant beyond the application to homogenised operators, we consider in this paper. In particular, for $p\geq 2$ the orthotropic $p$-Laplace operator satisfies Assumption \ref{ass:1}, but not \eqref{eq:standard}. In fact, Theorem \ref{thm:Lqx} applies to a class of degenerate elliptic operators known as Finsler-Laplacians, of which the orthotropic $p$-Laplace operator is one example. We comment further on this application and its relation to existing results in Remark \ref{rem:anistropic1} and Remark \ref{rem:anisotropic2}.


The proof of Theorem \ref{thm:Lqx} is given in Section~\ref{sec:thm:Lqx} and consists of two main steps. In a first step, we consider \eqref{eq:thm:Lqx} with $F\equiv0$ and $\bfa$ being autonomous and prove \eqref{est:thm:Lqx} for every $q<\infty$, see Proposition~\ref{thm:Lq}. This is done via a variant of Moser's iteration method introduced in \cite{Domokos} and extensively used in the regularity theory of non-local operators, starting in \cite{Brasco}. Recently, this method has also been used to study non-uniformly elliptic operators \cite{deFilippis2024a}. In a second step, we deduce Theorem~\ref{thm:Lqx} from the homogeneous and autonomous case by perturbation, following ideas of \cite{CP98} and related works.

\smallskip

The second main result is a version of Theorem~\ref{thm:Lqx} for operators $\bfa$ with rapid periodic oscillations in the space variable.
\begin{theorem}\label{thm:Lqx:uniformeps}
Suppose that, for given $1<p<\infty$, Assumption~\ref{ass:standard} is satisfied. Moreover, assume that there is a nondecreasing, continuous modulus of continuity $\omega\colon [0,\infty)\to [0,\infty)$ with $\omega(0)=0$ such that \eqref{ass:x} is valid and that for every $\xi\in\R^n$ the map $x\mapsto \bfa(x,\xi)$ is $\mathbb Z^n$-periodic. Let $F\in L^{p\prime}(\Omega,\R^n)$ and let $u\in W^{1,p}(\Omega)$ be a weak solution to 
\begin{equation}\label{eq:thm:Lqx:uniformeps}
\nabla \cdot \bfa(\tfrac{x}\e,\nabla u)=\nabla \cdot F\qquad\mbox{in $\Omega$.}
\end{equation}
Then, for all $q\in[p,\infty)$ there exists $C=C(\Lambda,n,p,q,\omega)$ such that it holds for every $B=B_R(x)\Subset\Omega$
\begin{align}\label{thm:Lqx:uniformeps:claimest}
\|\mu+|\nabla u|\|_{\underline L^q(\frac12B)}\leq  C\left(\|\mu+\lvert \nabla u\rvert\|_{\underline L^p(B)}+ \| F\|_{\underline L^\frac{q}{p-1}(B)}^\frac1{p-1}\right).
\end{align}
\end{theorem}

As already described above, Theorem~\ref{thm:Lqx:uniformeps} relies on Theorem~\ref{thm:Lqx} in combination with homogenization methods and a perturbation arguments in the spirit of \cite{AL91,CP98}, see Section~\ref{sec:largescaleCZ}. A key point is that in the situation of Theorem~\ref{thm:Lqx:uniformeps} the homogenized operator $\overline \bfa$ satisfies Assumption~\ref{ass:1} and thus we have -- thanks to Theorem~\ref{thm:Lqx} -- Calderon-Zygmund estimates for the homogenized equation. In Section~\ref{sec:largescaleCZ}, we first establish version of estimate \eqref{thm:Lqx:uniformeps:claimest} valid on large-scales, that is we replace the $L^q$ norm of $\mu+|\nabla u|$ on the left-hand side in \eqref{thm:Lqx:uniformeps:claimest} by quantities of the form $\biggl(\int_{\frac12 B}(\fint_{B_\e(x)}\mu+|\nabla u|^p\dy)^\frac{q}p\dx\biggr)^\frac1q$, see Theorem~\ref{L:nonlinear:largescalereg:lp}. This large-scale regularity is valid without any regularity of $x\mapsto \bfa(x,\xi)$ besides measurability and estimate \eqref{thm:Lqx:uniformeps:claimest} follows by standard regularity theory provided the coefficients are more regular.

 We emphasize that Theorem~\ref{thm:Lqx:uniformeps} is known in the case of linear equations \cite{CP98}  and  nonlinear equations with linear growth (that is $p=2$), see  \cite{AD16}, which even covers the case of random coefficients. However to the best of our knowledge, Theorem~\ref{thm:Lqx:uniformeps} is new already in the case of $p$-Laplacian type equations with oscillating coefficients, that is $\bfa(y,z)=A(y)|z|^{p-2}z$ with $p\in(1,\infty)\setminus\{2\}$ and an uniformly elliptic coefficient matrix $A\in L^\infty(Y,\R^{n\times n})$.

\smallskip

In the non-degenerate / non-singular case, that is Assumption~\ref{ass:standard} holds with $\mu=1$, we obtain Lipschitz estimates for solutions of \eqref{eq:thm:Lqx:uniformeps} with $F\equiv0$.

\begin{theorem}\label{thm:Linfty:uniformeps}
Consider the situation of Theorem~\ref{thm:Lqx:uniformeps} with $\mu=1$ in Assumption~\ref{ass:standard}, and that there exists $\gamma>0$ such that  \eqref{ass:x} is valid with $\mu=1$ and $\omega(s)\leq \Lambda s^\gamma$ for all $s>0$.  Fix $M\in[1,\infty)$. There exists $C_{M}=C_M(\Lambda,M,n,p)\in[1,\infty)$ such that the following is true. Let $u\in W^{1,p}(B)$ with $B=B_R(x_0)\subset \R^n$ be such that
\begin{equation*}
\nabla\cdot \bfa(\tfrac{x}\e,\nabla u)=0\qquad\mbox{in $B$}\quad\mbox{and}\quad\fint_{B}|V_{p}(\nabla u)|^2\dx\leq M,
\end{equation*}
where $V_p(z):=(1+|z|)^\frac{p-2} 2 z$ for all $z\in\R^n$. Then, 
\begin{equation*}
\|V_{p}(\nabla u)\|_{L^\infty(\frac12 B)}^2 \leq C_M  \fint_{B}|V_{p}(\nabla u)|^2\dx.
\end{equation*}
\end{theorem}
Theorem~\ref{thm:Linfty:uniformeps} extends the seminal results of Avellaneda \& Lin \cite{AL87} for linear elliptic systems to the case of nonlinear (non-degenerate \& non-singular) elliptic equations. In the case $p=2$ variations of Theorem~\ref{thm:Linfty:uniformeps} are already proven, see e.g.\  \cite{AM16,NS19}. However, to the best of our knowledge the result is new for $p\in(1,\infty)\setminus\{2\}$. A key ingredient in the proof of Theorem~\ref{thm:Linfty:uniformeps} (provided in Section~\ref{sec:largescaleLip}) is a Lipschitz- and $C^{1,\alpha}$-regularity result for autonomous equations under 
\begin{assumption}\label{ass:homonondeg} Let $1<p<\infty$ and $\Lambda\in[1,\infty)$ be given. Suppose that $\bfa:\R^n\times\R^{n}\to\R^{n}$ satisfies $\bfa(x,0)=0$ for all $x\in\R^n$ and for all $\xi_1,\xi_2\in \R^{n}$ 
\begin{align}\label{ass:1:monotone:strong}
\Lambda\langle \bfa(x,\xi_1)-\bfa(x,\xi_2),\xi_1-\xi_2\rangle \geq 
\begin{cases}
(1+|\xi_1-\xi_2|)^{p-2}|\xi_1-\xi_2|^2 &\text{ if } p\geq 2\\
(1+\lvert \xi_1\rvert +\lvert \xi_2\rvert)^{p-2}\lvert \xi_1-\xi_2\rvert^2 &\text{ if } p<2
\end{cases}
\end{align}
and 
\begin{align}\label{ass:acont}
|\bfa(x,\xi_1)-\bfa(x,\xi_2)|\leq \Lambda \begin{cases}(1+|\xi_1|+|\xi_2|)^{p-2} \lvert \xi_1-\xi_2\rvert&\mbox{if $p\geq 2$}\\(1+|\xi_1-\xi_2|)^{p-2} \lvert \xi_1-\xi_2\rvert&\mbox{if $p<2$}\end{cases}.
\end{align}
\end{assumption}

This is the content of Proposition~\ref{lem:HighRegularity} below and crucially relies on Theorem~\ref{thm:Lqx} and results for \textit{non-uniformly} elliptic equations \cite{BS21,T71}.  The latter is reflected in the appearance of $M$ in Theorem~\ref{thm:Linfty:uniformeps} which also prevents us from considering inhomogeneous equations, see Remark~\ref{rem:rhslipschitz}.    The importance of Assumption~\ref{ass:homonondeg} is that in the setting of Theorem~\ref{thm:Linfty:uniformeps} the homogenized coefficient satisfies Assumption~\ref{ass:homonondeg}. Moreover, it allows us to show that our large-scale regularity results also hold for certain anisotropic operators, including the orthotropic p-Laplace operator with rapidly varying coefficients, see Section~\ref{sec:anisotropic}.

\begin{remark}
The proofs of Theorems~\ref{thm:Lqx:uniformeps} and \ref{thm:Linfty:uniformeps}, explicitly use the periodicity assumption of $x\mapsto \bfa (x,\xi)$. However, the general strategy can also be applied in the case of random homogenization provided there are suitable estimates on the correctors in stochastic homogenization available. In \cite{Clozeau2023}, these estimates are established for specific models satisfying Assumption~\ref{ass:standard} with $\mu=1$ and $2<p<2\frac{n-1}{n-3}$ if $n>3$. Appealing to this, analogous results to Theorem~\ref{thm:Linfty:uniformeps} are established in \cite{CGS} in a random setting under the above restriction on $2<p$.
\end{remark}

\section{Notation and preliminary results}\label{sec:notation}

In this section we introduce our notation and recall some results regarding difference quotient characterisations of Sobolev and Besov spaces.

We denote the open ball of radius $r$ centered at $x\in \R^n$ by $B_r(x)$. Further we write $B_r= B_r(0)$. Given a ball $B=B_r(x_0)$, we denote for $t>0$, $t B = B_{tr}(x_0)$. We set $Y:=(0,1)^n\subset\R^n$. Throughout $\Omega\subset\R^n$ denotes a Lipschitz domain. Given $i\in \{1,\ldots,n\}$ and $h\in \R^n$, we denote $\Omega_{(-|h|)}=\{x\in\Omega \colon d(x,\partial\Omega)>|h|\}$. Given $u\colon \Omega \to \R$, $h\in \R^n$, we denote $\tau_h u(x) = u(x+h)-u(x)$ and $\tau_h^2 u(x) = \frac{u(x+h)-2 u(x)+u(x-h)} 2$, whenever these expressions make sense.

In any estimate, if the right-hand side is not well-defined, we understand it to be infinite and the statement to be empty.

\subsection{Function spaces}

Given $p\geq1$, we denote by $L^p(\Omega)$ and $W^{1,p}(\Omega)$ the standard Lebesgue and Sobolev spaces on $\Omega$. $W^{1,p}_{per,0}(Y)$ is the space of $W^{1,p}$ functions $u$ that are periodic on $Y$ with ${\fint_Y u \dx = 0}$.  We find it convenient to denote for $p\geq 1$,
\begin{align*}
\|f\|_{\underline L^p(U)} = \left(\fint_U \lvert f\rvert^{p}\dx\right)^\frac 1 p.
\end{align*}
Moreover, for $u\in L^1(\Omega)$, we write $(u)_\Omega = \fint_\Omega u\dx$. Finally, for $\alpha\in (0,1)$, $C^{0,\alpha}(\Omega)$ denotes the space of $\alpha$-H\"older-continuous functions on $\Omega$.

As we will work with difference quotients, the following semi-norm will be useful to us: Given a convex domain $\Omega$, $s\in (0,1]$ and $p\geq 1$, we set
\begin{align}\label{eq:defbesov}
[u]_{s,p,\Omega}:= \sup_{0\neq h\in \R^n} \left(\int_{\Omega_h} \frac{|\tau_h u|^p}{|h|^{sp}}\dx\right)^\frac 1 p\quad \text{where}\quad \Omega_h = \{x\in \Omega \colon x+th\in \Omega \text{ for all } t\in [0,1]\}.
\end{align}
If $p=+\infty$, the expression needs to be modified in the obvious way.
The quantity $[\cdot]_{s,p,\Omega}$ is well-known. In fact, for $s=1$, $[u]_{s,p,\Omega}$ is equivalent to $\|\nabla u\|_{L^p(\Omega)}$, while for $s<1$ it corresponds to the semi-norm in the Besov space $B^{s,p}_\infty(\Omega)$.  We recommend \cite{Triebel} as a reference to the interested reader for the general theory of these function spaces. Here, we recall only a number of relevant facts and provide elementary proofs for some of them.

We begin by recalling the relation between difference quotients and $L^q$-gradient norms.
\begin{lemma}[\texorpdfstring{\cite[Theorem 10.55]{Leoni}}{}]
Let $q>1$ and $B=B_r(x)$. If $u\in W^{1,q}(B)$, then for any $h\in \R^n$,
\begin{align}\label{eq:fact4:1}
\int_{B_{h}} \frac{\lvert \tau_h u\rvert^q}{\lvert h\rvert^q}\dx\leq \int_B \left\lvert \nabla u\cdot \frac{h}{|h|}\right\rvert^q\dx.
\end{align}
If $u\in L^q(B)$, then for any $h\in \R^n$,
\begin{align}\label{eq:gradientToDifferenceQuotient}
    \int_B |\nabla u\cdot \frac{h}{|h|}|^q\leq \liminf_{h\to 0^+} \int_{B_h} \frac{|\tau_h u|^q}{|h|^q}\dx.
\end{align}
\end{lemma}
The following lemma relates second-order difference quotients and first-order Sobolev spaces.
\begin{lemma}[\texorpdfstring{\cite[Lemma 2.4.]{Garain}}{}]\label{lem:differenceQuotient}
Let $q>1$, $\rho,h_0\in (0,1]$ and $\e\in (0,q)$. Assume $u\in L^q(B_{\rho+7h_0})$. Then, there is ${c=c(n,q)>0}$ such that
\begin{align}\label{eq:fact4}
\|\nabla u\|_{L^q(B_\rho)}\leq \frac{c q^2}{\e(q-\e)}\biggl(\frac{\|u\|_{L^q(B_{\rho+6h_0})}}{\lvert h_0\rvert^{1+\frac \e q}}+\biggl(\sup_{0<\lvert h\rvert<h_0} \int_{B_{\rho+5h_0}} \frac{\lvert \tau_h^2 u\rvert^q}{\lvert h\rvert^{q+\e}}\dx\biggr)^\frac1q\biggr).
\end{align}
\end{lemma}

The following embedding theorem for Besov spaces is well-known.
\begin{lemma}\label{L:besovembedding}
Let $p>1$, $s\in(0,1)$ and let $\Theta\in(1,\infty]$ be given by 
$$
\frac{1}{\Theta}=\max\{1-(1-s)\frac{p}{n},0\}
$$
with the understanding $\frac{1}{\infty}=0$. Then there is $c= c(n,p,s)>0$ such that for every $B_R(x_0)\subset \R^{n}$
\begin{align}\label{eq:fact2}
[u]_{s,\Theta p,B_R(x_0)}\leq cR^{1-s}\|\nabla u\|_{L^{p}(B_R(x_0))}.
\end{align}
\end{lemma}
Lemma~\ref{L:besovembedding} can be deduced from e.g.~\cite[Section 3.3.1]{Triebel}. However, since we did not find a reference for the exact estimate \eqref{eq:fact2}, we provide an elementary proof of this fact without appealing to the theory of Besov spaces.
\begin{proof}
Throughout the proof we write '$\lesssim$' if '$\leq$' holds up to a positive multiplicative constant that depends only $n,p$ and $\Theta$. By scaling and translation, it suffices to consider $B:=B_1(0)$. Set $v:=u-(u)_{B}\in W^{1,p}(B)$ and consider an extension $v\in W^{1,p}(\R^n)$ satisfying
\begin{equation}\label{est:vextension}
\|v\|_{W^{1,p}(\R^{n})}\lesssim\|v\|_{W^{1,p}(B_1)}\lesssim \|\nabla u\|_{L^p(B_1)},
\end{equation} 
where the last inequality follows from the definition of $v$ and Poincar\'e inequality.

\step 1 The case $\Theta<\infty$. For $r>0$ and $z\in\Z^{n}$, consider $B_{z,r}:=B_{\sqrt{n}r}(rz)$. It is easy to check that
\begin{equation}\label{ass:coveringbzr}
1\leq \sum_{z\in\Z^{d}}\mathds 1_{B_{z,r}}\lesssim1.
\end{equation}
Let $h\in\R^{n}$ be given. The definition of $v$ and triangle inequality yield 
\begin{align}\label{interpolatethu}
\|\tau_h u\|_{L^{\Theta p}((B_1(0))_h)}=&\|\tau_h v\|_{L^{\Theta p}((B_1(0))_h)}\notag\\
\leq& \biggl(\sum_{z\in\Z^{n}}\|\tau_h v-(\tau_h v)_{B_{z,r}}\|_{L^{\Theta p}(B_{z,r})}^{\Theta p}\biggr)^\frac{1}{\Theta p}+\biggl(\sum_{z\in\Z^{n}}|B_{z,r}||(\tau_h v)_{B_{z,r}}|^{\Theta p}\biggr)^\frac{1}{\Theta p}.
\end{align}
For the first term on the right-hand side of \eqref{interpolatethu}, we use Poincar\'e inequality 
$$
\|\tau_h v-(\tau_h v)_{B_{z,r}}\|_{L^{\Theta p}(B_{z,r})}\lesssim r^{1+\frac{n}{\Theta p}-\frac{n}{p}}\|\nabla \tau_h v\|_{L^p(B_{z,r})}
$$
together with the discrete $\|\cdot\|_{\ell^\Theta}\leq \|\cdot\|_{\ell^1}$ inequality and \eqref{ass:coveringbzr} in the form
\begin{align*}
\biggl(\sum_{z\in\Z^{n}}\|\tau_h v-(\tau_h v)_{B_{z,r}}\|_{L^{\Theta p}(B_{z,r})}^{\Theta p}\biggr)^\frac{1}{\Theta p}\lesssim& r^{1+\frac{n}{\Theta p}-\frac{n}{p}}\biggl(\sum_{z\in\Z^{n}}\|\tau_h\nabla v\|_{L^{p}(B_{z,r})}^{\Theta p}\biggr)^\frac{1}{\Theta p}\\
\lesssim&r^{1+\frac{n}{\Theta p}-\frac{n}{p}}\|\tau_h\nabla v\|_{L^p(\R^{n})}\leq 2r^{1+\frac{n}{\Theta p}-\frac{n}{p}}\|\nabla v\|_{L^p(\R^{n})}.
\end{align*}
For the second term on the right-hand side of \eqref{interpolatethu}, we use Jensen inequality, $\|\cdot\|_{\ell^\Theta}\leq \|\cdot\|_{\ell^1}$ and \eqref{ass:coveringbzr} 
\begin{align*}
\biggl(\sum_{z\in\Z^{n}}|B_{z,r}||(\tau_h v)_{B_{z,r}}|^{\Theta p}\biggr)^\frac{1}{\Theta p}\leq& \biggl(\sum_{z\in\Z^{n}}|B_{z,r}|^{1-\Theta}\|\tau_h v\|_{L^p(B_{z,r})}^{\Theta p}\biggr)^\frac{1}{\Theta p}\leq c(n)r^{\frac{n}{\Theta p}-\frac{n}{p}}\|\tau_h v\|_{L^p(\R^{n})}
\end{align*}
together with $\|\tau_h v\|_{L^p(\R^{n})}\lesssim|h|\|\nabla v\|_{L^p(\R^{n})}$. Altogether, we have shown for all $r>0$
\begin{align*}
\|\tau_h u\|_{L^{\Theta p}((B_1(0))_h)}\lesssim r^{\frac{n}{\Theta p}-\frac{n}{p}}(r+|h|) \|\nabla v\|_{L^p(\R^{n})}.
\end{align*}
Clearly, the choice $r=|h|$ together with the definition of $\Theta$ in the form $1+\frac{n}{\Theta p}-\frac{n}{p}=s$ and \eqref{est:vextension} imply $\|\tau_h u\|_{L^{\Theta p}((B_1(0))_h)}\lesssim |h|^s\|\nabla u\|_{L^p(B_1(0)}$ and the claim follows.
%

\step 2  The case $\Theta=\infty$, that is $p>n$ and $s\leq 1-\frac{n}{p}$. In this case the claim follows directly from Morrey's embedding. 
\end{proof}

%
%

Next, we state a Poincar\'e-type inequality involving the semi-norm defined in \eqref{eq:defbesov}. The statement is certainly well-known but for completeness we include a short proof.
\begin{lemma}\label{L:BesovPoincare} Let $Q\subset \R^n$ be a cube. For $p\in(1,\infty)$ and $s\in (0,1)$ there exists $c=c(n,p,s)>0$ such that for all $u\in L^p(Q)$ with $[u]_{s,p,Q}<\infty$, it holds that
\begin{align}\label{eq:BesovPoincareInequality}
\| u - (u)_Q\|_{L^p(Q)}\leq c\lvert Q\rvert^{\frac s n} [u]_{s,p,Q}.
\end{align}
\end{lemma}
\begin{proof}
By standard translation, we can assume that $Q$ is centered at the origin. A direct computation based on Jensen inequality and Fubini yield
\begin{align*}
\| u - (u)_Q\|_{L^p(Q)}^{p}=&\int_{Q}\biggl|\fint_{Q}u(x)-u(y)\,\mathrm dy\biggr|^p\,\mathrm dx\leq \frac{1}{|Q|}\int_{\R^{n}}\int_{\R^{n}}|u(x)-u(y)|^p\mathds 1_Q(x)\mathds 1_Q(y)\,\mathrm dy\mathrm dx\\
=& \frac{1}{|Q|}\int_{\R^{n}}\int_{\R^{n}}|u(x+h)-u(x)|^p\mathds 1_{Q}(x+h)\mathds 1_{Q}(x)\,\mathrm dh\mathrm dx\\
=&\frac{1}{|Q|}\int_{\R^{n}}\int_{(Q)_h}|u(x+h)-u(x)|^p\,\mathrm dx\,\mathrm dh.
\end{align*}
Using $(Q)_h=\emptyset$ if $|h|>c(n)|Q|^\frac{1}{n}=:r$ and the definition of $[u]_{s,p,Q}$, see \eqref{eq:defbesov}, we observe
\begin{align*}
\| u - (u)_Q\|_{L^p(Q)}^{p}\leq&\frac{1}{|Q|}\int_{B_{r}(0)}|h|^{sp}[u]_{s,p,Q}\,\mathrm dh\leq \tilde c(n)|Q|^\frac{sp}{n}[u]_{s,p,Q}
\end{align*}
and the claim follows.

\end{proof}



Finally, we note a basic estimate concerning periodic functions:
\begin{lemma}\label{lem:easyperiodic}
Let $a\in L^1_{\rm loc}(\R^n)$ be $Y$-periodic. For measurable $A\subset \R^n$ it holds%
\begin{equation}
\int_A |a(\tfrac{x}\e)|\dx\leq |(A)_{\sqrt{n}\e}|\|a\|_{L^1(Y)}
\end{equation}
where $(A)_{\sqrt{n}\e}:=A+B_{\sqrt{n}\e}$.
\end{lemma}

\begin{proof}
Set $T_\e:=\{k\in \mathbb Z^n\colon \e(k+Y)\cap A\neq \emptyset\}$. We have
\begin{align*}
\int_A|a(\tfrac{x}\e)|\dx\leq \sum_{\genfrac{}{}{0pt}{2}{z\in \Z^d}{\e(z+Y)\cap A\neq\emptyset}}\int_{\e(z+Y)}|a(\tfrac{x}\e)|\dx=\e^n\#T_\e \|a\|_{L^1(Y)}
\end{align*}
and the claim follows since for all $k\in T_\e$ it holds $\e(k+Y)\subset A+B_{\sqrt{n}\e}$. 
\end{proof}

\subsection{\texorpdfstring{V}{}-functions}
For $1<p<\infty$, we set for $z,z_1,z_2\in\R^n$
\begin{align}\label{def:Wpz}
V_p(z):=(1+|z|)^\frac{p-2} 2 z\quad\mbox{and}\quad W_p(z_1,z_2) := \begin{cases}
			V_{p}(z_1-z_2) \qquad&\text{ if } p\geq 2\\
			(1+|z_1|+|z_2|)^\frac{p-2}2(z_1-z_2) \qquad&\text{ if } p\leq 2.
			\end{cases}
\end{align}
 We collect a number of estimates regarding $V$-functions (and $W$).
\begin{lemma}\label{L:propV}
Let $p>1$. Then the following properties hold:
\begin{enumerate}[(i)]
\item $z\to |V_{p}(z)|^2$ is monotonic increasing in $|z|$
\item For any $z_1,z_2\in \R^n$, $|W_p(z_1,z_2)|^2 \leq |V_{p}(z_1-z_2)|^2$.
\item Scaling. It holds that for $z\in \R^n$ and $\lambda\in \R$,
\begin{align}\label{eq:VScaling}
\min\{\lambda^{p-2},1\}\lambda^2 |V_{p}(z)|^2\leq |V_{p}(\lambda z)|^2\leq  \max\{\lambda^{p-2},1\}\lambda^2 |V_{p}(z)|^2 
\end{align}
\item Triangle inequality. There is $c=c(p)>0$ such that for any $z_1,z_2,z_3\in \R^n$,
\begin{align}\label{eq:triangleV}
|V_{p}(z_1-z_2)|^2\leq c(p)\left(|V_{p}(z_1)|^2+|V_{p}(z_2)|^2\right),\\
|W_{p}(z_1,z_2)|^2\leq c(p)\left(|W_p(z_1,z_3)|^2+|W_p(z_2,z_3)|^2\right).\label{eq:triangleW}
\end{align}
\item There is $c=c(p)>0$ such that for any $z_1,z_2\in \R^n$,
\begin{align}\label{eq:VEquiv}
c(p)^{-1}|V_{p}(z_1)-V_{p}(z_2)|^2\leq (1+|z_1|+|z_2|)^{p-2}|z_1-z_2|^2\leq c(p)|V_{p}(z_1)-V_{p}(z_2)|^2.
\end{align}
\item Young's inequality. There exists $c=c(p)>0$ such that for any $\tau>0$ and any $z,w\in \R^n$,
\begin{align}
|z| |w| \leq \tau|V_{p}(z)|^2+c\max\{\tau^{-\frac1{p-1}},\tau^{-1}\}|V_{p^\prime}(w)|^2\label{eq:YoungV:0}\\
(1+|z|)^{p-2}|z||w| \leq \tau |V_{p}(z)|^2+c\max\{\tau^{-1},\tau^{-(p-1)}\}|V_{p}(w)|^2\label{eq:YoungV}.
\end{align}
\item Poincar\'e inequality. Let $R>r>0$. There is $c=c(n,p,r/R)>0$ such that for any bounded, convex domain $D$ with $B_r(x_0)\subset D \subset B_R(x_0)$,
\begin{align}\label{eq:PoincareV}
\left\|V_{p}\left(\frac{u-(u)_D}{R}\right)\right\|_{\underline L^2(D)}\leq c \|V_{p}(\nabla u)\|_{\underline L^2(D)}.
\end{align}
\item Sobolev-Poincar\'e inequality. Let $r>0$. There is $c=c(n,p)>0$ and $\theta\in(0,1)$ such that for $B=B(x,r)$,
\begin{align}\label{eq:PoincareSobolev}
\left\|V_p\left(u-(u)_B\right)\right\|_{\underline L^2(B)}\leq c \|V_p(\nabla u)\|_{\underline L^{2\theta}(B)}.
\end{align}
\end{enumerate}
\end{lemma}
\begin{proof}
The properties of $V$ we list are well known. \eqref{eq:YoungV:0} is contained in \cite[Lemma 2.7]{Irving2023}, \eqref{eq:PoincareV} can be found in \cite[Lemma 1]{Bhattacharya} and \eqref{eq:PoincareSobolev} is a consequence of \cite[Lemma 2.2]{Celada}. All other properties can be found in \cite[Section 6]{Schmidt} and \cite[Lemma 2.3]{Acerbi2002}. 

Regarding the properties of $W$, it is immediate that $|W_p(z_1,z_2)|^2\leq |V_{p}(z_1-z_2)|^2$ for any $z_1,z_2\in \R^n$. Moreover, for $p\leq 2$, using \eqref{eq:triangleV} (which in fact holds for $V_{\mu,p}(z)=(\mu+|z|)^\frac{p-2} 2 z$ and any $\mu\geq 0$ \cite[Lemma 2.3]{Acerbi2002}),
\begin{align*}
&|W_p(z_1,z_2)|^2\leq c(p) (1+\max\{|z_1|,|z_2|\}+|z_1-z_2|)^{p-2}|z_1-z_2|^2\\
\leq& c(p) (1+\max(|z_1|,|z_2|)+|z_1-z_3|)^{p-2}|z_1-z_3|^2+c(p)(1+\max\{|z_2|,|z_3|\}+|z_2-z_3|)^{p-2}|z_2-z_3|^2\\
\leq& c(p) (|W_p(z_1,z_3)|^2+|W_p(z_2,z_3)|^2.
\end{align*}
\end{proof}

Further, we need some elementary estimates.

\begin{lemma}
Let $p\geq 1$. There exists $c_0=c_0(p)>0$ such that for every $A,B\in\R$, it holds
\begin{equation}\label{eq:elementarybrasco}
|A-B|^p\leq c_0||A|^{p-1}A-|B|^{p-1}B|.
\end{equation}
Further, there exists $c_1=c_1(p)>0$ such that for all $z_1,z_2\in\R^n$ and all $\tau\in(0,1]$
\begin{equation}\label{est:V1pAB}
|V_{p'}((1+|z_1|+|z_2|)^{p-2}z_1)|^2\leq \mathds 1_{p>2}\tau |V_{p}(z_2)|^2 +c_1(1+\tau^{-\frac{p-2}p}) |V_{p}(z_1)|^2
\end{equation}
Finally, if $p\in (1,2]$, there is $c_2=c_2(p)>0$ such that for all $z_1,z_2\in\R^n$ and all $\tau>0$ it holds
\begin{equation}\label{est:WpVsub}
|V_{p}(z_1-z_2)|\leq c|W_{p}(z_1,z_2)|(1+\tau^{-\frac{2-p}p})+\tau (|V_{p}(z_1)|+|V_{p}(z_2)|).
\end{equation}

\end{lemma}

\begin{proof}
\eqref{eq:elementarybrasco} is \cite[Lemma A.3]{Brasco18}. Regarding \eqref{est:V1pAB}, the case $p=2$ is trivial. In the case $p\in(1,2)$, we have 
\begin{align*}
|V_{p'}((1+|z_1|+|z_2|)^{p-2}z_1)|^2=& (1+(1+|z_1|+|z_2|)^{p-2}|z_1|)^\frac{2-p}{p-1}(1+|z_1|+|z_2|)^{2(p-2)}|z_1|^2\\
\leq&(1+|z_1|)^{2-p}(1+|z_1|)^{2(p-2)}|z_1|^2=|V_{p}(z_1)|^2.
\end{align*}
It remains, to consider the case $p>2$. We have,
\begin{align*}
|V_{p'}((1+|z_1|+|z_2|)^{p-2}z_1)|^2=& (1+(1+|z_1|+|z_2|)^{p-2}|z_1|)^\frac{2-p}{p-1}(1+|z_1|+|z_2|)^{2(p-2)}|z_1|^2\\
\leq&2^{2(p-2)}\left(|z_1|^2+\left((|z_1|+|z_2|)^{p-2}|z_1|\right)^\frac{2-p}{p-1}(|z_1|+|z_2|)^{2(p-2)}|z_1|^2\right)\\
=&2^{2(p-2)}\left(|z_1|^2+(|z_1|+|z_2|)^{(p-2)\frac{p}{p-1}}|z_1|^\frac{p}{p-1}\right)\\
\leq&2^{2(p-2)}(|z_1|^2+2^{(p-2)\frac{p}{p-1}}(|z_1|^p+|z_2|^{(p-2)\frac{p}{p-1}}|z_1|^\frac{p}{p-1}))
\end{align*}
and the claimed inequality follows by Young's inequality and $|z|^2+|z|^p\leq 2 |V_{p}(z)|^2$ for $p\geq 2$. 

We now turn to \eqref{est:WpVsub}. The case $p=2$ is trivial and we suppose $p\in(1,2)$. Suppose that $|z_1|+|z_2|\leq1$. Then it holds
$$
|V_{p}(z_1-z_2)|^2\leq |z_1-z_2|^2\leq 2^{2-p}(1+|z_1|+|z_2|)^{p-2}|z_1-z_2|^2=2^{2-p}|W_{p}(z_1,z_2)|^2.
$$
Suppose now $|z_1|+|z_2|\geq1$. Then, we have
\begin{align*}
|V_{p}(z_1-z_2)|^2\leq |z_1-z_2|^p\leq& \tau(1-p/2)(1+|z_1|+|z_2|)^p+\tau^{-\frac{2-p}p}\tfrac{p}2 |W_{p}(z_1,z_2)|^2\\
\leq& 2^p\tau (|z_1|+|z_2|)^p+\tau^{-\frac{2-p}p}\tfrac p 2 |V_{p}(z_1-z_2)|^2.
\end{align*}
Combining the above two estimates with $|z|^p\leq 2^{2-p}|V_{p}(z)|^2$ for $|z|\geq1$, we obtain \eqref{est:WpVsub} (by redefining $\tau$).
\end{proof}

\subsection{Energy estimates} Here, we gather several standard regularity and energy estimates. The estimates collected in the following propositions are essentially standard, but we provide proofs in the appendix \ref{sec:appendixEnergyEstimates} for the sake of completeness. Proposition \ref{prop:regularity} provides estimates in the setting of natural growth, while Proposition \ref{prop:regularity:basiclip} considers the setting of controlled growth.
\begin{proposition}\label{prop:regularity}
Let $1<p<\infty$ and suppose that $\bfa,\overline \bfa\colon \R^n\times\R^n \to \R^n$ satisfy Assumption~\ref{ass:1} with joint constants $\Lambda\in[1,\infty)$ and $\mu\in[0,1]$. Let $B=B_r(x_0)\subset \R^n$ and suppose that $u,w\in W^{1,p}(B)$ and $F\in L^{p'}(B)$ satisfy 
\begin{align}\label{eq:basicRegularityEquation}
\nabla \cdot \bfa(x,\nabla u) = \nabla \cdot F \qquad \mbox{ in $B$}
\end{align}
and
\begin{align}\label{eq:basicRegularityEquation2}
w\in u+W_0^{1,p}(B)\quad\mbox{and}\quad \nabla \cdot \overline \bfa(x,\nabla w) = 0 \qquad \mbox{ in $B$.}
\end{align}
Then the following estimates hold:
\begin{enumerate}[(i)]
\item Energy inequality: There exists $c=c(\Lambda,p)\in[1,\infty)$ such that 
\begin{equation}\label{L:energyestimate:eq1}
\|\mu+|\nabla w|\|_{\underline L^p(B)}\leq c\|\mu+|\nabla u|\|_{\underline L^p(B)}.
\end{equation}
\item Comparison estimate: Suppose that there exists $\delta\geq0$ such that
\begin{equation}\label{prop:regularity:ass:comparison}
\forall z\in\R^n:\quad\sup_{x\in B_1}|\bfa(x,z)-\overline \bfa(x,z)|\leq \delta^{p-1}(\mu+|z|)^{p-2}|z|.
\end{equation}
Then, there exists $c=c(\Lambda,p)\in [1,\infty)$ such that
for all $\tau\in(0,1]$,
\begin{equation}\label{L:energyestimate:eq2}
\|\nabla u-\nabla w\|_{\underline L^p(B)}\leq (\tau+c\delta^{\min\{1,p-1\}}) \|\mu+|\nabla u|\|_{\underline L^p(B)}+c\|\max\{1,\tau^{p-2}\}F\|_{\underline L^{p'}(B)}^\frac1{p-1}.
\end{equation}
\item Caccioppoli inequality: There exists $c=c(p,\Lambda)>0$ such that for any $b\in \R$,
\begin{align}\label{eq:caccioppoli}
 \|\mu+|\nabla u|\|_{\underline L^p(\frac 1 2 B)}\leq c (\mu+r^{-1}\|u-b\|_{\underline L^p(B)}+\|F\|_{\underline L^{p^\prime}(B)}^\frac 1 {p-1}).
\end{align}
\item Higher differentiability: If $\overline \bfa$ is autonomous, there is $c=c(\Lambda,n,p)\in[1,\infty)$ such that
\begin{align}\label{est:reghom:bes}
\forall \rho\in [\tfrac12,1)\qquad[\nabla w]_{\min\{1,\frac 1 {p-1}\},p,\rho B}\leq c \left(r(1-\rho)\right)^{-\min\{1,\frac{1}{p-1}\}}\|\mu+|\nabla u|\|_{L^p(B)}.
\end{align}
\item Meyer's estimate: Suppose that $F\equiv0$. There exists $m=m(\Lambda,n,p)>1$ and $c=c(\Lambda,n,p)>0$ such that
\begin{align}
\forall \rho\in [\tfrac12,1)\qquad\|\nabla u\|_{\underline L^{mp}(\rho B)}\leq c(1-\rho)^{\frac n p\left(\frac 1 m-1\right)}\|\mu+|\nabla u\|_{\underline L^p(B)}\label{est:nonlinearmeyerloc:lp}.
\end{align}
Moreover, if $u\in W^{1,mp}(B)$, then there exists $c=c(\Lambda,p)>0$ such that
\begin{align}
\|\nabla w\|_{L^{mp}(B)}\leq c\|\mu+|\nabla u\|_{L^{mp}(B)}\label{est:nonlinearmeyerglob:lp}.
\end{align}
\end{enumerate}
\end{proposition}

\begin{proposition}\label{prop:regularity:basiclip}
Let $1<p<\infty$ and suppose that $\bfa,\overline \bfa\colon \R^n\times\R^n \to \R^n$ satisfy Assumption~\ref{ass:homonondeg} with joint constant $\Lambda\in[1,\infty)$. There exists $c=c(\Lambda,p)>0$ such that the following is true:
Let $B=B_r(x_0)\subset \R^n$ and suppose that $u,w\in W^{1,p}(B)$ and $F\in L^{p'}(B)$ satisfy 
\begin{align}\label{eq:basicRegularityEquation:basic:lip}
\nabla \cdot \bfa(x,\nabla u) = \nabla \cdot F \qquad \mbox{ in $B$}
\end{align}
and
\begin{align}\label{eq:basicRegularityEquation2:basic:lip}
w\in u+W_0^{1,p}(B)\quad\mbox{and}\quad \nabla \cdot \overline \bfa(x,\nabla w) = 0 \qquad \mbox{ in $B$.}
\end{align}
Then the following estimates hold:
\begin{enumerate}[(i)]
\item Energy inequality: There exists $c=c(\Lambda,p)\in[1,\infty)$ such that 
\begin{equation}\label{L:energyestimate:eq1:basic:lip}
\|V_{p}(\nabla w)\|_{\underline L^2(B)}\leq c\|V_{p}(\nabla u)\|_{\underline L^2(B)}.
\end{equation}
\item Comparison estimate: Suppose that $\bfa =\overline \bfa$. Then, there exists $c=c(\Lambda,p)\in [1,\infty)$ such that
for all $\tau\in(0,1]$,
\begin{align}\label{L:energyestimate:eq2:basic:lip}
\|V_{p}(\nabla u-\nabla w)\|_{\underline L^2(B)}\leq& \tau\|V_{p}(\nabla u)\|_{\underline L^2(B)}+c \left(1+\tau^\frac{p-2}{p-1}\right)\min_{a\in \R^n}\|V_{p'}(F-a)\|_{\underline L^{2}(B)}.
\end{align}
\item Caccioppoli inequality: There exists $c=c(\Lambda,n,p)$ such that for any $b\in \R$,
\begin{align}\label{est:cacc}
\|V_{p}(\nabla u)\|_{\underline L^2(\frac12B)}\leq c \left\|V_{p}\left(\frac{u-b}{r}\right)\right\|_{\underline L^2(B)}
\end{align}
\item Meyer's estimate: Suppose that $F\equiv0$. There exists $m=m(\Lambda,n,p)>1$ and $c=c(\Lambda,n,p)>0$ such that
\begin{align}
\forall \rho\in [\tfrac12,1)\qquad\|V_{p}(\nabla u)\|_{\underline L^{2m}(\rho B)}\leq c(1-\rho)^{\frac n 2\left(\frac 1 m-1\right)}\|V_{p}(\nabla u)\|_{\underline L^2(B)}\label{est:nonlinearmeyerloc:lip}.
\end{align}
Moreover, if $u\in W^{1,2m}(B)$, then there exists $c=c(\Lambda,p)>0$ such that
\begin{align}
\|V_{p}(\nabla w)\|_{L^{2m}(B)}\leq c\|V_{p}(\nabla u)\|_{L^{2m}(B)}\label{est:nonlinearmeyerglob:lip}.
\end{align}
\end{enumerate}
\end{proposition}

\section{Gradient integrability for autonomous equations and proof of Theorem~\ref{thm:Lqx}}\label{sec:thm:Lqx}

In this section, we first consider autonomous operators $\bfa:\R^n\to\R^n$ satisfying Assumption~\ref{ass:1}. We show that weak solutions of the equation
\begin{align}\label{equation}
\nabla \cdot\bfa(\nabla u)= 0\qquad\mbox{in $\Omega\subset\R^n$}
\end{align}
are locally $W^{1,q}$-regular for any $q>1$, see Theorem~\ref{thm:Lq} below. With help of an abstract perturbation, see Lemma~\ref{lem:CZ} below, we extend this result to inhomogeneous equations with non-autonomous operators $\bfa(x,z)$ with a small contrast in the $x$-variable, see Corollary~\ref{cor:Lqx}. From this Theorem~\ref{thm:Lqx} follows.

\begin{theorem}\label{thm:Lq}
Suppose that Assumption~\ref{ass:1} is satisfied for some $1<p<\infty$ and assume that $\bfa$ is autonomous. Let $u\in W^{1,p}(\Omega)$ be a weak solution of \eqref{equation}. Then, $u\in W^{1,q}_{\rm loc}(\Omega)$ for every $q<\infty$. Moreover, there exists for every $q<\infty$ a constant $c=c(n,\Lambda,p,q)\in[1,\infty)$ such that for all balls $B=B_R(x_0)\Subset\Omega$ it holds
\begin{equation}\label{P1:automomous:claim}
\|\mu+|\nabla u|\|_{\underline L^q(\frac12 B)}\leq c\|\mu+|\nabla u|\|_{\underline L^p(B)}.
\end{equation}
Finally, in the case $p\in(1,2]$ there exists $c=c(\Lambda,n,p)\in[1,\infty)$ such that
\begin{equation}\label{P1:automomous:claim:infty}
\|\mu+|\nabla u|\|_{L^\infty(\frac12 B)}\leq c\|\mu+|\nabla u|\|_{\underline L^p(B)}.
\end{equation}
\end{theorem}

\begin{remark}
The novelty of Theorem~\ref{thm:Lq} is in the superquadratic case $p>2$. In the case $p\in(1,2]$, Theorem~\ref{thm:Lq} follows by regularity theory for operators with natural growth, see e.g.\ \cite[Theorem 2.7]{Fonseca2002} in a variational setting. Since, we did not find a precise reference for \eqref{P1:automomous:claim:infty} in our setting we include a proof also for this case.   
\end{remark}

\begin{proof} 
By standard scaling and translation arguments it suffices to consider the case $B=B_1(0)\Subset\Omega$ and establish \eqref{P1:automomous:claim}, \eqref{P1:automomous:claim:infty}. Throughout the proof we write '$\lesssim$' if '$\leq$' holds up to a positive multiplicative constant that depends only $\Lambda,n$ and $p$.

The proof is divided in three steps. We begin by establishing in Step~1 a Caccioppoli inequality of the form
\begin{align*}
&\|\nabla (|\tau_h u|^{\frac{\sigma-1}{p}}\tau_h u)\|_{L^p}\leq C(\Lambda,n,p,\sigma)\|\mu+|\nabla u|\|_{L^{\sigma+p-1}}^{\frac{p+\sigma-1}p}|h|^{\frac{\sigma-1}p}\begin{cases}(\frac{|h|}{\rho-r})^{\frac{p'}{p}}&\mbox{if $p\in[2,\infty)$}\\(\frac{|h|}{\rho-r})&\mbox{if $p\in(1,2)$}\end{cases}.
\end{align*}
see \eqref{P1:pf:step1:claim0} below for the precise form. In the case $p\in(1,2]$, the Caccioppoli inequality implies bounds on $\nabla (|\frac{\tau_h u}{|h|}|^{\frac{\sigma-1}{p}}\frac{\tau_h u}{|h|})$ in $L^p$ that are uniform in $|h|$. In that case we can differentiate powers of $|\nabla u|$ and \eqref{P1:automomous:claim:infty} follows by a 'standard' Moser type iteration, see Step~2 below. In the superquadratic case the Caccioppoli inequality yield bounds on $\nabla (|\frac{\tau_h u}{|h|^\alpha}|^{\frac{\sigma-1}{p}}\frac{\tau_h u}{|h|^\alpha})$ with $\alpha=\frac{\sigma-1+p'}{\sigma-1+p}=1-\frac{p-2}{\sigma-1+p}\in(0,1)$ that are uniform in $|h|$. In this case we appeal to a fractional Moser type iteration heavily inspired by \cite{Brasco,deFilippis2024a} to obtain \eqref{P1:automomous:claim}, see Step~3 below.

\step 1 A Caccioppoli inequality. We claim that for all $\frac12\leq r<\rho\leq1$ and $h\in\R^n$ with $|h|\leq h_0:=(\rho-r)/32$ the following is true: For every $\sigma\geq1$ there exists $c=c(\Lambda,n,p)\in[1,\infty)$ such that 
\begin{align}\label{P1:pf:step1:claim0}
&\|\nabla (|\tau_h u|^{\frac{\sigma-1}{p}}\tau_h u)\|_{L^p(B_{r+\frac{\rho-r}4})}
\leq& c(\sigma-1+p)\|\mu+|\nabla u|\|_{L^{\sigma+p-1}(B_\rho)}^\frac{\sigma-1+p} p |h|^\frac{\sigma-1} p\begin{cases}(\frac{|h|}{\rho-r})^{\frac{p'}{p}}&\mbox{if $p\in[2,\infty)$}\\
(\frac{|h|}{\rho-r}) &\mbox{if $p\in(1,2)$}\end{cases}.
\end{align}
Without loss of generality, we can assume that $u\in W^{1,\sigma-1+p}(B_\rho)$ (otherwise  \eqref{P1:pf:step1:claim0} is an empty statement). Let $\eta\in C_c^\infty (B_\rho)$ be a smooth cut-off function satisfying
\begin{equation}\label{ass:eta}
0\leq\eta\leq1,\quad \supp \eta\subset B_{\rho-(\rho-r)/4},\quad \eta\equiv 1\mbox{ on $B_{r+(\rho-r)/4}$},\quad |\nabla^k \eta|\lesssim (\rho-r)^k.
\end{equation}
The choice of $\eta$ and the assumption $u\in W^{1,\sigma-1+p}(B_\rho)$ ensures that $\psi:=\tau_{-h}(\eta^{p_2} \tau_h u\lvert \tau_h u\rvert^{\sigma-1})$ where ${p_2:=\max\{p,2\}}$ is a valid test function for \eqref{equation}.  Indeed, using \eqref{ass:1:growth} and $\nabla u\in L^{\sigma-1+p}(B_\rho)$, we have $\bfa(\nabla u)\in L^{\frac{\sigma-1+p}{p-1}}(B_\rho)$. An application of H\"older's inequality shows that $\psi\in W^{1,\frac{\sigma-1+p}{\sigma}}(B_\rho)$ with compact support. A standard approximation argument (using $\frac{p-1}{\sigma-1+p}+\frac{\sigma}{\sigma-1+p}=1$) then gives that $\psi$ is a valid test function in \eqref{equation},  allowing us to deduce
\begin{align}
\sigma I:=&\sigma \int_{\R^n} \eta^{p_2}|\tau_hu|^{\sigma-1} \langle \tau_h \bfa(\nabla u),\nabla \tau_h u\rangle\dx\notag\\
 =&-\int_{\R^n} \langle \bfa(\nabla u),\tau_{-h}(|\tau_hu|^{\sigma-1}(\tau_h u) \nabla (\eta^{p_2}))\rangle\dx=:II. \label{P1:basic:IandII}
\end{align}
The growth conditions of $\bfa$, \eqref{ass:eta} and H\"older inequality imply
\begin{align}\label{pf:Thm4:s1:est:II:1}
|II|\leq& \Lambda \|(\mu+|\nabla u|)\|_{L^{\sigma-1+p}(B_\rho)}^{p-1}\|\tau_{-h}(|\tau_h u|^{\sigma-1}(\tau_h u)\nabla (\eta^{p_2}))\|_{L^\frac{\sigma-1+p}{\sigma}(\R^n)}.
\end{align}
Using \eqref{eq:fact4:1} (with $\Omega=\R^n$), the product rule and \eqref{ass:eta}, we find
\begin{align}\label{pf:Thm4:s1:est:II:2}
&\|\tau_h(|\tau_h u|^{\sigma-1}(\tau_h u)\nabla (\eta^{p_2}))\|_{L^\frac{\sigma-1+p}{\sigma}(\R^n)}\notag\\
\leq& |h|\|\nabla (|\tau_h u|^{\sigma-1}(\tau_h u)\nabla (\eta^{p_2}))\|_{L^{\frac{\sigma-1+p}\sigma}(\R^n)}\notag\\
 \lesssim&\frac{|h|}{\rho-r}\|\eta^{p_2-1}\nabla (|\tau_h u|^{\sigma-1}(\tau_h u))\|_{L^{\frac{\sigma-1+p}\sigma}(\R^n)}+\frac{|h|^{1+\sigma}}{(\rho-r)^2}\|\nabla u\|_{L^{\sigma-1+p}(B_\rho)}^\sigma,
\end{align}
where we use for the last inequality \eqref{eq:fact4:1}, \eqref{ass:eta} and $|h|\leq (\rho-r)/32$ in the form
$$\||\nabla^2 \eta| |\tau_h u|^\sigma\|_{L^\frac{\sigma-1+p}\sigma(\R^n)}\lesssim (\rho-r)^{-2}|h|^\sigma\|\nabla u\|_{L^{\sigma-1+p}(B_\rho)}^\sigma.$$ 
By chain rule and H\"older inequality with exponents $\frac{\sigma p}{\sigma-1+p}$ and $\frac{\sigma p}{(\sigma-1)(p-1)}$ , we find
\begin{align}\label{pf:Thm4:s1:est:II:3}
\|\eta^{p_2-1}\nabla (|\tau_h u|^{\sigma-1}(\tau_h u))\|_{L^{\frac{\sigma-1+p}\sigma}(\R^n)}\leq&\sigma \|\eta^{p_2-1}  |\tau_h u|^{\sigma-1}| \nabla \tau_h u|\|_{L^{\frac{\sigma-1+p}\sigma}(\R^n)}\notag\\
\leq&\sigma \|\eta  |\tau_h u|^{\frac{\sigma-1}{p}}| \nabla \tau_h u|\|_{L^{p}(\R^n)}\|\eta^{p_2-2} |\tau_h u|^\frac{\sigma-1}{p'}\|_{L^{\frac{p(\sigma-1+p)}{(\sigma-1)(p-1)}}(\R^n)}\notag\\
\leq&\sigma \|\eta  |\tau_h u|^{\frac{\sigma-1}{p}}| \nabla \tau_h u|\|_{L^{p}(\R^n)}(h\|\nabla u\|_{L^{\sigma-1+p}(B_\rho)})^{\frac{\sigma-1}{p'}}
\end{align}
(with the understanding $\infty=\frac10$). Combining, \eqref{pf:Thm4:s1:est:II:1}--\eqref{pf:Thm4:s1:est:II:3} we have
\begin{align}\label{pf:Thm4:s1:est:II:4}
|II|\lesssim&
\frac{\sigma|h|^{\frac{p'+\sigma-1}{p'}}\|\mu+|\nabla u|\|_{L^{\sigma-1+p}(B_\rho)}^{\frac{\sigma-1+p}{p'}}}{\rho-r}\|\eta  |\tau_h u|^{\frac{\sigma-1}{p}}| \nabla \tau_h u|\|_{L^{p}(\R^n)}+\frac{|h|^{1+\sigma}\|\mu+|\nabla u|\|_{L^{\sigma-1+p}(B_\rho)}^{\sigma+p-1}}{(\rho-r)^2}.
\end{align}
The coercivity condition \eqref{ass:1:monotone} and chain rule imply in the case $p\geq2$
\begin{align*}
I \geq& \frac1\Lambda \int_{\R^n} \eta^{p_2} \lvert \tau_h \nabla u\rvert^p \lvert \tau_h u\rvert^{\sigma-1}\,\dx\geq \frac1\Lambda\biggl(\frac{p}{\sigma-1+p}\biggr)^p \int_{\R^n} \eta^{p_2}\left\lvert \nabla (\left\lvert\tau_h u \right\rvert^\frac{\sigma-1} p \tau_h u)\right\rvert^p\dx.
\end{align*}
Hence, \eqref{pf:Thm4:s1:est:II:4} and a suitable application of Young's inequality yield
\begin{eqnarray*}
\int_B \eta^p  \left\lvert \nabla (\left\lvert\tau_h u \right\rvert^\frac{\sigma-1} p \tau_h u)\right\rvert^p\dx&\stackrel{\eqref{P1:basic:IandII}}\lesssim& \biggl(\frac{\sigma-1+p}{p}\biggr)^p\frac1{\sigma}|II|\\
&\lesssim&(\sigma-1+p)^p\biggl(\frac{|h|^{\sigma-1+p'}}{(\rho-r)^{p'}}+\frac{|h|^{1+\sigma}}{(\rho-r)^2}\biggr)\|\mu+|\nabla u|\|_{L^{\sigma-1+p}(B_\rho)}^{\sigma-1+p}.
\end{eqnarray*}
The claimed estimate \eqref{P1:pf:step1:claim0} (in the case $p\geq2$) follows using  $|h|\leq \rho-r$.

In the case $p\in(1,2]$, the coercivity condition \eqref{ass:1:monotone} implies 
\begin{align}\label{P1:pf:step1:monotoneleq}
I \geq& \frac1\Lambda \int_{\R^n} \eta^2 (\mu+|\nabla u|+|\nabla u(\cdot+h)|)^{p-2}\lvert \tau_h \nabla u\rvert^2 \lvert \tau_h u\rvert^{\sigma-1}\dx=:\frac1{\Lambda} I'.
\end{align}
Moreover, we have with help of H\"older inequality
\begin{align}\label{pf:Thm4:s1:est:II:5}
\|\eta  |\tau_h u|^{\frac{\sigma-1}{p}}| \nabla \tau_h u|\|_{L^{p}(\R^n)}\leq& \biggl(\int_{\R^n}\eta^2|\tau_hu|^{\sigma-1}|\nabla \tau_h u|^2(\mu+|\nabla u|+|\nabla u(\cdot+h)|)^{p-2}\dx\biggr)^\frac12\notag\\
&\times \biggl(\int_{B_{\rho-(\rho-r)/4}}|\tau_h u|^{\sigma-1}((\mu+|\nabla u|+|\nabla u(\cdot+h)|)^{p}\dx\biggr)^\frac{2-p}{2p}\notag\\
\lesssim&  (|h|^{\sigma-1}\|\mu+|\nabla u|\|_{L^{\sigma+p-1}(B_\rho)}^{\sigma-1+p})^{\frac{2-p}{2p}}\sqrt{I'}
\end{align}
Combining \eqref{P1:basic:IandII}, \eqref{P1:pf:step1:monotoneleq} and \eqref{pf:Thm4:s1:est:II:5} with a suitable application of Young's inequality, we find
\begin{align*}
\int_{\R^n}\eta^2|\tau_hu|^{\sigma-1}|\nabla \tau_h u|^2(\mu+|\nabla u|+|\nabla u(\cdot+h)|)^{p-2}\dx\lesssim&\frac{(\sigma-1+p)|h|^{1+\sigma}\|\mu+|\nabla u|\|_{L^{\sigma-1+p}(B_\rho)}^{\sigma+p-1}}{(\rho-r)^2}.
\end{align*}
Inserting the previous estimate for $I'$ into \eqref{pf:Thm4:s1:est:II:5}, we obtain the claimed inequality \eqref{P1:pf:step1:claim0} for $p\in(1,2)$.

\step 2 The case $p\in(1,2]$. Consider the situation of Step~1. Let $\Theta=\Theta(n,p)>1$ be given by
\begin{equation}\label{P1:pf:choiceThetasub}
\Theta =\frac{n}{n-p}\quad\mbox{if $p<n$ and}\quad  \Theta =2\quad \mbox{if $p=n=2$}.
\end{equation}
Sobolev embedding in combination with \eqref{P1:pf:step1:claim0} and $\frac{1}{2}\leq r<\rho\leq1$  implies
\begin{align*}
\||\tau_h u|^\frac{\sigma-1}{p}\tau_h u\|_{L^{\Theta p}(B_{r+(\rho-r)/4})}\lesssim& \||\nabla (|\tau_h u|^\frac{\sigma-1}{p}\tau_h u)\|_{L^{p}(B_{r+(\rho-r)/4})}+ \|| \tau_h u|^\frac{\sigma-1}{p}\tau_h u\|_{L^{p}(B_{r+(\rho-r)/4})}\\
\lesssim&(\sigma-1+p)\frac{|h|^\frac{\sigma+p-1}{p}}{\rho-r}\|\mu+|\nabla u|\|_{L^{\sigma-1+p}(B_\rho)}^\frac{\sigma-1+p}{p}.
\end{align*}
Hence, we have for all $h\in\R^n$ with $0<|h|\leq h_0=(\rho-r)/32$ that
\begin{align*}
\int_{B_{r+(\rho-r)/4}}\biggl|\frac{\tau_h u}{|h|}\biggr|^{\Theta (\sigma+p-1)}\dx\lesssim  \frac{(\sigma-1+p)^{\Theta p}}{(\rho-r)^{\Theta p}}\|\mu+|\nabla u|\|_{L^{\sigma-1+p}(B_\rho)}^{\Theta(\sigma-1+p)}.
\end{align*}
and thus by \eqref{eq:gradientToDifferenceQuotient} and the arbitrary direction of $h$,
\begin{align*}
\int_{B_r}|\nabla u|^{\Theta (\sigma+p-1)}\dx\lesssim \biggl(\frac{(\sigma-1+p)}{\rho-r}\biggr)^{\Theta p}\|\mu+|\nabla u|\|_{L^{\sigma-1+p}(B_\rho)}^{\Theta(\sigma-1+p)}.
\end{align*}

Setting $q=q(p,\sigma):=\sigma+p-1$, the above estimate implies that there exists $c_1 = c_1(\Lambda,n,p)\in[1,\infty)$ such that
\begin{align}\label{P1:pf:onestepimprovesub}
\|\mu+|\nabla u|\|_{L^{\Theta q }(B_r)}\leq \biggl( \frac{c_1}{(\rho-r)}\biggr)^{\frac{p}{q}}\|\mu+|\nabla u|\|_{L^{q}(B_\rho)}
\end{align}
Estimate \eqref{P1:pf:onestepimprovesub} can be iterated: Set  $\rho_i= \frac 1 2 \left(1+ 2^{-i}\right)$, $\sigma_0 = 1$ and 
\begin{align*}
(\sigma_i+p-1)\Theta = \sigma_{i+1}+p-1 ,
\end{align*}
where $\Theta=\Theta(n,p)>1$ is given in \eqref{P1:pf:choiceThetasub}. The choice of $\sigma_i$ ensures 
$$
q_i:=\sigma_i+p-1=\Theta^ip\quad\mbox{for all $i\in\mathbb N$}
$$
and repeated applications of \eqref{P1:pf:onestepimprovesub} imply for every $k\in\mathbb N$
\begin{equation*}
\|\mu+|\nabla u|\|_{L^{q_k}(B_{\rho_k})}\leq \prod_{i=0}^k \biggl( c_1 2^i\Theta^ip\biggr)^{\Theta^{-i}}\|\mu+|\nabla u|\|_{L^{p}(B_1)}
\end{equation*}
and thus
\begin{equation*}
\|\mu+|\nabla u|\|_{L^{\infty}(B_{\frac12})}\leq  \biggl( 2c_1\Theta p\biggr)^{\sum_{i=0}^\infty i \Theta^{-i}}\|\mu+|\nabla u|\|_{L^{p}(B_1)}.
\end{equation*}
Clearly, this proves the claim in the case $p\in(1,2]$.

\step 3 The case $p\in(2,\infty)$. Consider the setting of Step~1. Let $s=s(p)\in(0,1)$ and $\Theta=\Theta(n)>1$ be given by
\begin{equation}\label{P1:pf:choiceTheta}
s=1-\frac{1}{2p}\quad\mbox{and}\quad \Theta:=\frac{2n}{2n-1}.
\end{equation}
The fractional Sobolev embedding \eqref{eq:fact2} (note that $1-(1-s)\frac{p}{n}=\frac{1}{\Theta}$) implies 
\begin{align}\label{est:besovsobolevthm1}
\left[\left\lvert\tau_h u\right\rvert^\frac{\sigma-1} p \tau_h u \right]_{s,\Theta p,B_{r+(\rho-r)/4}}^p\lesssim& \|\nabla (\left\lvert\tau_h u\right\rvert^\frac{\sigma-1} p \tau_h u)\|_{L^p(B_{r+(\rho-r)/4})}^p,
\end{align}
where we use that $\frac{1}{2}\leq r<\rho\leq1$. The above estimate in combination with the pointwise inequality \eqref{eq:elementarybrasco} (with $q=\frac{\sigma-1+p}p$) implies
\begin{eqnarray}\label{P1:pf:estI}
&&\lvert h\rvert^{-sp}\left(\int_{B_r+(\rho-r)/8}\lvert \tau_h^2 u \rvert^{\Theta (\sigma-1+p)}\dx\right)^\frac 1 {\Theta}\notag\\&\stackrel{\eqref{eq:elementarybrasco}}{\leq}& C \lvert h\rvert^{-sp}\left\|\tau_{-h}\biggl(\left\lvert \tau_h u\right\rvert^\frac{\sigma-1} p \tau_h u\;\biggr)\right\|_{L^{\Theta p}(B_{r+(\rho-r)/8})}^p\notag\\
&\stackrel{\eqref{eq:defbesov}}{\leq}&C \left[\left\lvert\tau_h u\right\rvert^\frac{\sigma-1} p \tau_h u \right]_{s,\Theta p,B_{r+(\rho-r)/4}}^p\notag\\
&\stackrel{\eqref{est:besovsobolevthm1},\eqref{P1:pf:step1:claim0}}{\leq}& C \frac{|h|^{\sigma+p'-1}}{(\rho-r)^{p'}}\|\mu+|\nabla u|\|_{L^{\sigma+p-1}(B_\rho)}^{\sigma+p-1}
\end{eqnarray}
where here and below $C=C(\Lambda,n,p,\sigma)\in[1,\infty)$ which might change at each occurence. Clearly, \eqref{P1:pf:estI} implies
\begin{align*}
\left(\int_{B_{r+(\rho-r)/8}} \left\lvert\frac{\tau_{h}^2 u}{\lvert h\rvert^{\frac{sp+\sigma+p'-1}{\sigma-1+p}}}\right\rvert^{\Theta (\sigma-1+p)}\dx\right)^\frac 1 {\Theta}\leq& \frac{C}{(\rho-r)^{p'}} \|\mu+|\nabla u|\|_{L^{\sigma-1+p}(B_\rho)}^{\sigma-1+p}.
\end{align*}
The choice of $\Theta$ and $s$, see \eqref{P1:pf:choiceTheta}, implies
\begin{align*}
sp+\sigma + p'-1 = \sigma + p + p' - \frac 3 2.
\end{align*}
In particular, applying \eqref{eq:fact4} with  $q=\Theta(\sigma-1+p)$, $\e = \Theta (p'-\frac12)\in (0,q)$ and $h_0=(\rho-r)/40$ to $u-(u)_{B_\rho}$, we find, using also Poincar\'{e}'s inequality,
\begin{align}\label{eq:main}
\|\mu+\lvert\nabla u\rvert\|_{L^{\Theta(\sigma+p-1)}(B_r)}\leq \frac{C}{(\rho-r)^\gamma}\|\mu+\lvert\nabla u\rvert\|_{L^{\sigma+p-1}(B_\rho)}.
\end{align}
with $\gamma=1+\frac{2p'-\frac{1}{2}}{\sigma-1+p}$ 
 and $C=C(\Lambda,n,p,\sigma)>0$. The claim \eqref{P1:automomous:claim} follows from \eqref{eq:main} by a straightfoward iteration: As in Step~2, we set $\rho_i= \frac 1 2 \left(1+ 2^{-i}\right)$, $\sigma_0 = 1$ and 
\begin{align*}
(\sigma_i+p-1)\Theta = \sigma_{i+1}+p-1,
\end{align*}
where $\Theta>1$ is now given in \eqref{P1:pf:choiceTheta}. We observe that $q_k:=\sigma_k+p-1=\Theta^k p$ and by repeated application of \eqref{eq:main}, we find for every $k\in\mathbb N$ that 
\begin{equation*}
\|\mu+\lvert\nabla u\rvert\|_{L^{q_k}(B_{\rho_k})}\leq C(k,\Lambda,n,p)\|\mu+\lvert\nabla u\rvert\|_{L^{p}(B_1)}
\end{equation*}
Since $q_k\to\infty$ as $k\to\infty$, we obtain the claimed estimate \eqref{P1:automomous:claim} for $B=B_1(0)$ (using $\rho_k\geq\frac12$ for all $k\in\mathbb N$). 

\end{proof}

Next, we provide an adaption to $L^p$ of \cite[Lemma 7.2]{Armstrong2019}, which is useful in order to derive Calderon-Zygmund-type estimates.
\begin{lemma}\label{lem:CZ}
Let $p\in(1,\infty)$. For each $q\in(2,\infty)$, $q_0\in(q,\infty]$ and $A\geq 1$, there exist constants ${\delta_0(A,n,p,q,q_0)>0}$ and $C(A,n,p,q,q_0)<\infty$ such that for every $\delta\in(0,\delta_0]$ the following holds: Let $f\in L^p(B_1)$ and $g \in L^\frac{q}{p-1}(B_1) $ be such that for every $x\in B_\frac12$ and $r\in (0,1/4]$, there exists $f_{x,r}\in L^{q_0}(B_r(x))$ satisfying
\begin{align}
\|f_{x,r}\|_{\underline L^{q_0}(B_r(x))}\leq A\|f\|_{\underline L^p(B_{2r}(x))}+\|g\|_{\underline L^{p^\prime}(B_{2r(x)})}^\frac 1 {p-1},\label{lem:CZ:ass1}\\
\|f-f_{x,r}\|_{\underline L^p(B_r(x))}\leq \delta \|f\|_{\underline L^p(B_{2r}(x))}+\|g\|_{\underline L^{p^\prime}(B_{2r(x)})}^\frac 1 {p-1}\label{lem:CZ:ass2}
\end{align}
Then $f\in L^q(B_{1/2})$ and we have the estimate
\begin{align}\label{lem:CZ:claim}
\|f\|_{L^q(B_{1/2})}\leq C \left(\|f\|_{L^p(B_1)}+\|g\|_{L^\frac q {p-1}(B_1)}^\frac 1 {p-1}\right).
\end{align}
\end{lemma}
\begin{proof}
We reduce the case $p\in(1,\infty)$ to the case $p=2$ for which the statement is proven in \cite[Lemma 7.2]{Armstrong2019}. 

Set $\tilde f = \lvert f\rvert^\frac{p-2} 2 f$, $\tilde f_{x,r} = \lvert f_{x,r}\rvert^\frac{p-2} 2 f_{x,r}$ and $\tilde g = \lvert  g\rvert^\frac{p^\prime-2} 2 g$. We claim that there exists $c_p\in[1,\infty)$ depending only on $p$ and $d$ such that for any $x\in B_1$ and $r\in (0,1/4]$, we obtain
\begin{align}
\|\tilde f_{x,r}\|_{\underline L^\frac{2q_0}{p}(B_r(x))}\leq& c_pA^\frac{p}2\|\tilde f\|_{\underline L^p(B_{2r}(x)}+\|c_p\tilde g\|_{\underline L^{2}(B_{2r(x)})}\label{lem:CZ:pf:1}\\
\|\tilde f-\tilde f_{x,r}\|_{\underline L^2(B_r(x))}\leq& c_p\delta^{\min\{1,p/2\}} \|\tilde f\|_{\underline L^2(B_{2r}(x))}+\|c_p\delta^{\min\{0,1-\frac{p}2\}}\tilde g\|_{\underline L^{2}(B_{2r(x)})}\label{lem:CZ:pf:2}
\end{align}
Estimates \eqref{lem:CZ:pf:1} and \eqref{lem:CZ:pf:2} in combination with \cite[Lemma 7.2]{Armstrong2019}, imply that there exists a constant $C=C(A,n,p,q,q_0)\in[1,\infty)$ such that
\begin{align*}
\|\tilde f\|_{L^\frac{2q}p(B_{1/2})}\leq C \left(\|\tilde f\|_{L^2(B_1)}+\|\tilde g\|_{L^{\frac 2p q}(B_1)}\right).
\end{align*}
Clearly, this implies \eqref{lem:CZ:claim} by the definition of $\tilde f$ and $\tilde g$. The claimed estimate \eqref{lem:CZ:pf:1} follows by the definition of $\tilde f$, $\tilde f_{x,r}$, $\tilde g$ and \eqref{lem:CZ:ass1} in the form
\begin{eqnarray*}
\|\tilde f_{x,r}\|_{\underline L^{\frac{2 q_0} p}(B_r(x))}^\frac 2 p &=& \|f_{x,r}\|_{\underline L^{q_0}(B_r(x))}\stackrel{\eqref{lem:CZ:ass1}}\leq A\|f\|_{\underline L^p(B_{2r}(x))}+\|g\|_{\underline L^{p^\prime}(B_{2r}(x))}^\frac 1 {p-1}\\
&=& A \|\tilde f\|_{\underline L^2(B_{2r}(x))}^\frac 2 p+ \|\tilde g\|_{\underline L^2(B_{2r}(x))}^\frac 2 p.
\end{eqnarray*}
Raising the above inequality to the $\frac{p}2$-th power, we obtain \eqref{lem:CZ:pf:1}. Next, we show \eqref{lem:CZ:pf:2}. In what follows $c_p\in[1,\infty)$ might change in each occurrence but only depends on $p$ and $n$. For $p\in(2,\infty)$, we estimate using \eqref{lem:CZ:ass2},
\begin{eqnarray*}
\|\tilde f_{x,r}-\tilde f\|_{\underline L^2(B_r(x))}^2&\leq& c_p\fint_{B_r(x)} (\lvert f\rvert^{p-2}+\lvert f_{x,r}\rvert^{p-2})\lvert f-f_{x,r}\rvert^2\dx\\
&\leq& c_p\left(\|f\|_{\underline L^p(B_r(x))}^{p-2}+\|f_{x,r}\|_{\underline L^p(B_r(x)}^{p-2}\right) \|f-f_{x,r}\|_{\underline L^p(B_r(x))}^2\\
&\leq& \delta^2\|f\|_{\underline L^p(B_r(x))}^p+c_p\delta^{2-p} \|f-f_{x,r}\|_{\underline L^p(B_r(x))}^p\\
&\stackrel{\eqref{lem:CZ:ass2}}\leq&\delta^2 \|f\|_{\underline L^p(B_r(x))}^p+c_p\delta^{2-p}(\delta^p\|f\|_{\underline L^p(B_r(x))}^p+\|g\|_{\underline L^{p'}(B_{2r}(x))}^{p'})\\
&\leq& c_p\delta^2\|f\|_{\underline L^p(B_r(x))}^p+c_p\delta^{2-p}\|g\|_{\underline L^{p'}(B_{2r}(x)}^{p'}\\
&=&c_p\delta^2\|\tilde f\|_{\underline L^2(B_{2r}(x))}^2+c_p\delta^{2-p}\|\tilde g\|_{L^2(B_{2r}(x))}^2,
\end{eqnarray*}
which yields \eqref{lem:CZ:pf:2} in the case $p\in(2,\infty)$. In the case $p\in(1,2]$, we have
\begin{align*}
\|\tilde f_{x,r}-\tilde f\|_{\underline L^2(B_r(x))}^2\leq& c_p\|f-f_{x,r}\|_{\underline L^p(B_r(x))}^p
\leq c_p\delta^p \|f\|_{L^p(B_{2r}(x))}^p+c_p \|g\|_{L^{p^\prime}(B_{2r}(x))}^{p^\prime}\\
=& c_p\delta^p \|\tilde f\|_{L^2(B_{2r}(x))}^2+c_p \|\tilde g\|_{L^2(B_{2r}(x))}^2
\end{align*}
and thus \eqref{lem:CZ:pf:2} follows in the remaining case $p\in(1,2]$.
\end{proof}

Appealing to Lemma~\ref{lem:CZ}, we extend Theorem~\ref{thm:Lq} to a Calderon-Zygmund theory for inhomogeneous equations with non-autonomous operators with sufficiently small oscillations in the $x$-variable.

\begin{corollary}\label{cor:Lqx}
Assume $1< p<\infty$ and suppose that $\bfa(x,z)$ satisfies Assumption~\ref{ass:1}. For every $q\in(p,\infty)$ there exists $\delta_0=\delta_0(\Lambda,n,p,q)>0$ and $C=C(\Lambda,n,p,q)>0$ such that the following is true: Suppose that
\begin{equation}\label{thm:lpqx:ass1}
\forall z\in\R^n:\quad\sup_{x\in B_1}|\bfa(x,z)-\bfa(0,z)|\leq\Lambda \delta^{p-1}(\mu+|z|)^{p-2}|z|
\end{equation}
for some $\delta\in(0,\delta_0]$ and $u\in W^{1,p}(B_1)$ is a weak solution to
\begin{equation}\label{thm:lpqx:ass2}
\nabla \cdot\bfa(x,\nabla u)=\nabla\cdot F\quad\mbox{in $B_1$}.
\end{equation}
Then it holds
\begin{align*}
\left(\fint_{B_{1/2}} (\mu+\lvert \nabla u\rvert)^q\dx\right)^\frac 1 q\leq C\left(\left(\fint_{B_1} (\mu+\lvert \nabla u\rvert)^p\dx\right)^\frac 1 p+ \left(\fint_{B_1} \lvert F\rvert^\frac{q}{p-1}\dx\right)^\frac 1 {q}\right).
\end{align*}
\end{corollary}
\begin{proof}
Fix $1<p<q<\infty$. In view of Lemma~\ref{lem:CZ} it suffices to show that there exists ${A=A(\Lambda,n,p,q)>0}$ such that for all $x\in B_\frac12$ and $r\in(0,\frac14]$ there exists $u_{x,r}\in W^{1,p}(B_{2r}(x))$ satisfying for every $1\geq\delta_1\geq\delta$, where  $\delta>0$ is the constant in \eqref{thm:lpqx:ass1},
\begin{align}
\|\mu+\lvert \nabla  u_{x,r}\rvert\|_{\underline L^{2q}(B_r(x))}\leq& A\|\mu+\lvert \nabla u\rvert\|_{\underline L^p(B_{2r}(x))}\label{thm:lpqx:pf:claim1}\\
\|\nabla \overline u-\nabla u\|_{\underline L^p(B_{2r}(x))}\leq& A\delta_1^{\min\{1,p-1\}} \|\mu+|\nabla u|\|_{\underline L^p(B_{2r}(x))}+A\delta_1^{\min\{0,p-2\}}\|F\|_{\underline L^{p'}(B_{2r}(x))}^\frac1{p-1}\label{thm:lpqx:pf:claim2}.
\end{align}
Indeed, for $\delta>0$ sufficiently small, we can choose $\delta_1\geq\delta$ such that $A\delta_1^{\min\{1,p-1\}}=\delta_0$, where ${\delta_0=\delta_0(A,n,p,q,2q)}$ is the corresponding constant in Lemma~\ref{lem:CZ} and the claim follows from Lemma~\ref{lem:CZ} with the choice $g:=(A\delta_1^{\min\{0,p-2\}})^\frac 1 pF$. 

Fix $x\in B_\frac12$ and $r\in(0,\frac14]$ and define $u_{x,r}\in W^{1,p}(B_{2r}(x))$ by
\begin{equation*}
u_{x,r}\in u+W_0^{1,p}(B_{2r}(x))\quad\mbox{and}\quad \nabla \cdot\bfa(0,\nabla u_{x,r})=0\mbox{ in $B_{2r}(x)$}.
\end{equation*} 
We show that $u_{x,r}$ satisfies \eqref{thm:lpqx:pf:claim1} and \eqref{thm:lpqx:pf:claim2}. Appealing to Theorem~\ref{thm:Lq}, we have 
\begin{equation*}
\|\mu+\lvert \nabla  u_{x,r}\rvert\|_{\underline L^{2q}(B_r(x))}\leq C \|\mu+\lvert \nabla u_{x,r}\rvert\|_{\underline L^p(B_{2r}(x))},
\end{equation*}
where $C=C(\Lambda,n,p,q)>0$. The above estimate in combination with the energy estimate \eqref{L:energyestimate:eq1} of Proposition~\ref{prop:regularity},
\begin{equation*}
\|\mu+\lvert \nabla u_{x,r}\rvert\|_{\underline L^p(B_{2r}(x))}\leq c \|\mu+\lvert \nabla u\rvert\|_{\underline L^p(B_{2r}(x))}
\end{equation*}
with $c=c(\Lambda,p)>0$ yield \eqref{thm:lpqx:pf:claim1}. The remaining estimate \eqref{thm:lpqx:pf:claim2} follows directly from \eqref{L:energyestimate:eq2} of Proposition~\ref{prop:regularity} with $\overline\bfa(x,z):=\bfa(0,z)$ and $\tau=\delta_1^{\min\{1,p-1\}}$.

\end{proof}

We end this section by noting
\begin{proof}[Proof of Theorem~\ref{thm:Lqx}]
This follows directly from Corollary~\ref{cor:Lqx} and a suitable scaling and translation argument. Indeed, let $q\in(p,\infty)$ be fixed and let $\delta_0=\delta_0(\Lambda,n,p,q)$ be as in Corollary~\ref{cor:Lqx}. By the continuity of $\omega$ and $\omega(0)=0$, we find $r_q=r_q(\delta_0,\Lambda,\omega,p)\in(0,R/2)$ such that $\omega(2r_q)\leq \Lambda \delta_0^{p-1}$. Hence, for every ball $B_{r_q}(x_0)\subset \Omega$, after translating and rescaling to the unit ball, assumption \eqref{thm:lpqx:ass1} of Corollary~\ref{cor:Lqx} is satisfied. Applying the corollary on each such ball and covering $B_{R/2}(x)$ by finitely many balls of radius $r_q/2$, we obtain \eqref{est:thm:Lqx}. The resulting constant depends on $\Lambda,n,p,q,\omega$ and $R$ as claimed.
\end{proof}

\begin{remark}\label{rem:anistropic1}
The orthotropic $p$-Laplace operator for $p\geq 2$ satisfies Assumption~\ref{ass:1} and hence Theorem \ref{thm:Lqx} applies to this operator. While for this particular example, Lipschitz regularity of solutions has been known since \cite{Uraltseva} and is by now known under a variety of assumptions on a non-divergence form right-hand side \cite{Demengel,Bousquet2018}, Assumption \ref{ass:1} applies to a more general class of degenerate elliptic operators. In fact, the Finsler $p$-Laplacian, $\bfa(z)=\|z\|^{p-1}\nabla(\|z\|)$ of any $\min\{p,2\}$-smooth and $\max\{p,2\}$-convex norm $\|\cdot\|$ satisfies Assumption~\ref{ass:1},  see \cite{Roach1991} both for the concept of $\tau$-smoothness and $\sigma$-convexity of norms and proof of this. 
To the interested reader we also recommend \cite{Goering2023} where the abstract convex analysis results of \cite{Roach1991} are translated to the context of monotone operators. For general Finsler $p$-Laplacian operators, Theorem~\ref{thm:Lqx} seems to be new even if $F=0$.
\end{remark}

\section{The homogenized equation and correctors}\label{sec:hom}
In this section, we introduce the homogenized equation and correctors. We also collect a number of structural properties that these objects satisfy, as well as some basic regularity results concerning the homogenized equation. The results essentially follow from standard theory and methods, but for the sake of completeness, we include proofs.

Throughout this section, we write $\lesssim$  if $\leq $ holds up to a multiplicative positive constant depending only on $p$ and $\Lambda$.

\subsection{Corrector and flux corrector}
We begin by introducing the corrector and flux corrector.
\begin{lemma}\label{L:phi}
Suppose Assumption~\ref{ass:1}  is satisfied. For every $\xi\in \R^{n}$ there exist $\phi_\xi\in W_{\per,0}^{1,p}(Y)$ and $\sigma_\xi\in W_{\per,0}^{1,p'}(Y,\R^{n\times n})$ that are uniquely defined by: 
\begin{equation}\label{eq:convexcorrector}
 \nabla \cdot \bfa(y,\xi+\nabla \phi_\xi(y))=0\qquad\mbox{in $Y$}
\end{equation}
and
\begin{align}
 -\Delta \sigma_{\xi,jk} =& \partial_k J(\xi)\cdot e_j-\partial_j J(\xi)\cdot e_k \qquad\mbox{ in $Y$}\label{eq:sigma},\end{align}
where
$$
J(\xi):=\bfa(\cdot,\xi+\nabla \phi_\xi)-\fint_Y \bfa(y,\xi+\nabla \phi_\xi)\dy.
$$
Moreover, it holds
 \begin{align}\label{prop:sigma2}
            \sigma_{\xi,jk}=-\sigma_{\xi,kj}\quad\mbox{and}\quad  -\sum_{k=1}^n\partial_k \sigma_{\xi,jk}=J(\xi)\cdot e_{j}.
    \end{align}
\end{lemma}
\begin{proof}
The existence and uniqueness of the corrector follows from general theory of monotone operators, e.g. \cite{Zeidler}. The growth condition \eqref{ass:1:growth} and $\phi\in W_{\operatorname{per},0}^{1,p}(Y)$ imply $\bfa(\cdot,\xi+\nabla\phi_\xi)\in L^{p'}(Y)$ and thus the existence of a unique $\sigma_\xi\in W_{\per,0}^{1,p'}(Y,\R^{d\times d}) $ satisfying \eqref{eq:sigma} follows by the $L^p$-theory for Poisson equation. The properties \eqref{prop:sigma2} follow exactly as in the linear case (see e.g.\ \cite[Proposition~2.1.1]{Shen}).
%
%
\end{proof}

Next, we establish needed estimates for the corrector $\phi_\xi$ and flux corrector $\sigma_\xi$.
\begin{lemma}\label{L:Dphi}
 Suppose $\bfa$ satisfies Assumption~\ref{ass:1}. For every $\xi\in\R^{n}$ we have: There is $c=c(p,\Lambda)>0$ such that
\begin{equation}\label{eq:correctorNatural}
\|\phi_\xi\|_{W^{1,p}(Y)}^p+\|\sigma_\xi\|_{W^{1,p^\prime}(Y)}^{p^\prime}\leq c (\mu+|\xi|)^p.
\end{equation}

If $\bfa$ satisfies Assumption~\ref{ass:homonondeg}, then we may sharpen the estimates to
\begin{align}\label{eq:correctorControlled}
\|V_{p}(\phi_\xi)\|_{W^{1,2}(Y)}^2+\|V_{p^\prime}(\sigma_\xi)\|_{W^{1,2}(Y)}^2\leq c |V_{p}(\xi)|^2.
\end{align}
\end{lemma}

\begin{proof}[Proof of Lemma~\ref{L:Dphi}]

\step 1 Estimates for \eqref{eq:correctorNatural}. For $\xi\in\R^n$, we set $F=\xi+\nabla \phi_\xi$. It suffices to show
\begin{align}\label{eq:calc11}
\|\mu+|F|\|_{L^p(Y)}\lesssim \mu+|\xi|.
\end{align}
Indeed, \eqref{eq:calc11} and the growth conditions of $\bfa$ imply
\begin{align}\label{eq:calc111}
\|\bfa(\cdot,F)\|_{L^{p'}(Y)}\leq\Lambda \|\mu+|F|\|_{L^p(Y)}^{p-1}\lesssim (\mu+|\xi|)^{p-1}.
\end{align}
Estimate \eqref{eq:correctorNatural} follows from \eqref{eq:calc11} in combination with triangle inequality and Poincar\'e's inequality with $\fint_Y \phi_\xi = 0$ (for the corrector $\phi_\xi$) and \eqref{eq:calc111} together with the fact that the map $T:L_{\per,0}^{p'}(Y,\R^n)\to W_{\per,0}^{1,p'}(Y,\R)$ given by $TF:=u$ with $-\Delta u=\nabla\cdot F$, is linear and bounded. Hence, it remains to show \eqref{eq:calc11}. For $p\in[2,\infty)$, we have
\begin{align}\label{eq:calc1}
\Lambda^{-1}\int_Y |F|^p\dy\leq& \int_Y \langle \bfa(y,F)-\bfa(y,0),F\rangle\dy\stackrel{\eqref{eq:convexcorrector}}{=} \int_Y \langle \bfa(y,F),\xi\rangle-\langle \bfa(y,0),F\rangle\dy\notag\\
 \leq& \Lambda\int_Y (\mu+|F|)^{p-1}|\xi|+\mu^{p-1}|F|\dy
\end{align}
and \eqref{eq:calc11} follows by Young's inequality. The case $p\in(1,2)$ follows in the same way using the pointwise inequality \eqref{ineq:pointwisesubq} (with $\nabla w$ replaced by $F$).

\step 2 Estimates for \eqref{eq:correctorControlled}. Fix $\xi\in \R^n$ and set $F=\xi+\nabla \phi_\xi$. Similar computations as in Step~1 in combination with Lemma~\ref{L:propV} yield
\begin{eqnarray*}
\frac1\Lambda\int_Y |V_{p}(F)|^2\dy&\stackrel{\eqref{ass:1:monotone:strong}}\leq& \int_Y \langle \bfa(y,F),F\rangle\dy\stackrel{\eqref{eq:convexcorrector}}= \int_Y \langle \bfa(y,F),\xi\rangle\dy \stackrel{\eqref{ass:acont}}\lesssim \int_Y (1+|F|)^{p-2}|F||\xi|\dy\notag\\
&\stackrel{\eqref{eq:YoungV}}\leq& \frac1{2\Lambda} \int_Y |V_{p}(F)|^2\dy + C(\Lambda,p) |V_{p}(\xi)|^2
\end{eqnarray*}
and thus
\begin{equation}\label{eq:calc2}
\|V_{p}(F)\|_{L^{2}(Y)}^2\lesssim  |V_{p}(\xi)|^2.
\end{equation}
Estimate \eqref{eq:calc2}, together with the triangle type inequality \eqref{eq:triangleV} and \eqref{eq:PoincareV} yield $\|V_{p}(\phi_\xi)\|_{W^{1,2}(Y)}^2\lesssim  |V_{p}(\xi)|^2$.
Moreover, \eqref{eq:calc2} imply $\|V_{p'}(J(\xi))\|_{L^{2}(Y)}^2\lesssim  |V_{p}(\xi)|^2$ and the claimed estimate for $\sigma_\xi$ follows from Lemma~\ref{lem:sigmaDecomposition} below.
\end{proof}

Next, we provide needed estimates for the Laplace equation on the torus:
\begin{lemma}\label{lem:sigmaDecomposition}
Let $f\in L^1(Y,\R^n)$. Suppose $-\Delta u = \nabla \cdot f$ in $Y$ with periodic boundary conditions. For all $p\in(1,\infty)$ there exists $c=c(n,p)>0$ such that
\begin{align*}
\|V_{p}(\nabla u)\|_{L^2(Y)}\leq c \|V_{p}(f)\|_{L^2(Y)}.
\end{align*}
\end{lemma}

\begin{proof}
Throughout the proof we write $\lesssim$ if $\leq$ holds up to a positive multiplicative constant depending only on $n$ and $p$. For $p\geq2$, the claim follows directly from the usual $L^q$-estimates for $\nabla u$:
$$
\|V_{p}(\nabla u)\|_{L^2(Y)}^2\lesssim \int_Y|\nabla u|^2+|\nabla u|^p\dx\lesssim\int_Y|f|^2+|f|^p\dx\lesssim \|V_{p}(f)\|_{L^2(Y)}^2.
$$
Suppose that $p\in(1,2)$. Clearly, we have $u=v+w$, where $v,w$ are the unique periodic solutions of
$$
-\Delta v=\nabla \cdot (\mathds 1_{\{|f|\leq1\}}f)\quad\mbox{and}\quad -\Delta w=\nabla \cdot (\mathds 1_{\{|f|>1\}}f)
$$
Combining $L^2$ estimates for $\nabla v$ with the elementary inequality $1\leq 2^{2-p}(1+|f|)^{p-2}$ if $0<|f|\leq1$, we have
\begin{align*}
\|V_{p}(\nabla v)\|_{L^2(Y)}^2\leq \int_Y |\nabla v|^2\dx\leq \int_Y| (\mathds 1_{\{|f|\leq1\}}f)|^2\dx\leq 2^{2-p} \int_Y(1+|f|)^{p-2}|f|^2\dy=\|V_{p}(f)\|_{L^2(Y)}^2.
\end{align*}
Similarly, we obtain with $L^p$-estimates for $\nabla w$ and $|f|^{p-2}\leq 2^{2-p}(1+|f|)^{p-2}$ if $|f|\geq1$ that
\begin{align*}
\|V_{p}(\nabla w)\|_{L^2(Y)}^2\leq \int_Y |\nabla w|^p\dx\lesssim \int_Y| (\mathds 1_{\{|f|>1\}}f)|^p\dx\leq 2^{2-p} \int_Y(1+|f|)^{p-2}|f|^2\dy=\|V_{p}(f)\|_{L^2(Y)}^2.
\end{align*}
The claimed estimate for $V_{p}(\nabla u)$ now follows by \eqref{eq:triangleV}.
\end{proof}

\subsection{The homogenized operator}\label{sec:operator}

For $\bfa$ satisfying Assumption~\ref{ass:1}, we introduce the homogenized operator $\overline\bfa\colon \R^n\to \R^n$  by setting for $\xi\in \R^n$,
\begin{align}\label{def:ahom}
\overline\bfa(\xi)= \int_Y \bfa(y,\xi+\nabla \phi_\xi)\,\dy,
\end{align}
where $\phi_\xi$ is defined via \eqref{eq:convexcorrector}. The main aim of this section is to prove that if $\bfa$ satisfies Assumption~\ref{ass:standard}, then $\overline \bfa$ satisfies Assumption~\ref{ass:1} and some continuity condition, see \eqref{eq:continuityHomogenised}. Moreover, if $\bfa$ satisfies Assumption~\ref{ass:standard} with $\mu=1$, then $\overline \bfa$ satisfies Assumption~\ref{ass:homonondeg}. This results essentially follow from \cite{Fusco1986} and in particular \cite[Section 7]{DalMaso1990}. However, the framework of \cite{DalMaso1990} is significantly more general than our set-up and so we choose to present a self-contained presentation including proofs.



\begin{lemma}\label{lem:stability}
Suppose $\bfa\colon \R^n\times \R^n\to \R^n$ satisfies Assumption~\ref{ass:standard} for some $1<p<\infty$, $\mu\in[0,1]$ and $\Lambda\in[1,\infty)$. Then $\overline \bfa\colon \R^n\to\R^n$ satisfies Assumption~\ref{ass:1} with the same $1<p<\infty$, $\mu\in[0,1]$, and some $\overline\Lambda=\overline\Lambda(\Lambda,p)\in[1,\infty)$. Moreover, there is $c=c(\Lambda,p)\in[1,\infty)$ such that
\begin{align}\label{eq:continuityHomogenised}
|\overline \bfa(\xi_1)-\overline \bfa(\xi_2)|\leq c \begin{cases}
(\mu+|\xi_1|+|\xi_2|)^{p-2}|\xi_1-\xi_2| \quad& \text{ if } p\geq 2\\
(\mu+|\xi_1-\xi_2|)^{p-2}|\xi_1-\xi_2| \quad& \text{ if } p\leq 2
\end{cases}.
\end{align}
In the case $\mu=1$, the operator $\overline \bfa$ satisfies Assumption~\ref{ass:homonondeg}.
\end{lemma}

\begin{proof}
The first part of this result is well-known, see e.g. \cite{Fusco1986}, but we provide a proof for completeness. 
For $\xi_j\in\R^n$, we set  $\xi_j+\nabla\phi_{\xi_j}=F_j$ and recall $\int_Y F_j\dx=\xi_j$. 

\step 1 Proof of \eqref{ass:1:monotone} and \eqref{ass:1:monotone:strong}: Suppose $p\geq2$. The convexity of $\R^n\ni z\mapsto (\mu+|z|)^{p-2}|z|^2$ and Jensen inequality imply
\begin{eqnarray*}
(\mu+|\xi_1-\xi_2|)^{p-2}|\xi_1-\xi_2|^2&\leq& \int_Y (\mu+|F_1-F_2|)^{p-2}|F_1-F_2|^2\,\dy\\
&\stackrel{p\geq2, \eqref{eq:coer:intro}}\leq&\Lambda\int_Y \langle \bfa(y,F_1)-\bfa(y,F_2),F_1-F_2 \rangle\,\dy\\
&\stackrel{\eqref{eq:convexcorrector},\eqref{def:ahom}}=&\Lambda \langle \overline \bfa(\xi_1)-\overline \bfa(\xi_2),\xi_1-\xi_2\rangle.
\end{eqnarray*}
%
In case $p\leq 2$, we instead find using H\"older inequality,
\begin{eqnarray*}
|\xi_1-\xi_2|^2&\leq& \biggl(\int_Y|F_1-F_2|^p\dy\biggr)^\frac2p\\
&\leq& \biggl(\int_Y|F_1-F_2|^2(\mu+|F_1|+|F_2|)^{p-2}\dy\biggr)\biggl(\int_Y (\mu+|F_1|+|F_2|)^p\biggr)^\frac{2-p}p\\
&\leq&\biggl(\Lambda \int_Y \langle \bfa(y,F_1)-\bfa(y,F_2),F_1-F_2\rangle\,\dy\biggr)\biggl(\int_Y (\mu+|F_1|+|F_2|)^p\biggr)^\frac{2-p}p\\
&\stackrel{\eqref{eq:correctorNatural}}\lesssim&\langle \overline\bfa(\xi_1)-\overline\bfa(\xi_2),\xi_1-\xi_2\rangle(\mu+|\xi_1|+|\xi_2|)^{2-p}.
\end{eqnarray*}
Clearly, the above two estimates imply that $\overline \bfa$ satisfies \eqref{ass:1:monotone} and in the case $\mu=1$ even \eqref{ass:1:monotone:strong}

\step 2 Proof of \eqref{ass:1:growth}: We estimate
\begin{align*}
\lvert \overline \bfa(\xi_1)\rvert \leq \int_Y \lvert \bfa(y,F_1)\rvert \,\dy\leq \Lambda \int_Y (\mu+\lvert F_1\rvert)^{p-1}\,\dy\leq\Lambda \left(\int_Y (\mu+\lvert F_1\rvert)^p\dx\right)^\frac{p-1} p\stackrel{\eqref{eq:correctorNatural}}\lesssim (\mu+\lvert \xi_1\rvert)^{p-1}.
\end{align*}

\step 3 Proof of \eqref{eq:continuityHomogenised}: We note that as $\bfa(\cdot,0)=0$, it is immediate that $\overline \bfa(0)=0$ from the definition.

In case $p\geq 2$, the claim is shown in \cite[Proposition 2.7]{Fusco1986}. In the case $p\leq 2$, we first observe
\begin{eqnarray*}
|\overline \bfa(\xi_1)-\overline \bfa(\xi_2)|^{p^\prime}&\leq& \int_Y |\bfa(y,F_1)-\bfa(y,F_2)|^{p^\prime}\dy\lesssim \int_Y (\mu+|F_1|+|F_2|)^{(p-2)p^\prime} |F_1-F_2|^{p^\prime}\dy\\
&\stackrel{p^\prime\geq 2}{\leq}& \int_Y (\mu+|F_1|+|F_2|)^{p-2}|F_1-F_2|^2
\lesssim \int_Y \langle \bfa(y,F_1)-\bfa(y,F_2),F_1-F_2\rangle\\
&\leq&  |\overline \bfa(\xi_1)-\overline \bfa(\xi_2)| |\xi_1-\xi_2|
\end{eqnarray*}
and thus
\begin{align}\label{eq:standard1}
|\overline \bfa(\xi_1)-\overline \bfa(\xi_2)|\lesssim  |\xi_1-\xi_2|^{p-1}.
\end{align}
Moreover, we have
\begin{align*}
|\overline \bfa(\xi_1)-\bfa(\xi_2)|^2 \leq& \int_Y |\bfa(y,F_1)-\bfa(y,F_2)|^2\dy
\lesssim  \int_Y (\mu+|F_1|+|F_2|)^{2(p-2)}|F_1-F_2|^2\dy\\
\leq& \mu^{p-2}\int_Y (\mu+|F_1|+|F_2|)^{(p-2)}|F_1-F_2|^2\dy\\
\lesssim& \mu^{p-2}\int_Y \langle \bfa(y,F_1)-\bfa(y,F_2),F_1-F_2\rangle\dy\\
\leq& \mu^{p-2}|\overline \bfa(\xi_1)-\bfa(\xi_2)| |\xi_1-\xi_2|.
\end{align*}
Combining the above inequality with \eqref{eq:standard1}, we obtain
\begin{align*}
|\overline \bfa(\xi_1)-\overline \bfa(\xi_2)|\lesssim \min\{|\xi_1-\xi_2|^{p-2},\mu^{p-2}\} |\xi_1-\xi_2|\lesssim (\mu+|\xi_1-\xi_2|)^{p-2}|\xi_1-\xi_2|.
\end{align*}
This concludes the proof.
\end{proof}

Finally, we state an estimate for the difference of correctors assuming that $\bfa $ and $\overline \bfa$ satisfy  Assumption~\ref{ass:homonondeg}.

\begin{corollary}
 Suppose $\bfa:\R^n\times\R^n\to\R^n$ and the homogenized operator $\overline \bfa:\R^n\to\R^n$ satisfies Assumption~\ref{ass:homonondeg} for some joint constant $1<p<\infty$ and $\Lambda\in[1,\infty)$. Then there exists $c=c(p,\Lambda)>0$ such that for any $\xi_1,\xi_2\in \R^n$ and $\tau>0$,

\begin{align}\label{eq:correctorDifference}
\|V_{p}(\nabla\phi_{ \xi_1}-\nabla \phi_{\xi_2})\|_{L^2(Y)}^2 \leq c\begin{cases}\left(1+|\xi_1|+|\xi_2|\right)^{p-2}|\xi_1-\xi_2|^2& \text{ if } p\geq 2\\ (|V_{p}(\xi_1-\xi_2)|^2(1+\tau^{-\frac{2-p}p})+\tau(|V_{p}(\xi_1)|+|V_{p}(\xi_2)|)^2& \text{ if } p\leq 2
\end{cases}
\end{align}
\end{corollary}
\begin{proof}
As before, we denote for $i=1,2$, $F_i = \xi_i + \nabla \phi_{\xi_i}$. Then
\begin{align*}
\Lambda^{-1}\int_Y |W_p(F_1,F_2)|^2\dy\leq& \int_Y \langle \bfa(y,F_1)-\bfa(y,F_2),F_1-F_2\rangle\dy= \int_Y \langle \bfa(y,F_1)-\bfa(y,F_2),\xi_1-\xi_2\rangle\dy\\
=& \langle \overline \bfa(\xi_1)- \overline \bfa(\xi_2),\xi_1-\xi_2\rangle.
\end{align*}
In the case $p\geq2$, the claim follows with help of \eqref{ass:1:monotone:strong}, $|W_p(F_1,F_2)|^2=|V_{p}(F_1-F_2)|^2$ and the triangle inequality \eqref{eq:triangleV} in the form
$$
|V_{p}(\nabla \phi_{\xi_1}-\nabla \phi_{\xi_2})|^2\lesssim |V_{p}(F_1-F_2)|^2+|V_{p}(\xi_1-\xi_2)|^2\lesssim  |V_{p}(F_1-F_2)|^2+(1+|\xi_1|+|\xi_2|)^{p-2}|\xi_1-\xi_2|^2.
$$
In the case $p\in(1,2)$, we use in addition \eqref{est:WpVsub} and the estimates \eqref{eq:correctorControlled}.
\end{proof}

\subsection{A stable assumption under homogenization}\label{sec:anisotropic} In this section, we recall some (dual) monotonicity conditions, see Assumption~\ref{ass:dualAssumption} below, that are a special case of \cite[Section 7]{{DalMaso1990}}. These assumptions are satisfied by examples that are not covered by Assumption~\ref{ass:standard}, see Remark~\ref{rem:anisotropic2}. Moreover, they have the interesting property that they are closed under homogenization, see Lemma~\ref{lem:ahom} below, and ensure enough regularity on $\bfa$ and $\overline \bfa$ such that the large-scale regularity theory developed in Section~\ref{sec:largescaleCZ} and Section~\ref{sec:largescaleLip} applies.

\begin{assumption}\label{ass:dualAssumption}
Let $1<p<\infty$ and $\Lambda\in[1,\infty)$ be given. Suppose that $\bfa:\R^{n}\times \R^n\to\R^{n}$ satisfies for all $\xi_1,\xi_2\in \R^{n}$ that $\bfa(x,0)=0$,
\begin{align}\label{ass:2:monotone:strong}
\Lambda\langle \bfa(x,\xi_1)-\bfa(x,\xi_2),\xi_1-\xi_2\rangle \geq 
\begin{cases}
(\mu+|\xi_1-\xi_2|)^{p-2}|\xi_1-\xi_2|^2 &\text{ if } p\geq 2\\
(\mu+\lvert \xi_1\rvert +\lvert \xi_2\rvert)^{p-2}\lvert \xi_1-\xi_2\rvert^2 &\text{ if } p<2
\end{cases}
\end{align}
and 
\begin{align}\label{ass:2:dualmonotone}
\Lambda\langle \bfa(x,\xi_1)-\bfa(x,\xi_2),\xi_1-\xi_2\rangle \geq 
\begin{cases}
\left(\mu^{p-1}+\lvert \bfa(x,\xi_1)\rvert+\lvert \bfa(x,\xi_2)\rvert\right)^{p'-2}|\bfa(x,\xi_1)-\bfa(x,\xi_2)|^2 &\text{ if } p\geq 2\\
(\mu^{p-1}+\lvert \bfa(x,\xi_1)-\bfa(x,\xi_2)\rvert)^{p^\prime-2}|\bfa(x,\xi_1-\bfa(x,\xi_2)|^2 &\text{ if } p\leq 2
\end{cases}.
\end{align}
\end{assumption}

\begin{remark}\label{rem:anisotropic2}
Similar to our comments in Remark \ref{rem:anistropic1}, Assumption~\ref{ass:dualAssumption} and the results of this section show that our homogenization results in Section~\ref{sec:largescaleCZ} cover certain oscillatory Finsler $p$-Laplacian operators.  As a specific example the orthotropic p-Laplace operator with oscillating coefficients, that is $\bfa(y,z)= A(y)\sum_{i=1}^n|z_i|^{p-2}z_i$ with $p\in (1,\infty)$ and a uniformly elliptic coefficient matrix $A\in L^\infty(Y,\R^{n\times n})$ satisfies Assumption~\ref{ass:dualAssumption}.  In fact, whenever $\|\cdot\|\colon \R^n\to \R$ is a $\min\{p,2\}$-smooth and $\max\{p,2\}$-convex norm, then the Finsler $p$-Laplacian $\bfa(z)=\|z\|^{p-1}\nabla(\|z\|)$ satisfies \eqref{ass:dualAssumption}, see \cite{Roach1991}. 
 \end{remark}

We think of \eqref{ass:2:dualmonotone} as a coercivity assumption on $\bfa^{-1}$, even when $\bfa^{-1}$ is not formally well-defined. This point of view is present in the classical literature regarding the simpler setting of quadratic growth $p=2$ \cite{Zhikov1994} and has gained growing interest for proving regularity properties of convex minimisation problems \cite{Carozza2015,Cianchi2019,deFilippis2024}. Here we implement these ideas in the context of homogenization. It is instructive to consider the special case of the $p$-Laplacian $\bfa(\xi)=|\xi|^{p-2}\xi$. Note that in that case
\begin{align*}
\langle |\xi_1|^{p-2}\xi_1-|\xi_2|^{p-2}\xi_2,\xi_1-\xi_2\rangle = \langle \bfa(\xi_1)-\bfa(\xi_2),|\bfa(\xi_1)|^{p^\prime-2}\bfa(\xi_1)-|\bfa(\xi_2)|^{p^\prime-2}\bfa(\xi_2)\rangle.
\end{align*}
Hence the dual nature of the relationship betwen \eqref{ass:2:monotone:strong} and \eqref{ass:2:dualmonotone} is apparent in this case. It is therefore natural to expect that when $\bfa(x,z)=\partial_z F(x,z)$ for $F$, convex as a function of $z$, \eqref{ass:2:dualmonotone} is equivalent to a growth condition on $\bfa$.
\begin{lemma}\label{lem:Ass3vs4}
Suppose $\bfa\colon \R^n\times \R^n\to \R^n$ satisfies Assumption~\ref{ass:dualAssumption} for some $1<p<\infty$, $\mu\in[0,1]$ and $\Lambda\in [1,\infty)$. Then there exists $\Lambda'=\Lambda'(\Lambda,p)\in[1,\infty)$ such that 
\begin{align}\label{eq:continuityHomogenised2}
|\bfa(x,\xi_1)-\bfa(x,\xi_2)|\leq \Lambda' \begin{cases}
(\mu+|\xi_1|+|\xi_2|)^{p-2}|\xi_1-\xi_2| \quad& \text{ if } p\geq 2\\
(\mu+|\xi_1-\xi_2|)^{p-2}|\xi_1-\xi_2| \quad& \text{ if } p\leq 2
\end{cases}.
\end{align}
In particular, Assumption~\ref{ass:dualAssumption} with $1<p<\infty$, $\mu=1$ and $\Lambda\in [1,\infty)$ implies Assumption~\ref{ass:homonondeg} with the same $1<p<\infty$ and a suitable $\Lambda'=\Lambda'(\Lambda,p)\in[1,\infty)$. Moreover, if $\bfa(x,z)= \partial_z F(x,z)$ for some $F\colon\R^n\times \R^n \to \R^n$, then Assumption~\ref{ass:homonondeg} implies Assumption~\ref{ass:dualAssumption} with $\mu=1$ upon changing the value of  $\Lambda$.
\end{lemma}
\begin{proof}

Clearly, it suffices to show that \eqref{ass:2:dualmonotone} implies \eqref{eq:continuityHomogenised2} and that Assumption~\ref{ass:homonondeg} implies \eqref{ass:2:dualmonotone} with $\mu=1$ in the case $\bfa(x,z)=\partial_z F(x,z)$.

\step 1 We show that \eqref{ass:2:dualmonotone} implies \eqref{eq:continuityHomogenised2}.  Condition \eqref{ass:2:dualmonotone} and Cauchy-Schwarz inequality imply
\begin{align}\label{ass:2:dualmonotone2}
\Lambda |\xi_1-\xi_2| \geq 
\begin{cases}
\left(\mu^{p-1}+\lvert \bfa(x,\xi_1)\rvert+\lvert \bfa(x,\xi_2)\rvert\right)^{p^\prime-2}|\bfa(x,\xi_1)-\bfa(x,\xi_2)| \quad&\text{ if } p\geq 2\\
(\mu^{p-1}+| \bfa(x,\xi_1)-\bfa(x,\xi_2)|)^{p'-2}|\bfa(x,\xi_1)-\bfa(x,\xi_2)| \quad&\text{ if } p\leq 2
\end{cases}.
\end{align}
Estimate \eqref{ass:2:dualmonotone2} for $\xi_1=\xi\in\R^n$ and $\xi_2=0$ implies for $p\in(1,2)$ 
\begin{align*}
|\bfa (x,\xi)|^\frac1{p-1}=|\bfa(x,\xi)|^\frac{2-p}{p-1}|\bfa(x,\xi)|\leq (\mu^{p-1}+|\bfa(x,\xi)|)^{p'-2}|\bfa(x,\xi)|\leq \Lambda|\xi|
\end{align*}
and for $p\geq2$
\begin{align*}
(\mu^{p-1}+|\bfa (x,\xi)|)^\frac1{p-1}=(\mu^{p-1}+|\bfa(x,\xi)|)^\frac{2-p}{p-1}(\mu^{p-1}+|\bfa(x,\xi)|)\leq \mu+ \Lambda|\xi|.
\end{align*}
Hence,
\begin{align}\label{eq:easyUpper}
|\bfa(x,\xi)|\leq 2\Lambda^{p-1}(\mu+|\xi|)^{p-1}.
\end{align}
In the case $p\in[2,\infty)$, \eqref{ass:2:dualmonotone2} in combination with \eqref{eq:easyUpper}, $0\leq 2-p'=\frac{p-2}{p-1}$ and the fact that $t\mapsto t^\frac1{p-1}$ is subadditive imply the desired claim. In the case $p\in(1,2)$, inequality \eqref{ass:2:dualmonotone2} implies
\begin{align*}
|\bfa (x,\xi_1)-\bfa (x,\xi_2)|\leq \Lambda |\xi_1-\xi_2|(\max\{\mu,\Lambda|\xi_1-\xi_2|\})^{p-2}\leq 2^{2-p}\Lambda (\mu+|\xi_1-\xi_2|)^{p-2}|\xi_1-\xi_2|,
\end{align*}
where we use in the last inequality $\frac12(a+b)\leq \max\{a,b\}$ and $\Lambda\geq1$.

\step 2 If $\bfa(x,z)=\partial_z F(x,z)$ for some $F\colon \R^n\times \R^n \to \R^n$, denote the Fenchel conjugate of $F$ by $F^\ast$. Then by a straightforward adaption of the discussion in \cite[Section 1.1.]{KochOT}, for any $\sigma_1,\sigma_2\in \R^n$ it holds that
\begin{align}\label{eq:dualLower}
\langle (F^\ast)^\prime(\sigma_1)-(F^\ast)^\prime(\sigma_2),\sigma_1-\sigma_2\rangle\gtrsim& \begin{cases}
(\mu^{p-1}+|\sigma_1|+|\sigma_2|)^{p^\prime-2}|\sigma_1-\sigma_2|^2 \quad &\text{ if } p\geq 2\\
(\mu^{p-1}+|\sigma_1-\sigma_2|)^{p^\prime-2}|\sigma_1-\sigma_2|^2 \quad &\text{ if } p\leq 2.
\end{cases}
\end{align}
As $F$ is superlinear, strictly convex and $C^1$, $F^\prime$ and $(F^\ast)^\prime$ are homeomorphisms of $\R^n$. Moreover, it holds that $(F^\ast)^\prime(F^\prime(z))=z$ for all $z\in \R^n$ \cite[Theorem 26.5, 26.6]{Rockafellar}. In particular, with the choice $\sigma_i = F^\prime(z_i)$ for $i=1,2$, \eqref{eq:dualLower} gives \eqref{ass:2:dualmonotone}.
\end{proof}

The advantage of Assumption~\ref{ass:dualAssumption} is that it is stable under homogenization.
\begin{lemma}\label{lem:ahom}
Suppose $\bfa\colon \R^n\times \R^n\to \R^n$ satisfies Assumption~\ref{ass:dualAssumption} for some $1<p<\infty$, $\mu\in[0,1]$ and $\Lambda\in[1,\infty)$. Then $\overline \bfa\colon \R^n \to \R^n$ satisfies Assumption~\ref{ass:dualAssumption} with the same $1<p<\infty$, $\mu\in[0,1]$ and for some $\overline \Lambda=\overline \Lambda(\Lambda,p)\in[1,\infty)$.
\end{lemma}

\begin{proof}

Clearly, we have $\overline \bfa(0)=0$ and Step~1 of the proof of Lemma \ref{lem:stability} shows that $\overline \bfa$ satisfies \eqref{ass:2:monotone:strong} replacing $\Lambda$ with $\overline \Lambda(\Lambda,p)\in[1,\infty)$. Hence, it remains to prove that $\overline \bfa$ satisfies \eqref{ass:2:dualmonotone} (upon changing $\Lambda$).\\
Fix $p\in[2,\infty)$. We have
\begin{eqnarray*}
|\overline \bfa(\xi_1)-\overline \bfa(\xi_2)|^2&\leq& \biggl(\int_Y|\bfa(y,F_1)-\bfa(y,F_2)|^{p'}\dy\biggr)^\frac2{p'}\\
&\leq& \biggl(\int_Y\left(\mu^{p-1}+\lvert \bfa(x,F_1)\rvert+\lvert \bfa(x,F_2)\rvert\right)^{p'-2}|\bfa(x,F_1)-\bfa(x,F_2)|^2\dy\biggr)\\
&&\times\biggl(\int_Y (\mu^{p-1}+|\bfa(x,F_1)|+|\bfa(x,F_2)|)^{p'}\biggr)^\frac{p-2}p\\
&\stackrel{\eqref{ass:2:dualmonotone},\eqref{eq:easyUpper}}\lesssim&\biggl( \int_Y \langle \bfa(y,F_1)-\bfa(y,F_2),F_1-F_2\rangle\,\dy\biggr)\biggl(\int_Y (\mu+|F_1|+|F_2|)^p\biggr)^\frac{p-2}p\\
&\stackrel{\eqref{eq:correctorNatural}}\lesssim&\langle \overline\bfa(\xi_1)-\overline\bfa(\xi_2),\xi_1-\xi_2\rangle(\mu+|\xi_1|+|\xi_2|)^{p-2}.
\end{eqnarray*}
The desired inequality follows from the facts that $\overline \bfa$ satisfies $\overline \bfa(0)=0$ and \eqref{ass:2:monotone:strong} (with $\overline \Lambda$) which imply $|\xi_i|^{p-1}\leq \overline \Lambda |\overline\bfa (\xi_i)|$. For $p\in(1,2)$, we  use convexity of $z\mapsto (\mu^{p-1}+|z|)^{p'-2}|z|^2$, \eqref{def:ahom} and Jensen inequality to obtain
\begin{eqnarray*}
& &(\mu^{p-1}+\lvert \overline\bfa(\xi_1)-\overline\bfa(\xi_2)\rvert)^{p^\prime-2}\lvert \overline \bfa(\xi_1)-\overline \bfa(\xi_2)\rvert^{2}\\
&\leq & \int_Y (\mu^{p-1}+|\bfa(x,F_1)-\bfa(y,F_2)|)^{p'-2}|\bfa(y,F_1)-\bfa(y,F_2)|^2\,\dy\\
&\stackrel{\eqref{ass:2:dualmonotone}}\leq &\Lambda \int_Y \langle \bfa(y,F_1)-\bfa(y,F_2),F_1-F_2\rangle\,\dy\\
&\stackrel{\eqref{eq:convexcorrector},\eqref{def:ahom}}=& \Lambda\langle \overline \bfa(\xi_1)-\overline \bfa(\xi_2),\xi_1-\xi_2\rangle,
\end{eqnarray*}
which shows that $\overline \bfa$ satisfies \eqref{ass:2:dualmonotone} (with the same $\Lambda$) in this case.
\end{proof}

\section{Calderon-Zygmund theory on large scales}\label{sec:largescaleCZ}
The main result of this section is a \textit{large-scale} Calderon-Zygmund estimate from which Theorem~\ref{thm:Lqx:uniformeps} follows. For $f\in L^1(\R^n)$, we recall the definition of the maximal function and define a truncated maximal function as follows
\begin{equation}
M(f)(x):=\sup_{\rho>0}\fint_{B_\rho(x)}f(y)\dy\quad\mbox{and}\quad M_r(f)(x):=\sup_{\rho\geq r}\fint_{B_{\rho}(x)}f(y)\dy.
\end{equation}

\begin{theorem}\label{L:nonlinear:largescalereg:lp}
Fix $1<p<\infty$, $\mu\in[0,1]$ and $\Lambda\in[1,\infty)$. Suppose $\bfa:\R^n\times\R^n\to\R^n$ is periodic in the first variable. Assume that $\bfa$ and the corresponding homogenized operator $\overline \bfa:\R^n\to\R^n$, defined by \eqref{def:ahom}, satisfy for all $x,\xi,\xi_1,\xi_2\in\R^n$ \eqref{ass:1:monotone}, \eqref{ass:1:growth} and
\begin{align}\label{eq:continuityHomogenised1}
|\bfa(x,\xi_1)-\bfa(x,\xi_2)|+|\overline \bfa(\xi_1)-\overline \bfa(\xi_2)|\leq \Lambda \begin{cases}
(\mu+|\xi_1|+|\xi_2|)^{p-2}|\xi_1-\xi_2| \quad& \text{ if } p\geq 2\\
(\mu+|\xi_1-\xi_2|)^{p-2}|\xi_1-\xi_2| \quad& \text{ if } p\leq 2.
\end{cases}
\end{align}

 For every $q\in(p,\infty)$ there exists $C=C(\Lambda,n,p,q)\in[1,\infty)$ such that the following is true. Let $u\in W^{1,p}(B)$ and $F\in L^\frac{p}{p-1}(B)$ with $B=B_R(x_0)\subset\R^n$ and $\e>0$ be such that
\begin{equation}\label{eq:T1:Lp}
\nabla\cdot \bfa(\tfrac{x}\e,\nabla u)=\nabla \cdot F\qquad\mbox{in $B$}
\end{equation}
Then, it holds
\begin{equation}\label{L1:lp}
\|M_\e(\mathds 1_B (\mu+|\nabla u|)^p)\|_{\underline L^\frac{q}p(\frac12 B)}\leq C\biggl(\|\mu+|\nabla u|\|_{\underline L^p(B)}^p+\|M_\e(\mathds 1_{B}|F|^{p'})\|_{\underline L^{\frac{q}{p}}(B)}\biggr).
\end{equation}
\end{theorem}

Next, we prove that Theorem~\ref{L:nonlinear:largescalereg:lp} in combination with Theorem~\ref{thm:Lqx} imply Theorem~\ref{thm:Lqx:uniformeps}.

\begin{proof}[Proof of Theorem~\ref{thm:Lqx:uniformeps}]
Note that due to Lemma~\ref{lem:stability}, we are in a situation where we may apply Theorem~\ref{L:nonlinear:largescalereg:lp}.

We distinguish between small and large balls starting with the smaller ones.

\step 1 Suppose that $B=B_R(x_0)$ with $R\leq 4\e$. By rescaling, we have that
\begin{equation*}
\nabla \cdot a(y,\nabla u_\e)=\nabla \cdot F_\e\qquad \mbox{in $B_{R/\e}(x_0/\e)$,}
\end{equation*}
where $u_\e:=\e^{-1}u(\e x)$ and $F_\e(x):=F(\e x)$. Since $R/\e\leq 4$, we deduce from Theorem~\ref{thm:Lqx} that
\begin{equation*}
\|\mu+|\nabla u_\e|\|_{\underline L^q(B_{\frac{R}{2\e}}(x_0/\e))}^q\leq C\biggl(\|\mu+|\nabla u_\e|\|_{\underline L^p(B_{R/\e}(x))}^q+\|F_\e\|_{\underline L^\frac{q}{p-1}(B_{R/\e}(x_0/\e))}^\frac{q}{p-1}\biggr),
\end{equation*}
where $C=C(\Lambda,n,p,q,\omega)>0$ and thus by scaling
\begin{align*}
\|\mu+|\nabla u|\|_{\underline L^q(B_{R/2}(x_0))}^q\leq& C\biggl(\|\mu+|\nabla u|\|_{\underline L^p(B_{R}(x_0))}^q+\|F\|_{\underline L^\frac{q}{p-1}(B_{R}(x_0))}^\frac{q}{p-1}\biggr),
\end{align*}
which proves the statement in that case.

\step 2 Suppose that $B=B_R(x_0)$ with $R> 4\e$. For every $x\in \frac12B$, we use the estimate from Step~1 applied to the ball $B_{2\e}(x)\subset B$ and obtain
\begin{align*}
\|\mu+|\nabla u|\|_{\underline L^q(B_\e(x))}^q\leq& C_1\biggl(\|\mu+|\nabla u|\|_{\underline L^p(B_{2\e}(x))}^q+\|F\|_{\underline L^\frac{q}{p-1}(B_{2\e}(x))}^\frac{q}{p-1}\biggr)\notag\\
\leq&C_1\biggl((M_\e(\mathds 1_{B}(\mu^p+|\nabla u|^p))(x))^\frac{q}p+\|F\|_{\underline L^\frac{q}{p-1}(B_{2\e}(x))}^\frac{q}{p-1}\biggr),
\end{align*}
where $C_1=C_1(\Lambda,n,p,q,\omega)>0$. Moreover, the assumption $R>4\e$ ensures by Fubini
\begin{align*}
\int_{\frac{1}{2}B}\fint_{B_\e(x)}(\mu+|\nabla u(y)|)^q\dy\dx=\frac{1}{|B_\e|}\int_{B}(\mu+|\nabla u(y)|)^q\int_{\frac{1}{2}B}\mathds 1_{B_\e(x)}(y)\dx\dy\geq \int_{\frac{1}{4} B}(\mu+|\nabla u|)^q\dy,
\end{align*}
where we use for the inequality $\int_{\frac{1}{2}B}\mathds 1_{B_\e(x)}(y)\dx=|\frac{1}{2} B\cap B_{\e}(y)|=|B_\e|$ for $y\in\frac{1}{4} B$.

 The above two estimates in combination with \eqref{L1:lp} yield,
\begin{align*}
\|\mu+|\nabla u|\|_{\underline L^q(\frac14B)}^q\leq&2^n\fint_{\frac12B}\fint_{B_\e(x)}(\mu+|\nabla u|)^q\dy\dx\\
\leq& 2^nC_1\biggl(\fint_{\frac12B}(M_\e(\mathds 1_B(\mu^p+|\nabla u|^p))(x))^\frac{q}p+\fint_{B_{2\e}(x)}|F|^\frac{q}{p-1}(y)\dy\dx\biggr)\\
\leq& 2^nC_1C\biggl(\|\mu+|\nabla u|\|_{\underline L^p(B)}^q+\|M_\e(\mathds 1_{B}|F|^{p'})\|_{\underline L^{\frac{q}{p}}(B)}^\frac{q}{p}+\|F\|_{\underline L^\frac{q}{p-1}(B)}^\frac{q}{p-1}\biggr)\\
\leq&2^nC_1Cc(n,\tfrac{q}p)\biggl(\|\mu+|\nabla u|\|_{\underline L^p(B)}^q+\|F\|_{\underline L^\frac{q}{p-1}(B)}^\frac{q}{p-1}\biggr),
\end{align*}
where $C=C(\Lambda,n,p,q)>0$ and we use in the last estimate the strong maximal function inequality. The claimed inequality \eqref{thm:Lqx:uniformeps:claimest} follows by a simple covering argument.
\end{proof}

The proof of Theorem~\ref{L:nonlinear:largescalereg:lp} relies on a combination of Lemma~\ref{lem:CZ} and the following comparison estimate
\begin{lemma}\label{Preg:1:lp}
Consider the situation of Theorem~\ref{L:nonlinear:largescalereg:lp}. There exist constants $c=c(n,\Lambda,p)\in[1,\infty)$ and ${\beta=\beta(\Lambda,n,p)>0}$ such that the following is true: Let $u\in W^{1,p}(B)$ with $B=B_R(x_0)$ and  $\e\in(0,\frac{R}8]$ be such that
\begin{equation}\label{eq:P1}
\nabla\cdot \bfa(\tfrac{x}\e,\nabla u)=0\qquad \text{ in }B.
\end{equation}
Then there exists $v\in u+ W_0^{1,p}(B)$ satisfying,
\begin{equation}\label{P1:est:xi-zeta}
\|\nabla u-\nabla v\|_{\underline L^p(B_{R})}\leq c(\e/R)^\beta\|\mu+|\nabla u|\|_{\underline L^p(B_R)}.
\end{equation}
Moreover, for every $q\in[p,\infty)$ there exists $c_q=c_q(\Lambda,n,p,q)\in[1,\infty)$ such that
\begin{align}\label{P1:est:claim}
\biggl(\fint_{B_{\frac{R}2}}\biggl(\max_{r'\in[\e,\frac{R}4]}\fint_{B_{r'}(x)}|\nabla v|^p\dy\biggr)^\frac{q}p\,dx\biggr)^\frac1q\leq c_q\|\mu+|\nabla u|\|_{\underline L^p(B_R)}.
\end{align}
\end{lemma}
The proof of Lemma~\ref{Preg:1:lp} relies on homogenization estimates in combination with the higher gradient integrability result of Theorem~\ref{thm:Lq} applied to a suitable homogenized equation. We postpone the proof of Lemma~\ref{Preg:1:lp} to the end of the section and start with the proof of Theorem~\ref{L:nonlinear:largescalereg:lp} by following the reasoning of \cite[Lemma 7.6]{Armstrong2019} in the linear case.

\begin{proof}[Proof of Theorem~\ref{L:nonlinear:largescalereg:lp}]

We consider $1<p<q<2q$. Throughout the proof we write $\lesssim$ if $\leq$ holds up to a positive multiplicative constant depending only on $\Lambda,n$ and $p$. We divide the proof in several steps. In Step~1--Step~4, we prove the statement in the case that $B=B_1$ and in the final Step~4, we argue that by scaling this implies the claim for all balls $B$.  

\step 1 Preparations. We consider the case that $B=B_1$ and $u\in W^{1,p}(B_1)$ solves \eqref{eq:T1:Lp}. Without loss of generality, we assume that $(u)_{B_1}=0$ and by the Sobolev extension theorem we can assume that $u\in W^{1,p}(\R^n)$ and
\begin{equation}\label{pf:Th:largescaleCZ:usctend}
\|\nabla u\|_{L^p(\R^n)}\lesssim \|u\|_{W^{1,p}(B_1)}\lesssim \|\nabla u\|_{L^{p}(B_1)},
\end{equation}
where the last inequality is a consequence of Poincar\'e inequality and $(u)_{B_1}=0$. Fix $\delta\in(0,1]$ and let $1\leq N=N(\delta,\Lambda,n,p)\in\mathbb N$ be the smallest integer satisfying $cN^{-\beta}\leq \delta$, where $c=c(\Lambda,n,p)\in[1,\infty)$ and $\beta=\beta(\Lambda,n,p)>0$ are as in Lemma~\ref{Preg:1:lp}.  For $z\in \R^{n}$, define
\begin{equation}\label{pf:largescaleCZ:deff}
f(z):=\mu+\max_{r'\in[N\e,\frac14]}\frac 1 {r'}\|u-(u)_{B_{r'}(z)}\|_{\underline L^p(B_{r'}(z))}+\|\nabla u\|_{\underline L^p(B_1)}.
\end{equation}
We claim that there exist $C=C(\Lambda,n,p)>0$ and $C_q=C(\Lambda,n,p,q)>0$ such that for every $x\in B_\frac12$ and $r\in(0,\frac14]$ there exists $f_{x,r}$ satisfying
\begin{align}
\|f_{x,r}\|_{\underline L^{2q}(B_r(x))}\leq C_q\|f\|_{\underline L^p(B_{4r}(x))}+C\|g_\delta \|_{\underline L^{p'}(B_{4r}(x))}^\frac1{p-1}\label{T:largelp:pf:claim1}\\
\|f-f_{x,r}\|_{\underline L^p(B_r(x))}\leq C\delta^{\min\{1,p-1\}} \|f\|_{\underline L^p(B_{4r}(x))}+C\|g_\delta \|_{\underline L^{p'}(B_{4r}(x))}^\frac1{p-1},\label{T:largelp:pf:claim2}
\end{align}
where 
\begin{equation}\label{L:nonlinear:largescalereg:lp:defg}
g_\delta:=M_\e(\mathds 1_{B_1}|F_\delta|^{p'})^\frac1{p'}\qquad\mbox{with}\qquad F_\delta:=\max\{1,\delta^{p-2}\}F.
\end{equation}
This will be the content of the following two steps.

\step 2 Suppose that $x\in B_\frac12$ and $N\e\geq r$. In this case, we have for all $r'\in[N\e,\frac18]$ that $B_{r'}(z)\subset B_{2r'}(z')$ for all $z,z'\in B_r(x)$ and thus
$$
\sup_{z\in B_{r}(x)}\frac 1 {r'}\|u-(u)_{B_{r'}(z)}\|_{\underline L^p(B_{r'}(z))}\leq 2^{1+\frac{n}p}\inf_{z\in B_{r}(x)}\frac 1 {2r'}\|u-(u)_{B_{2r'}(z)}\|_{\underline L^p(B_{2r'}(z))}.
$$
For $r'\in[\frac18,\frac14]$, we obtain by Poincar\'e inequality
$$
\sup_{z\in B_{r}(x)}\frac 1 {r'}\|u-(u)_{B_{r'}(z)}\|_{\underline L^p(B_{r'}(z))}\lesssim\sup_{z\in B_{r}(x)} \|\nabla u\|_{\underline L^p(B_{r'}(z))}\leq 8^\frac{n}p\|\nabla u\|_{\underline L^p(B_1)}.
$$
Combining the previous two displayed estimates with the definition of $f$, see \eqref{pf:largescaleCZ:deff}, we obtain
\begin{equation}
\sup_{z\in B_r(x)} f(z)\lesssim \inf_{z\in B_r(x)}f(z)\leq \|f\|_{\underline L^p(B_r(x))}.
\end{equation} 
Hence, the choice $f_{x,r}=f$ satisfies \eqref{T:largelp:pf:claim1} and \eqref{T:largelp:pf:claim2}.

\step 3 Suppose $N\e<r$. Define $w_{r}\in W^{1,p}(B_{2r}(x))$ by
$$
w_r-u\in W_0^{1,p}(B_{2r}(x)),\qquad \nabla \cdot\bfa(\tfrac{x}\e,\nabla w_r)=0\quad \mbox{in $B_{2r}(x)$},
$$
and extend $w_{r}$ by $u$ outside of $B_{2r}(x)$. In what follows, we frequently use that for all $t\in(1,\infty)$, $s\geq r$ and $h\in L^t(B)$, we have
\begin{equation}\label{est:strangepointwisetomaximal}
\|h\|_{\underline L^t(B_s(x))}^t\leq 2^n\|M_\e(\mathds 1_B|h|^t)\|_{\underline L^1(B_s(x))}.
\end{equation}
Note that \eqref{est:strangepointwisetomaximal} is a consequence of $B_{s}(x)\subset B_{2s}(y)$ for all $y\in B_{s}(x)$ in the form
\begin{equation*}
\|h\|_{\underline L^{t}(B_{s}(x))}^{t}\leq 2^n\fint_{B_{s}(x)}\|\mathds 1_Bh\|_{\underline L^{t}(B_{2s}(y))}^{t}\dy\leq 2^n\|M_{2s}(\mathds 1_{B}|h|^{t})\|_{\underline L^1(B_{2s}(x))}
\end{equation*}
and $\e\leq 2s$. The comparison estimate \eqref{L:energyestimate:eq2} in Proposition~\ref{prop:regularity} with $\tau=\delta$ together with \eqref{est:strangepointwisetomaximal} and the definitions of $F_\delta$, $g_\delta$ (see \eqref{L:nonlinear:largescalereg:lp:defg}) give
\begin{align}\label{T:CZ:eneuw}
\|\nabla u-\nabla w_{r}\|_{\underline L^p(B_{2r}(x))}\leq& \delta \|\mu+|\nabla u|\|_{\underline L^p(B_{2r}(x)}+c \|F_\delta\|_{\underline L^{p'}(B_{2r}(x))}^\frac1{p-1}\notag\\
\leq& \delta \|\mu+|\nabla u|\|_{\underline L^p(B_{2r}(x)}+c2^\frac np \|g_\delta\|_{\underline L^{p'}(B_{2r}(x))}^\frac1{p-1},
\end{align}
where $c=c(\Lambda,n,p)>0$.

In view of Lemma~\ref{Preg:1:lp} there exists $v_{r}\in w_{r}+W_0^{1,p}(B_{2r}(x))=u+W_0^{1,p}(B_{2r}(x))$ satisfying
\begin{equation*}
\|\nabla w_{r}-\nabla v_{r}\|_{\underline L^p(B_{2r}(x))}\leq (c(\e/r)^\beta)\|\mu+|\nabla w_{r}|\|_{\underline L^p(B_{2r}(x))}\leq \delta\|\mu+|\nabla w_r|\|_{\underline L^p(B_{2r}(x))},
\end{equation*}
where we use $r^{-1}<(\e N)^{-1}$ and $c N^{-\beta}\leq \delta$ in the last inequality. Moreover, there exists a constant ${c_{q}=c_{q}(\Lambda,n,p,q)\in[1,\infty)}$ such that
\begin{align*}
\biggl(\fint_{B_{r}(x)}\biggl(\max_{r'\in[\e,\frac{r}4]}\fint_{B_{r'}(z)}|\nabla v_r|^p\dy\biggr)^\frac{2q}p\,{\mathrm d}z\biggr)^\frac1{2q}\leq c_{q}\|\mu+|\nabla w_r|\|_{\underline L^p(B_{2r}(x))}.
\end{align*}
We extend $v_r$ by $u$ outside of $B_{2r}(x)$. Combining the above two estimates with triangle inequality and \eqref{T:CZ:eneuw}, we find
\begin{equation}\label{das:2}
\|\nabla u-\nabla v_{r}\|_{\underline L^p(B_{2r}(x))}\leq 2\delta\|\mu+|\nabla u|\|_{\underline L^p(B_{2r}(x))}+c\|g_\delta\|_{\underline L^{p^\prime}(B_{2r})}^\frac 1 {p-1},
\end{equation}
and 
\begin{align}\label{est:largelqvr}
\biggl(\fint_{B_{r}(x)}\biggl(\max_{r'\in[\e,\frac{r}4]}\fint_{B_{r'}(z)}|\nabla v_{r}|^p\dy\biggr)^\frac{2q}p\,{\mathrm d}z\biggr)^\frac1{2q}\leq c_q\|\mu+|\nabla u|\|_{\underline L^p(B_{2r}(x))}+c\|g_\delta\|_{\underline L^{p^\prime(B_{2r})}}^\frac 1 {p-1}.
\end{align}
Set
\begin{equation*}
f_{x,r}(z):=\mu+\max_{r'\in[N\e,\frac r4]}\frac 1 {r'}\|v_r-(v_r)_{B_{r'}(z)}\|_{\underline L^p(B_{r'}(z))}+\|\nabla u\|_{\underline L^p(B_1)}.
\end{equation*}
By the Poincar\'{e} inequality, we have,
\begin{align*}
f_{x,r}(z)\lesssim \mu+\max_{r'\in [N\e,\frac r 4]} \|\nabla v_r\|_{\underline L^p(B_{r'}(z))}+\|\nabla u\|_{\underline L^p(B_1)},
\end{align*}
and thus with help of \eqref{est:largelqvr} 
\begin{align}
\biggl(\fint_{B_{r}(x)}|f_{x,r}|^{2q}\,{\mathrm d}z\biggr)^\frac1{2q}\lesssim c_q\|\mu+|\nabla u|\|_{\underline L^p(B_{2r}(x))}+\|g_\delta\|_{L^{p^\prime}(B_{2r})}^\frac 1 {p-1}+\|\nabla u\|_{L^p(B_1)}.
\end{align}
The above inequality and the (yet to be proven) estimate
\begin{equation}\label{est:pf:largescalelp:nablalp}
\|\mu+|\nabla u|\|_{\underline L^p(B_{2r}(x))}\lesssim\|f\|_{\underline L^p(B_{4r}(x))}+\|g_\delta\|_{\underline L^{p'}(B_{4r}(x))}^\frac1{p-1}
\end{equation}
imply
\begin{equation*}
\|f_{x,r}\|_{\underline L^{2q}(B_r(x))}\lesssim c_q \|f\|_{\underline L^p(B_{4r}(x))}+\|g_\delta\|_{L^{p^\prime}(B_{4r})}^\frac 1 {p-1}.
\end{equation*}
Hence, \eqref{T:largelp:pf:claim1} is proven in this case. Next, we show \eqref{est:pf:largescalelp:nablalp}. Suppose first that $r\in[\frac{1}{24},\frac{1}{4}]$. In this case, we trivially have
$$
\|\mu+|\nabla u|\|_{\underline L^p(B_{2r}(x))}\lesssim \|\mu+|\nabla u|\|_{\underline L^p(\R^{n}))}\stackrel{\eqref{pf:Th:largescaleCZ:usctend}}{\lesssim} \mu+\|\nabla u\|_{\underline L^p(B_1)}\stackrel{\eqref{pf:largescaleCZ:deff}}{\leq}  \|f\|_{\underline L^p(B_{4r}(x))}
$$ 
and thus \eqref{est:pf:largescalelp:nablalp} is proven in this case. Next, we consider $r\in(0,\frac{1}{24})$.  Caccioppoli inequality \eqref{eq:caccioppoli} yields
\begin{align*}
\|\mu+|\nabla u|\|_{\underline L^p(B_{2r}(x))}\lesssim& \mu+(4r)^{-1}\|u-(u)_{B_{4r}(x)}\|_{\underline L^p(B_{4r}(x))}+\|F\|_{\underline L^{p'}(B_{4r}(x))}^\frac1{p-1}\,.
\end{align*}
The claimed estimate \eqref{est:pf:largescalelp:nablalp} follows with help of \eqref{est:strangepointwisetomaximal} and the fact that for $r\leq \frac1{24}$ it holds
$$
r^{-1}\|u-(u)_{B_{4r}(x)}\|_{\underline L^p(B_{4r}(x))}\lesssim \inf_{y\in B_{4r}(x)}r^{-1}\|u-(u)_{B_{8r}(y)}\|_{\underline L^p(B_{8r}(y))}^p\lesssim \inf_{y\in B_{4r}(x)}f(y)\leq\|f\|_{\underline L^p(B_{4r}(x))}
$$
%
%
It is left to show \eqref{T:largelp:pf:claim2}. Combining the triangle inequality and Sobolev-Poincar\'{e} inequality with $p_\ast = \frac{np}{n+p}$ for $p<n$ and $p_\ast=\frac {p+1} 2$ if $n\geq p$, we obtain
\begin{align*}
|f(z)-f_{x,r}(z)|^{p_\ast}\leq& \max_{r'\in [N\e,\frac 1 4]}\frac 1 {(r')^{p_\ast}}\|v_r-u-( v_r-u)_{B_{r'}(z)}\|_{\underline L^p(B_{r'}(z))}^{p_\ast}\\
\lesssim& \max_{r'\in [N\e,\frac 1 4]}\| \nabla (v_r-u)\|_{\underline L^{p_\ast}(B_{r'}(z))}^{p_\ast}\\
\leq& M(|\nabla (v_r-u)|^{p_\ast})(z)
\end{align*}
where $M$ denotes the maximal function. Hence, by the strong maximal function estimate and the fact that ${v_r-u\in W_0^{1,p}(B_{2r}(x))}$, we obtain
\begin{eqnarray*}
\|f_{x,r}-f\|_{\underline L^p(B_r(x))}&\leq&|B_r|^{-\frac1p}\| |f_{x,r}-f|^{p_\ast}\|_{L^\frac p {p_\ast}(B_r(x))}^\frac 1 {p_\ast}\\
&\lesssim& |B_r|^{-\frac 1 p} \|M(|\nabla (v_r-u)|^{p_\ast})\|_{L^\frac p {p_\ast}(B_r(x))}^\frac 1 {p_\ast}\\
&\lesssim&  |B_r|^{-\frac 1 p} \| |\nabla (v_r-u)|^{p_\ast} \|_{L^\frac p {p_\ast}(\R^n)}^\frac 1 {p_\ast}\\
&\lesssim& \|\nabla (v_r-u)\|_{\underline L^p(B_{2r}(x))}.
\end{eqnarray*}
Combining the above estimate with \eqref{das:2}, we obtain
\begin{eqnarray*}
\|f_{x,r}-f\|_{\underline L^p(B_r(x))}&\lesssim& \delta\|\mu+|\nabla u|\|_{\underline L^p(B_{2r}(x))}+\|g_\delta\|_{\underline L^{p^\prime}(B_{2r}(x))}^\frac 1 {p-1}\\
&\stackrel{\eqref{est:pf:largescalelp:nablalp}}\lesssim& \delta \|f\|_{\underline L^p(B_{4r}(x)}+\|g_\delta\|_{\underline L^{p^\prime}(B_{4r}(x))}^\frac 1 {p-1},
\end{eqnarray*}
and thus \eqref{T:largelp:pf:claim2} is proven.

\step 4 Conclusion. By an elementary covering argument, we can upgrade the conclusion of the Steps~1--3, to the following: There exist $C=C(\Lambda,n,p)>0$ and $C_q=C(\lambda,n,p,q)>0$ such that for every $x\in B_\frac12$ and $r\in(0,\frac14]$ there exists $f_{x,r}$ satisfying
\begin{align}
\|f_{x,r}\|_{\underline L^{2q}(B_r(x))}\leq C_q\|f\|_{\underline L^p(B_{2r}(x))}+C\|g_\delta \|_{\underline L^{p'}(B_{2r}(x))}^\frac1{p-1}\label{T:largelp:pf:claim1a}\\
\|f-f_{x,r}\|_{\underline L^p(B_r(x))}\leq C\delta \|f\|_{\underline L^p(B_{2r}(x))}+C\|g_\delta \|_{\underline L^{p'}(B_{4r}(x))}^\frac1{p-1}\label{T:largelp:pf:claim2a}
\end{align}
Next, we can apply Lemma~\ref{lem:CZ} with $q_0=2q$, $A=C_q$ and obtain for $\delta=C^{-1}\delta_0(n,p,q,2q,C_q)$ that
\begin{align}\label{lem:CZ:claim:pf}
\|f\|_{L^q(B_{1/2})}\leq C \left(\|f\|_{L^p(B_1)}+\|g_\delta\|_{L^\frac q {p-1}(B_1)}^\frac 1 {p-1}\right).
\end{align}
We now show that \eqref{lem:CZ:claim:pf} implies \eqref{L1:lp} for $B=B_1$. By Caccioppoli inequality \eqref{eq:caccioppoli}, we have
$$
\|\nabla u\|_{\underline L^p(B_{r/2}(x))}\lesssim \mu+r^{-1}\|u-(u)_{B_r(x)}\|_{\underline L^p(B_r(x))}+\|F\|_{\underline L^{p'}(B_r(x))}^\frac1{p-1}
$$
from which we deduce the pointwise bound
\begin{equation}\label{pf:pointwise:nablauN}
N^{-n}M_{\e}(\mathds 1_{B_1} |\nabla u|^p)(x)\leq M_{N\e}(\mathds 1_{B_1} |\nabla u|^p)(x)\lesssim f(x)^p+M_{\e}(\mathds 1_{B_1} |F|^{p'})(x)\quad\mbox{for a.e.\ $x\in B_\frac12$}.
\end{equation}
Next, we show that $\|f\|_{L^p(B_1)}$ can be estimated in terms of $\|\nabla u\|_{L^p(B_1)}$. For this, we argue as in the estimate for $\|f_{x,r}-f\|_{\underline L^p(B_r(x))}$ in Step~3. Set $p_\ast = \frac{np}{n+p}$ for $p<n$ and $p_\ast=\frac {p+1} 2$ if $n\geq p$. The Poincar\'e-Sobolev inequality implies
\begin{align*}
f(z)^{p_\ast}\lesssim&\mu^{p_\ast}+ \max_{r'\in [N\e,\frac 1 4]}\frac 1 {(r')^{p_\ast}}\|u-(u)_{B_{r'}(z)}\|_{\underline L^p(B_{r'}(z))}^{p_\ast}+\|\nabla u\|_{\underline L^p(B_1)}^{p_\ast}\\
\lesssim&\mu^{p_\ast}+ M(|\nabla u|^{p_\ast})(z)+\|\nabla u\|_{\underline L^p(B_1)}^{p_\ast}
\end{align*}
and thus
\begin{eqnarray}\label{pf:largeCZ:finalf}
\|f\|_{\underline L^p(B_1)}&=&\| |f|^{p_\ast}\|_{\underline L^\frac p {p_\ast}(B_1)}^\frac 1 {p_\ast}\lesssim\mu+\|M(|\nabla u|^{p_\ast})\|_{\underline L^\frac p {p_\ast}(B_1)}^\frac 1 {p_\ast}+\|\nabla u\|_{\underline L^p(B_1)}\notag\\
&\lesssim& \mu+\|\nabla u\|_{L^p(\R^n)}\stackrel{\eqref{pf:Th:largescaleCZ:usctend}}\lesssim\mu+\|\nabla u\|_{\underline L^p(B_1)}.
\end{eqnarray}
Combining \eqref{lem:CZ:claim:pf} with \eqref{pf:pointwise:nablauN}, \eqref{pf:largeCZ:finalf} and the fact that $N=N(\Lambda,n,p,q)\in\mathbb N$, we obtain the desired estimate \eqref{L1:lp} for $B=B_1$.

\step 5 Scaling. We briefly describe how the claim for $B=B_1$ implies the statement for $B=B_R$. Suppose that $u\in W^{1,p}(B_R)$ solves \eqref{eq:T1:Lp} on $B=B_R$. By scaling, we deduce that $u_R\in W^{1,p}(B_1)$ given by $u_R(x)=R^{-1}u(Rx)$ satisfies
\begin{equation}
\nabla \cdot\bfa(\tfrac{Rx}{\e},\nabla u_R)=\nabla \cdot F_R\quad\mbox{in $B_1$}
\end{equation}
with $F_R(x):=F(Rx)$. Hence, we can apply the result for $B=B_1$ with $\e$ replaced by $\e/R$ and obtain
\begin{equation*}
\|M_\frac\e{R}(\mathds 1_{B_1}|\nabla u_R|^p)\|_{\underline L^\frac{q}p(B_\frac12)}\leq C_q\biggl(\|\mu+|\nabla u_R|\|_{\underline L^p(B_1)}^p+\|M_\frac\e R(\mathds 1_{B}|F_R|^{p'})\|_{\underline L^{\frac{q}{p}}(B_1)}\biggr).
\end{equation*}
The desired estimate follows by scaling and the following identity for $h\in L^1(B_R)$
\begin{align*}
M_\frac\e{R}(\mathds 1_{B_1}h(R\cdot))(z)=\sup_{r\geq\frac\e R}\fint_{B_r(z)}\mathds 1_{B_1}h(Ry)\dy=\sup_{r\geq\frac\e R}\fint_{B_{Rr}(Rz)}\mathds 1_{B_R} h\dx=M_\e(\mathds 1_{B_R}h)(Rz).
\end{align*}
\end{proof}

We finish this section with the proof of Lemma~\ref{Preg:1:lp}.
\begin{proof}[Proof of Lemma~\ref{Preg:1:lp}]

We prove the statement for $R=2$. The general statement follows by a simple scaling argument. Throughout the proof we write $\lesssim$ if $\leq$ holds up to a positive multiplicative constant depending only on $\Lambda,n$ and $p$.

\step 1 Harmonic approximation. Let $\overline u\in u+W^{1,p}_0(B_{\frac 3 2})$ be the unique function satisfying
\begin{equation}\label{eq:uhomlp}
\nabla \cdot \overline \bfa(\nabla \overline u)=0\qquad \text{ in } B_{\frac 3 2}.
\end{equation}
Next, we gather the needed regularity properties for $\overline u$ which are direct consequences of Proposition~\ref{prop:regularity}. There exists $c=c(\Lambda,n,p)\in[1,\infty)$ and $m=m(\Lambda,n,p)>1$ such that
\begin{align}
\|\mu+|\nabla \overline u|\|_{L^p(B_\frac32)}\leq& c\|\mu+|\nabla u|\|_{L^p(B_2)}\label{est:eneubarzeta:lp}\\
\forall r\in (0,\tfrac32)\qquad[\nabla \overline u]_{\min\left\{1,\frac 1 {p-1}\right\},p,B_r}\leq& c(\tfrac32-r)^{-\min\left\{1,\frac 1 {p-1}\right\}}\|\mu+|\nabla u|\|_{L^p(B_2)}\label{est:besovubar:lp}\\
\|\mu+|\nabla \overline u|\|_{L^{mp}(B_\frac32)}\leq& c\|\mu+|\nabla u|\|_{L^p(B_2)}\label{est:nonlinearmeyerglob:lp1}.
\end{align}
Clearly, \eqref{est:eneubarzeta:lp} and \eqref{est:besovubar:lp} follow directly from Proposition~\ref{prop:regularity} (see estimates \eqref{L:energyestimate:eq1} and \eqref{est:reghom:bes}). For \eqref{est:nonlinearmeyerglob:lp1}, we first use interior higher gradient integrability, see \eqref{est:nonlinearmeyerloc:lp}, of $u$, that is there exists $m=m(\Lambda,n,p)>1$ such that
$$
\|\nabla u\|_{\underline L^{mp}(B_{\frac32}))}\lesssim\|\mu+|\nabla u\|_{\underline L^p(B_2)} 
$$
and then we use the global higher gradient integrability (see \eqref{est:nonlinearmeyerglob:lp}) for $\overline u$ to deduce \eqref{est:nonlinearmeyerglob:lp1}.

\step 2 Definition of two-scale expansion and equation for residuum.

Fix $\rho>0$, $\ell \in\mathbb N$ with $\ell\geq2$ and define
\begin{equation}\label{def:mathcalQ}
\mathcal Q:=\{Q_{\ell \e}(z)\,|\,z\in \e\ell\Z^d,\,Q_{\ell \e}(z)\cap B_{\frac32-2\rho}\neq\emptyset\},\quad\mbox{with}\quad Q_{r}(z):=z+(-\tfrac{r}2,\tfrac{r}2)^d.
\end{equation}
For every $Q\in\mathcal Q$, we consider a cut-off function
\begin{equation}\label{def:etaQ}
\eta_Q\in C_c^\infty(Q),\quad \eta_Q=1\quad\mbox{in $Q_{(-\e)}:=\{y\in Q\,|\,{\rm dist}(y,Q^c)>\e\}$},\quad |\nabla \eta_Q|\leq \frac8\e 
\end{equation}
and set
\begin{equation}\label{def:overlineu2s}
\overline u^{2s}(x):=\overline u_{\e,\ell,\rho}^{2s}(x):=\overline u(x)+\sum_{Q\in\mathcal Q}\eta_Q\phi_{\xi_Q,\e}
\end{equation}
where for every $Q\in\mathcal Q$, 
\begin{equation}\label{def:xik:lp}
\xi_Q:=(\nabla \overline u)_Q \quad\text{ and }\quad \phi_{\xi_Q,\e}=\e \phi_{\xi_Q}(\cdot/\e).
\end{equation}
We claim that
\begin{equation}\label{eq:error}
-\nabla \cdot (\bfa( \tfrac{x}\e,\nabla \overline u^{2s})-\bfa( \tfrac{x}\e,\nabla u))=\nabla \cdot \sum_{i=1}^3R^{(i)},
\end{equation}
where
\begin{align*}
R^{(1)}:=&\overline \bfa(\nabla \overline u)-\sum_{Q\in\mathcal Q}\eta_Q\overline \bfa(\xi_Q),\\
R^{(2)}:=&-\sum_{Q\in\mathcal Q} \sigma_{\xi_{Q},\e}\nabla \eta_{Q},\\
R^{(3)}:=&-\biggl(\bfa(\tfrac{x}\e,\nabla \overline u^{2s}) - \sum_{Q\in\mathcal Q}\eta_Q\bfa(x,\xi_Q+\nabla \phi_{\xi_{Q},\e}))\biggr),
\end{align*}
where $\sigma_{\xi,\e}:=\e\sigma_\xi(\frac{\cdot}\e)$ and $\sigma_\xi$ for $\xi\in\R^n$ is introduced in Lemma~\ref{L:phi}. Equation \eqref{eq:error} follows from 
$$
\nabla \cdot\bfa(\tfrac{x}\e,\nabla u)=\nabla \cdot \overline\bfa(\nabla \overline u),
$$
 and the properties of $\sigma_\xi$ in the form
$$
\nabla \cdot \sum_{Q\in\mathcal Q}\eta_{Q}(\bfa(\tfrac{x}\e,\xi_{Q}+\nabla \phi_{\xi_{Q}})-\overline \bfa(\xi_{Q}))=\nabla \cdot \sum_{Q\in\mathcal Q} \eta_{Q}\nabla \cdot \sigma_{\xi_{Q},\e}=\nabla \cdot \sum_{Q\in\mathcal Q}\sigma_{\xi_{Q},\e}\nabla \eta_{Q}
$$
for all $Q\in\mathcal Q$, where in last equality the skew-symmetry \eqref{prop:sigma2} of $\sigma_\xi$ is  used.

\step 3 Control by continuity. We set
\begin{equation}
z:=z_{\e,\ell,\rho}:=u-\overline u^{2s}.
\end{equation}
We claim that for $0<\e<\rho\leq \frac1{10}$, $\ell\in\mathbb N$ with $\sqrt{n}\ell \e<\rho$ and every $\tau\in(0,1]$ it holds
\begin{align}\label{P:1step:claim:s2:lp}
&\int_{B_\frac32}|\nabla z|^p\dx\notag\\
\lesssim&\biggl(\tau+(1+\tau^{\frac{p-2}{p-1}})(\rho^{1-\frac1m}+((\ell\e\rho^{-1})^\frac{p}{p-1}+\ell^{-1})^\frac1{p-1}+(\ell\e\rho^{-1})^{p}+\ell^{-1})\biggr)\int_{B_2}(\mu+|\nabla u|)^p\dx.
\end{align}
where $m=m(\Lambda,n,p)>1$ is as in \eqref{est:nonlinearmeyerglob:lp1}. In order to prove \eqref{P:1step:claim:s2:lp}, we use the comparison estimate \eqref{L:energyestimate:eq2} of Proposition~\ref{prop:regularity} with $u=\overline u^{2s}$, $\bfa=\overline\bfa$ (and thus $\delta=0$), $F=\sum_{j=1}^3|R^{(j)}|$ and $w=u$ (and $\tau$ replaced by $2^{-1}\tau^\frac1p$) in the form
%
\begin{align*}
\fint_{B_\frac 3 2}|\nabla z|^p\leq \tau\fint_{B_\frac32}(\mu+|\nabla \overline u^{2s}|)^p\dx+c(1+\tau^\frac{p-2}{p-1})\sum_{i=1}^3\fint_{B_\frac32}|R^{(i)}|^{p^\prime}\dx,
\end{align*}
where $c=c(\Lambda,n,p)\in[1,\infty)$. The claimed estimate follows from
\begin{equation}\label{est:energyestu2s:lp}
\int_{B_\frac32}|\nabla \overline u^{2s}|^p\dx\lesssim \int_{B_2}(\mu+|\nabla u|)^p\dx
\end{equation} 
(which is a direct consequence of \eqref{est:eneubarzeta:lp} combined with \eqref{eq:decomp3} proven below) and suitable estimates \eqref{claim:est:R1}, \eqref{claim:est:R2} and \eqref{claim:est:R3} on $\int_{B_\frac32}|R^{(i)}|^{p^\prime}$ which we provide below.

\substep{3.1} The term $R^{(1)}$: We claim that 
\begin{align}\label{claim:est:R1}
\int_{B_\frac32}|R^{(1)}|^{p^\prime}\dx\lesssim&\biggl(\rho^{1-\frac1m}+(\ell\e\rho^{-1})^\frac{p}{(p-1)^2}+(\ell\e\rho^{-1})^{p}+\ell^{-1}\biggr)\int_{B_2}(\mu+|\nabla u|)^p\dx.
\end{align}
Indeed,  the decomposition
\begin{align}\label{eq:decomp}
|R^{(1)}|\leq |\overline\bfa(\nabla \overline u)-\sum_{Q\in\mathcal Q}\mathds 1_Q\overline \bfa(\nabla \overline u)|+\sum_{Q\in\mathcal Q}\mathds 1_Q|\overline \bfa(\xi_Q)-\overline\bfa (\nabla \overline u)|+\sum_{Q\in\mathcal Q}(\mathds 1_{Q}-\eta_Q)|\overline\bfa (\xi_Q)|
\end{align}
and the continuity of $\overline \bfa$ (see \eqref{eq:continuityHomogenised1}) yields
\begin{align}\label{est:pw:R1:lp}
|R^{(1)}_{\delta,\rho}|\lesssim&  (\mathds 1_{B_{\frac32}}-\mathds 1_{B_{\frac32-2\rho}})(\mu+|\nabla \overline u|)^{p-1}+ \mathds 1_{\{p\geq2\}}\sum_{Q\in\mathcal Q}\mathds 1_Q(\mu+|\xi_Q|+|\nabla \overline u|)^{p-2}|\nabla \overline u-\xi_Q|\notag\\
&+\mathds 1_{\{p<2\}}\sum_{Q\in\mathcal Q}\mathds 1_Q|\nabla \overline u-\xi_Q|^{p-1}+\sum_{Q\in\mathcal Q} (\mathds 1_Q-\eta_Q)(\mu+|\xi_Q|)^{p-1},
\end{align}
and we use the definition of $\mathcal Q$, see \eqref{def:mathcalQ}, and the assumption $\sqrt{n}\ell \e<\rho$ in the form
\begin{equation}\label{est:sumiq}
\mathds 1_{B_{\frac32-2\rho}}\leq \sum_{Q\in\mathcal Q}\mathds 1_Q\leq\mathds 1_{B_{\frac32-\rho}}.
\end{equation}
The higher integrability \eqref{est:nonlinearmeyerglob:lp1} and H\"older inequality yields
\begin{align}\label{est:boundarylayer:3rho:lp}
&\int_{B_\frac32\setminus B_{\frac32-2\rho}}(\mu+|\nabla \overline u|)^{(p-1)p^\prime}\dx= \int_{B_\frac32\setminus B_{\frac32-2\rho}}(\mu+|\nabla \overline u|)^{p}\dx\notag\\
\leq& |B_\frac32\setminus B_{(1-2\rho)}|^{1-\frac1{ m}}\biggl(\int_{B_\frac32}(\mu+|\nabla \overline u|)^{mp}\dx\biggr)^\frac{1}{ m}\stackrel{\eqref{est:nonlinearmeyerglob:lp1}}\lesssim \rho^{1-\frac1m}\int_{B_2}(\mu+|\nabla u|)^p\dx.
\end{align}
To estimate the last term on the right-hand side in \eqref{est:pw:R1:lp}, we use
\begin{equation}\label{measure:1qetaq}
\forall Q\in\mathcal Q:\quad|\{{\rm supp}(\mathds 1_Q-\eta_Q)\}|\lesssim |Q\setminus Q_{(-\e)}|\lesssim \ell^{-1}|Q|
\end{equation}
and Jensen inequality with the definition of $\xi_Q$, see \eqref{def:xik:lp} in the form
\begin{equation}\label{est:sumxiQtou:lp}
\int_{B_\frac32}\sum_{Q\in\mathcal Q}\mathds 1_Q|\xi_Q|^p\dx=\sum_{Q\in\mathcal Q}|Q||\xi_Q|^p\leq \sum_{Q\in\mathcal Q}\int_Q|\nabla \overline u|^p\dx\stackrel{\eqref{est:sumiq}}\leq \int_{B_{\frac32}}|\nabla \overline u|^p\dx\stackrel{\eqref{est:eneubarzeta:lp}}\lesssim \int_{B_2}(\mu+|\nabla u|)^p\dx.
\end{equation}
Combining \eqref{measure:1qetaq} and \eqref{est:sumxiQtou:lp}, we obtain
\begin{align}\label{est:internalboundarylayer:lp}
\sum_{Q\in\mathcal Q}\int_{B_\frac32}|(\mathds 1_Q-\eta_Q)(\mu+|\xi_Q|)^{p-1}|^{p'}\dx\lesssim& \sum_{Q\in\mathcal Q}\int_{B_\frac32}(\mathds 1_Q-\eta_Q)^{p^\prime}(\mu+|\xi_Q|)^p\dx\notag\\
\lesssim&\ell^{-1}\sum_{Q\in\mathcal Q}|Q|(\mu+|\xi_Q|)^p\stackrel{\eqref{est:sumxiQtou:lp}}\lesssim\ell^{-1}\int_{B_2}(\mu+|\nabla u|)^p\dx.
\end{align}
To estimate the second and third term coming from \eqref{est:pw:R1:lp}, we use Lemma~\ref{L:BesovPoincare} and \eqref{est:besovubar:lp} in the form
\begin{align}\label{ast:xikovu:lp}
\int_{B_s}\sum_{Q\in\mathcal Q}\mathds 1_Q| \nabla \overline u-\xi_Q|^p\dx\lesssim&  (\ell\e)^{\min\{p,\frac p {p-1}\}}\|\nabla \overline u\|_{\dot B_\infty^{{\min\{1,\frac1{p-1}\}},p}(B_{\frac32-\rho})}^p\notag\\
\lesssim& (\ell\e\rho^{-1})^{\min\{p,\frac p{p-1}\}}\int_{B_2}(\mu+|\nabla u|)^p\dx.
\end{align}
In the case $p\geq2$, we obtain 
\begin{eqnarray}\label{ast:xikovufullR1:lp}
& &\int_{B_\frac32}\sum_{Q\in\mathcal Q} \mathds 1_Q \left((\mu+|\xi_Q|+|\nabla \overline u|)^{p-2}|\nabla \overline u-\xi_Q|\right)^{p^\prime}\dx\notag\\
&\leq&\biggl(\int_{B_\frac32}\sum_{Q\in\mathcal Q}\mathds 1_Q(\mu+|\xi_Q|+|\nabla \overline u|)^p\dx\biggr)^\frac{(p-2)}{p-1}\biggl(\int_{B_\frac32}\sum_{Q\in\mathcal Q}\mathds 1_Q|\nabla \overline u-\xi_Q|^p\dx\biggr)^\frac{1}{p-1}\notag\\
&\stackrel{\eqref{est:eneubarzeta:lp},\eqref{est:sumxiQtou:lp},\eqref{ast:xikovu:lp}}\lesssim&(\ell\e\rho^{-1})^\frac {p}{(p-1)^2} \int_{B_2} (\mu+|\nabla u|)^p\dx\qquad (\mbox{for }p\geq2).
\end{eqnarray}
Similarly, we obtain for $p\in(1,2)$
\begin{align}\label{ast:xikovufullR1:lp:sub}
\int_{B_\frac32}\sum_{Q\in\mathcal Q} \mathds 1_Q \left(|\nabla \overline u-\xi_Q|^{p-1}\right)^{p^\prime}\dx=&\int_{B_\frac32}\sum_{Q\in\mathcal Q} \mathds 1_Q |\nabla \overline u-\xi_Q|^p\dx\notag\\
\lesssim&(\ell\e\rho^{-1})^{p} \int_{B_2} (\mu+|\nabla u|)^p\dx\qquad (\mbox{for }1<p<2).
\end{align}
Combining the decomposition \eqref{est:pw:R1:lp} with the estimates \eqref{est:boundarylayer:3rho:lp}, \eqref{ast:xikovufullR1:lp}, \eqref{ast:xikovufullR1:lp:sub} and \eqref{est:internalboundarylayer:lp}, we obtain the claimed inequality \eqref{claim:est:R1}.

\substep{3.2} The term $R^{(2)}$: Using $|\nabla \eta|\leq \frac 8 \e(\mathds 1_Q-\mathds 1_{(Q)_{-\e}})$, Lemma \ref{lem:easyperiodic}, \eqref{measure:1qetaq}, \eqref{eq:correctorNatural} and \eqref{est:sumxiQtou:lp}, we estimate
\begin{align}\label{claim:est:R2}
\int_{B_\frac32}|R^{(2)}|^{p^\prime}\dx\lesssim&\sum_{Q\in\mathcal Q} \int_Q|\sigma_{\xi_Q,\e}\nabla \eta_Q|^{p^\prime}\dx\lesssim \ell^{-1}\sum_{Q\in \mathcal Q}|Q|(\mu+|\xi_Q|)^p\notag\\
 \lesssim& \ell^{-1}\int_{B_2}(\mu+|\nabla u|)^p\dx.
\end{align}
\substep{3.3} The term $R^{(3)}$: We claim
\begin{equation}\label{claim:est:R3}
\int_{B_\frac32}|R^{(3)}|^{p^\prime}\dx\lesssim\biggl(\rho^{1-\frac1m}+(\ell\e\rho^{-1})^\frac{p}{(p-1)^2}+\ell^{-\frac1{p-1}}+(\ell\e\rho^{-1})^{p}+\ell^{-1}\biggr)\int_{B_2}(\mu+|\nabla u|)^p\dx.
\end{equation}
We note that
\begin{align}\label{eq:decomp2}
\sum_{Q\in\mathcal Q}\int_{B_\frac32}(\mathds 1_Q-\eta_Q)|\nabla \phi_{\xi_Q,\e}|^p\dx\lesssim& \ell^{-1}\sum_{Q\in\mathcal Q}|Q||\xi_Q|^p\lesssim \ell^{-1}\int_{B_2}(\mu+|\nabla u|)^p\dx,
\end{align}
and by similar computations
\begin{align}\label{eq:decomp3}
\sum_{Q\in\mathcal Q}\int_{B_\frac32}\mathds 1_Q |\nabla \phi_{\xi_Q,\e}|^p\dx+\ell \sum_{Q\in\mathcal Q}\int_{B_\frac32}(\mathds 1_Q-\mathds 1_{(Q)_{-\e}})|\phi_{\xi_Q,\e}\e^{-1}|^p\dx\lesssim\int_{B_2}(\mu+|\nabla u|)^p\dx. 
\end{align}

 The continuity of $\bfa$  implies
\begin{align*}
|R^{(3)}|\lesssim& \mathds 1_{B_\frac32\setminus B_{\frac32-2\rho}}(\mu+|\nabla\overline u|)^{p-1}+ \mathds 1_{\{p\geq2\}} \sum_{Q\in\mathcal Q}\mathds 1_Q(\mu+|\nabla \overline u^{2s}|+|\xi_Q+\nabla \phi_{\xi_{Q}}|)^{p-2}|\nabla \overline u^{2s}-(\xi_Q+\nabla \phi_{\xi_{Q}})|\\
&+ \mathds 1_{\{p<2\}} \sum_{Q\in\mathcal Q}\mathds 1_Q|\nabla \overline u^{2s}-(\xi_Q+\nabla \phi_{\xi_{Q}})|^{p-1}+\sum_{Q\in\mathcal Q}(\mathds 1_Q-\eta_Q)(\mu+|\xi_Q+\nabla \phi_{\xi_{Q}}|)^{p-1}\\
&=: I + \mathds 1_{\{p\geq2\}}II +\mathds 1_{\{p<2\}}II'+ III.
\end{align*}
The terms $I$ and $III$ can be estimated similarly as in Substep 3.1 (using in addition \eqref{eq:decomp2}), to the effect that
\begin{equation}
\int_{B_\frac32}|I+III|^{p^\prime}\dx\lesssim(\rho^{1-\frac1m}+\ell^{-1})\int_{B_2}(\mu+|\nabla u|)^p\dx.
\end{equation}
Next, we focus on the remaining terms $II$ and $II'$. In the case $p\geq2$, we find by H\"older inequality
\begin{align*}
\int_{B_\frac32} |II|^{p^\prime}\dx\lesssim&\biggl(\int_{B_\frac32}\mu^p+|\nabla \overline u^{2s}|^p+\sum_{Q\in\mathcal Q}\mathds 1_Q|\xi_Q+\nabla \phi_{\xi_Q}|^p\dx\biggr)^\frac{p-2}{p-1}\\
&\quad\times\biggl(\int_{B_\frac32}\sum_{Q\in\mathcal Q}\mathds 1_Q|\nabla \overline u^{2s}-(\xi_Q+\nabla \phi_{\xi_Q})|^p\dx\biggr)^\frac 1{p-1}\\
\lesssim&\biggl(\int_{B_2}(\mu+|\nabla u|)^p\dx\biggr)^\frac{p-2}{p-1}\biggl(\int_{B_\frac32}\sum_{Q\in\mathcal Q}\mathds 1_Q|\nabla \overline u^{2s}-(\xi_Q+\nabla \phi_{\xi_Q})|^p\dx\biggr)^\frac 1{p-1},
\end{align*}
where we use \eqref{est:energyestu2s:lp}, \eqref{est:sumxiQtou:lp} and \eqref{eq:decomp3} in the second estimate. Instead, we have in the case $p\in(1,2)$
\begin{align*}
\int_{B_\frac32} |II'|^{p^\prime}\dx=&\int_{B_\frac32}\sum_{Q\in\mathcal Q}\mathds 1_Q|\nabla \overline u^{2s}-(\xi_Q+\nabla \phi_{\xi_Q})|^p\dx.
\end{align*}
Hence, it is left to show
\begin{align}\label{est:bigclaimR3:lp}
&\int_{B_\frac32}\sum_{Q\in\mathcal Q} \mathds 1_Q|\nabla \overline u^{2s}-(\xi_Q+\nabla \phi_{\xi_{Q},\e})|^p\dx\lesssim((\ell \e \rho^{-1})^{\min\{p,\frac p {p-1}\}} +\ell^{-1})\int_{B_2}(\mu+|\nabla u|)^p\dx\end{align}
For this we use the identity 
\begin{align}\label{eq:u2s-xik}
(\nabla \overline u^{2s}-(\xi_Q+\nabla \phi_{\xi_{Q},\e}))=
(\nabla \overline u-\xi_Q)+(\eta_Q-1)\nabla \phi_{\xi_Q,\e}+\phi_{\xi_{Q},\e}\nabla\eta_{Q}\qquad\mbox{on $Q\in\mathcal Q$}.
\end{align}
Combining \eqref{ast:xikovu:lp}, \eqref{eq:decomp2} and \eqref{eq:decomp3}
we obtain
\begin{align*}
&\int_{B_\frac32}\sum_{Q\in\mathcal Q} \mathds 1_Q|\nabla \overline u^{2s}-(\xi_Q+\nabla \phi_{\xi_{Q},\e})|^p\dx\notag\\
\stackrel{\eqref{eq:u2s-xik}}\lesssim&\sum_{Q\in\mathcal Q}\int_{B_\frac32}\mathds 1_Q|\nabla \overline u-\xi_Q|^p+(\mathds 1_Q-\eta_Q)|\nabla \phi_{\xi_Q,\e}|^p+\mathds 1_Q |\phi_{\xi_Q,\e}\nabla \eta_Q|^p\dx\notag\\
\stackrel{\eqref{ast:xikovu:lp}}\lesssim&((\ell \e \rho^{-1})^{\min\{p,\frac p {p-1}\}}+\ell^{-1})\int_{B_2}(\mu+|\nabla u|)^p\dx,
\end{align*}
and thus \eqref{est:bigclaimR3:lp} is proven.

\step 4 Conclusion. Steps~1--3 and a suitable choice for $\rho,\tau\in(0,1]$ and $\ell$ depending only on $\Lambda,n$ and $p$ imply that there exists $\beta=\beta(\Lambda,n,p)>0$ and $c=c(\Lambda,n,p)\in[1,\infty)$ such that for every solution of \eqref{eq:P1} with $R=2$ and $0<\e\leq \e_0=\e_0(\Lambda,n,p)\leq1$, we have
\begin{equation}\label{P:reg:lp:final:error}
\|\nabla u-\nabla \overline u^{2s}\|_{L^p(B_\frac32)}\leq c\e^\beta \|\mu+|\nabla u|\|_{L^p(B_2)},
\end{equation}
and thus \eqref{P1:est:xi-zeta} is satisfied for $R=2$. By Theorem~\ref{thm:Lq}, we find for every $q\in[p,\infty)$ a constant ${c_q=c_q(\Lambda,n,p,q)\in[1,\infty)}$ such that
\begin{equation}\label{est:ubarlq}
\|\mu+|\nabla \overline u|\|_{L^q(B_\frac{11}8)}\leq c_q\|\mu+|\nabla \overline u|\|_{L^p(B_\frac{3}2)}\lesssim c_q \|\mu+|\nabla u|\|_{L^p(B_2)}.
\end{equation}
We claim that for every $x\in B_1$ and every $\e\leq r\leq \frac14$ it holds
\begin{equation}\label{est:nablau2sBrvsmaxfunction}
\fint_{B_{r}(x)}|\nabla \overline u^{2s}|^p\dy\lesssim M(\mathds 1_{B_\frac{11}8}(\mu+|\nabla \overline u|)^p)(x),
\end{equation}
where as usual $M(f)$ denotes the maximal function of $f$. By triangle inequality, we have
$$
|\nabla \overline u^{2s}|\lesssim |\nabla \overline u|+\sum_{Q\in\mathcal Q}(|\phi_{\e,\xi_Q}\nabla \eta_Q|+\eta_Q|\nabla \phi_{\e,\xi_Q}|).
$$
If $r\geq \ell\e$, we obtain with a similar computations as in Step 3.2, using Lemma \ref{lem:easyperiodic}
\begin{align*}
\fint_{B_r(x)}|\sum_{Q\in\mathcal Q}(|\phi_{\e,\xi_Q}\nabla \eta_Q|+\eta_Q|\nabla \phi_{\e,\xi_Q}|)|^p\dy\lesssim&\frac1{|B_r|}\sum_{\genfrac{}{}{0pt}{2}{Q\in\mathcal Q}{Q\cap B_r(x)\neq\emptyset}}\int_Q(|\phi_{\e,\xi_Q}\nabla \eta_Q|+\eta_Q|\nabla \phi_{\e,\xi_Q}|)|^p\dy\\
\lesssim&\frac1{|B_r|}\sum_{\genfrac{}{}{0pt}{2}{Q\in\mathcal Q}{Q\cap B_r(x)\neq\emptyset}}|Q|(\mu+|\xi_Q|)^p\\
\lesssim&\fint_{B_{r+\sqrt{n}\e\ell}(x)}|\nabla \overline u|^p\dy\lesssim M(\mathds 1_{B_\frac{11}8}|\nabla \overline u|^p)(x).
\end{align*}
For $\e\leq r<\ell \e$, we can use that $B_r(x)\cap Q\neq\emptyset$ implies $B_r(x)\subset 2Q$ in the form
  \begin{align*}
\fint_{B_r(x)}|\sum_{Q\in\mathcal Q}(|\phi_{\e,\xi_Q}\nabla \eta_Q|+\eta_Q|\nabla \phi_{\e,\xi_Q}|)|^p\dy\lesssim&\max_{\genfrac{}{}{0pt}{2}{Q\in\mathcal Q}{ Q\cap B_r(x)\neq\emptyset}}(\mu+|\xi_Q|)^p\lesssim M(\mathds 1_{B_\frac{11}8}|\nabla \overline u|^p)(x).
\end{align*}

Using the strong maximal function estimate, we have
\begin{eqnarray}\label{P:reg:lp:final:lqest}
\fint_{B_1} \biggl(\max_{r\in [\e,\frac14]} \fint_{B_r(x)}|\nabla \overline u^{2s}|^p\dy\biggr)^\frac q p\dx
&\stackrel{\eqref{est:nablau2sBrvsmaxfunction}}\lesssim& \fint_{B_1} (M(\mathds 1_{B_\frac{11}8}(\mu+|\nabla \overline u|)^p)(x))^\frac{q}p\dx\notag\\
&\lesssim& \int_{\R^n} \mathds 1_{B_\frac{11}8}\left(\mu+|\nabla \overline u|\right)^q\dx\notag\\
&\stackrel{\eqref{est:ubarlq}}\lesssim& \left(\fint_{B_2} (\mu+|\nabla u|)^p\dx\right)^\frac q p.
\end{eqnarray}
%

Choosing $v=u^{2s}\in u+W_0^{1,p}(B_\frac32)$, extending $v$ by $u$ outside $B_\frac32$ and using \eqref{P:reg:lp:final:error} and \eqref{P:reg:lp:final:lqest}, we have proven the lemma for $R=2$ and $0<\e\leq \e_0(\Lambda,n,p)\in(0,1]$. The case $\e\in(\e_0,\frac12]$ is trivial as we can choose $v=u$ which trivially implies \eqref{P1:est:xi-zeta} and we have
\begin{align*}
\fint_{B_1} \biggl(\max_{r\in [\e,\frac14]} \fint_{B_r(x)}|\nabla u|^p\dy\biggr)^\frac q p\dx
\leq& (2\e_0)^{-\frac{q n}p}\biggl(\fint_{B_{2}}|\nabla u|^p\dy\biggr)^\frac{q}p.
\end{align*}

\end{proof}

A direct corollary of Theorem~\ref{L:nonlinear:largescalereg:lp} is the following
\begin{corollary}\label{cor:largeScaleHolder}
Consider the situation of Theorem \ref{L:nonlinear:largescalereg:lp}. For $1<p<q<\infty$, there exists $c=c(\Lambda,n,p,q)$ such that the following is true: Suppose that \eqref{eq:T1:Lp} holds with $B=B_R(x_0)$ and $F\in L^\frac q {p-1}(B)$. Then, 
\begin{align}\label{eq:largeScaleHolder}
\forall r\in[\e,R]:\quad \|\nabla u\|_{\underline L^p(B_r(x_0))}\leq c\left(\frac R r\right)^{\frac{n}{q}} \biggl(\|\mu+|\nabla u|\|_{\underline L^p(B_R(x_0))}+\|F\|_{\underline L^\frac q {p-1}(B(x_0,R))}^\frac 1 {p-1}\biggr).
\end{align}
\end{corollary}
\begin{proof}
Without loss of generality, we assume $x_0=0$ and $\e\leq r<\frac{R}8$ (estimate \eqref{eq:largeScaleHolder} trivially holds for $r\in[\frac{R}8,R]$). Since $B_r:=B_r(0)\subset B_{2r}(x)$ for all $|x|<r$, we have 
$$
\fint_{B_r}|\nabla u|^p\dy\leq 2^n\fint_{B_r}\fint_{B_{2r}(x)}|\nabla u|^p\dy\dx
$$
and thus in combination with Jensen inequality for all $q\in(p,\infty)$
\begin{equation}\label{pf:largescaleHolder:est1}
\fint_{B_r}|\nabla u|^p\dy\leq 2^n\biggl(\frac{R}{2r}\biggr)^\frac{np}q\biggl(\fint_{B_{R/2}}\biggl(\fint_{B_{2r}(x)}|\nabla u|^p\dy\biggr)^\frac{q}p\dx\biggr)^\frac{p}q.
\end{equation}
Appealing to Theorem~\ref{L:nonlinear:largescalereg:lp} and $\e\leq r\leq \frac{R}8$, we obtain
\begin{align}\label{pf:largescaleHolder:est2}
\biggl(\fint_{B_{R/2}}\biggl(\fint_{B_{2r}(x)}|\nabla u|^p\dy\biggr)^\frac{q}p\dx\biggr)^\frac{p}q\leq& \|M_\e(\mathds 1_{B_R}|\nabla u|^p)\|_{\underline L^\frac{q}p(B_{\frac{R}2})}\notag\\
\leq& C\biggl(\|\mu+|\nabla u|\|_{\underline L^p(B_R)}^p+\|M_\e(\mathds 1_{B}|F|^{p'})\|_{\underline L^{\frac{q}{p}}(B_{R})}\biggr).
\end{align}
Combining \eqref{pf:largescaleHolder:est1}, \eqref{pf:largescaleHolder:est2} together with the maximal function estimate 
$$\|M_\e(\mathds 1_{B_R}|F|^{p'})\|_{\underline L^{\frac{q}{p}}(B_{R})}\leq |B_R|^{-\frac{p}q}\|M((\mathds 1_{B_R}|F|^{p'}))\|_{L^{\frac{q}{p}}(\R^n)}\leq c(n,\tfrac{q}p)\|F\|_{\underline L^{\frac{q}{p-1}}(B_R)}^\frac{p}{p-1},
$$
we obtain the desired estimate \eqref{eq:largeScaleHolder}.
\end{proof}

\begin{remark}\label{rem:wang2019}
As mentioned in the introduction, \cite[Theorem 1.2]{Wang2019} with $\kappa=\frac n q$ is Corollary~\ref{cor:largeScaleHolder} in the case $F\equiv0$ supposing that $\bfa$ and $\overline \bfa$ satisfy Assumption~\ref{ass:standard} with $\mu=0$. As mentioned above, the stability of Assumption~\ref{ass:standard} under homogenization is in general questionable, see however \cite[Section 5]{Wang2019} where it is verified for some specific examples.
\end{remark}

\section{Large-scale Lipschitz estimates}\label{sec:largescaleLip}

The main goal in this section is to prove a \textit{large-scale Lipschitz} estimate, see Theorem~\ref{L:nonlinear:largescalereg}, from which Theorem~\ref{thm:Linfty:uniformeps} easily follows. The main ingredient for Theorem~\ref{L:nonlinear:largescalereg} are a regularity result for operators satisfying Assumption~\ref{ass:homonondeg}, see Proposition~\ref{lem:HighRegularity}, and a  two-scale homogenization estimate, see Proposition~\ref{Preg:1} below, which is very similar to Lemma~\ref{Preg:1:lp}.

\begin{theorem}\label{L:nonlinear:largescalereg}
Suppose $\bfa$ is periodic in the first variable. Assume $\bfa$ and $\overline \bfa$ satisfy Assumption~\ref{ass:homonondeg} for some $1<p<\infty$. Fix $M\in[1,\infty)$. There exists $C_{M}=C_M(\Lambda,M,n,p)\in[1,\infty)$ such that the following is true. Let $u\in W^{1,p}(B)$ with $B=B_R(x_0)\subset \R^n$ be such that
\begin{equation}\label{eq:T1}
\nabla\cdot \bfa(\tfrac{x}\e,\nabla u)=0\qquad\mbox{in $B$}
\end{equation}
and
\begin{equation}\label{eq:T1:b}
\fint_{B}|V_{p}(\nabla u)|^2\dx\leq M.
\end{equation}
Then, for all $r\in[\e,R]$
\begin{equation}\label{L1:lipschitz}
\fint_{B_r(x_0)}|V_{p}(\nabla u)|^2\dx \leq C_M  \fint_{B_R(x_0)}|V_{p}(\nabla u)|^2\dx.
\end{equation}
\end{theorem}

We postpone the proof of Theorem \ref{L:nonlinear:largescalereg} to Section~\ref{sec:prooflargescalelip} and instead show how to deduce Theorem~\ref{thm:Linfty:uniformeps} from it.

\begin{proof}[Proof of Theorem~\ref{thm:Linfty:uniformeps}]
 Note that due to Lemma \ref{lem:stability}, $\bar \bfa$ satisfies Assumption \ref{ass:homonondeg}. As Assumption \ref{ass:standard} with $\mu=1$ is stronger than Assumption \ref{ass:homonondeg}, $\bfa$ also satisfies Assumption \ref{ass:homonondeg}. Consequently, we may apply Theorem \ref{L:nonlinear:largescalereg}.

We distinguish between small and large balls starting with the larger ones.

\step 1 Suppose that $B=B_R(x_0)$ with $R\geq 2\e$. By Theorem \ref{L:nonlinear:largescalereg}, we find $C_M=C_M(\Lambda,M,n,p)>0$ such that for any $z\in \frac 1 2 B$,
\begin{align}
\fint_{B_\e(z)} |V_p(\nabla u)|^2 \leq C_M \fint_{B_{R/2}(z)} |V_p(\nabla u)|^2 \leq 2^d C_M \fint_{B} |V_p(\nabla u)|^2.
\end{align}

\step 2 In light of Step 1, for any $z\in \frac 1 2 B$, we may find a point $x(z)\in\{z,x_0\}$ and a radius $r(z)\leq 2 \e$ such that $z\in \frac 1 2 B_{r(z)}(x(z))$ and
\begin{align}\label{eq:smallBallAssumption}
\fint_{B_{r(z)}(x(z))} |V_p(\nabla u)|^2 \leq 2^d C_M \fint_B |V_p(\nabla u)|^2.
\end{align}
Indeed, if $R\geq 2 \e$, due to Step 1 we may choose $x(z)=z$ and $r(z)=\e$, while if $R<2\e$, we choose $x(z)=x_0$ and $r(z)=R$.

By rescaling, we have that
\begin{align*}
    \nabla \cdot \bfa(y,\nabla u_\e) = 0 \quad \text{ in } B_{r(z)/\e}(x(z)/\e),
\end{align*}
where $u_\e:=\e^{-1}u(\e x)$. As $r(z)/\e\leq 2$, we deduce from essentially standard estimates, see Proposition \ref{eq:LipschitzStandard} in Appendix \ref{sec:appendixB}, that
\begin{align*}
\|V_p(\nabla u_\e)\|_{L^\infty(B_{r(z)/(2\e)}(x(z))}\leq C \fint_{B_{r(z)/\e}(x(z)/\e)} |V_p(\nabla u_\e)|^2
\end{align*}
with $C=C(\gamma,n,\Lambda,p)$ and thus by scaling
\begin{align*}
\|V_p(\nabla u)\|_{L^\infty(B_{r(z)/2}(x(z)))}\leq C \fint_{B_{r(z)}(x(z))} |V_p(\nabla u)|^2,
\end{align*}
which in combination with \eqref{eq:smallBallAssumption} proves the statement.
\end{proof}

\subsection{Auxiliary results}\label{sec:auxiliary}
 In this section, we collect a regularity result for autonomous operators satisfying Assumption \ref{ass:homonondeg} and a two-scale homogenization estimate in the situation of Theorem~\ref{L:nonlinear:largescalereg}.

 We begin with the regularity results.
\begin{proposition}\label{lem:HighRegularity}
Suppose that $\bfa:\R^n\to\R^n$ satisfies Assumption \ref{ass:homonondeg} for some $1<p<\infty$, $\Lambda\in[1,\infty)$. There exist $c=c(\Lambda,n,p)>0$ and $\Gamma=\Gamma(n,p)>0$ such that the following is true: Suppose that $B=B(x,r)\Subset\Omega$, $w\in W^{1,p}(B)$ and $\zeta\in\R^n$ satisfy
\begin{equation}\label{lem:HighRegularity:eq}
\nabla \cdot  \bfa (\zeta+\nabla w)=0
\end{equation}
and
\begin{equation*}
|V_{p}(\zeta)|^2+\fint_{B}|V_{p}(\nabla w)|^2\dx\leq M\in[1,\infty).
\end{equation*}
Then
\begin{align}\label{reg:ahom:power}
r\|\nabla^2 w\|_{\underline L^2(\frac12 B)}+\|\nabla w\|_{L^\infty(\frac12B)}\leq& cM^\Gamma \|V_{p}(\nabla w)\|_{\underline L^2(B)}
\end{align}
Moreover, there exist $\alpha=\alpha(\Lambda,M,n,p)>0$ and $C_M=C_M(\Lambda,M,n,p)>0$ such that
\begin{equation}\label{reg:ahom:holder}
r^\alpha\|\nabla w\|_{C^{0,\alpha}(\frac12 B)}\leq C_M   \|V_{p}(\nabla w)\|_{\underline L^2(B)}.
\end{equation}
\end{proposition}

\begin{proof}
The crucial part of the proof is to obtain the $L^\infty$ bound on $\nabla w$ stated in \eqref{reg:ahom:power}. This can be done in a self-contained way by revisiting the Moser type iteration in the proof of Theorem~\ref{thm:Lq}, using the stronger monotonicity and continuity assumptions on $\bfa$ and the higher integrability of Theorem~\ref{thm:Lq}. However, we present here a shorter argument using known results for linear non-uniformly elliptic equations. Throughout the proof we write $\lesssim$ if $\leq$ holds up to a multiplicative positive constant depending only on $\Lambda,n$ and $p$.

\step 1 Suppose that $\bfa\in C^1(\R^n,\R^n)$ and $B=B_1(0)$. We claim that
\begin{equation}\label{est:L:reg:ahom:a2:linfty}
\|\nabla^2 w\|_{L^2(B_{1/8})}+\|\nabla w\|_{L^\infty(B_{1/4})}\lesssim M^\Gamma\biggl(\fint_{B_{1/2}} |V_{p}(\nabla w)|^2\dx\biggr)^\frac1{2}
\end{equation}
and there exists $C_M=C_M(\Lambda,M,n,p)>0$ and $\alpha=\alpha(\Lambda,M,n,p)\in(0,1)$ such that
\begin{equation}\label{est:L:reg:ahom:a2:holder}
[\nabla w]_{C^{0,\alpha}(B_{1/8})}\leq C_M\biggl(\fint_{B_{1/2}} |V_{p}(\nabla w)|^2\dx\biggr)^\frac1{2}.
\end{equation}
The differentiability of $\bfa$ and the bounds \eqref{ass:1:monotone:strong} and \eqref{ass:acont}  imply  
\begin{equation}\label{ass:Dbfa}
 \min\{1,(1+|z|)^{p-2}\}|\xi|^2\lesssim \partial\overline \bfa (z)\xi\cdot\xi,\qquad  |\partial\overline \bfa(z)\xi|\lesssim  \max\{1,(1+|z|)^{p-2})\}|\xi|.
\end{equation}
Note that by Theorem~\ref{thm:Lq}, we have that $\nabla w\in L_{\rm loc}^q(B)$ for every $q<\infty$. With help of the higher gradient integrability and the difference quotient method, it is easy to see that the partial derivatives $\partial_i w$, $i=1,\dots,n$ are weakly differentiable satisfying
\begin{equation}\label{pf:L:reg:hom:eq1}
\nabla \cdot (A\nabla \partial_iw)=0\quad\mbox{where}\quad A(x):=\partial \bfa(\zeta+\nabla w(x))
\end{equation} 
and $\min\{1,(1+|\zeta+\nabla w|)^{\frac{p-2}2}\}\nabla \partial_i w\in L_{\rm loc}^2(B)$. Now we are in position to apply local boundedness results for linear non-uniformly elliptic equation, e.g. \cite{T71,BS21}, to \eqref{pf:L:reg:hom:eq1}. In order to apply the results from \cite{BS21}, we observe that \eqref{ass:Dbfa} implies that the quantities
$$
\lambda(x):=\inf_{\xi\in\R^n\setminus\{0\}}\frac{\langle \xi,A(x)\xi\rangle}{|\xi|^2},\quad\mu(x):=\sup_{\xi\in\R^n\setminus\{0\}}\frac{|A(x)\xi|^2}{\langle \xi,A(x)\xi\rangle}
$$
satisfy
$$
1\lesssim \lambda,\quad \mu\lesssim (1+|\zeta+\nabla w|)^{2(p-2)}\mbox{ if $p\geq2$ and }(1+|\zeta+\nabla w|)^{p-2}\lesssim \lambda,\quad \mu\lesssim 1\mbox{ if $p\in(1,2)$.}
$$
Since $\nabla w\in L_{\rm loc}^q(B)$ for every $q<\infty$, we have $\lambda^{-1},\mu\in L^q_{\rm loc}(B)$ for every $q<\infty$ and thus we can apply \cite[Theorem~1.1]{BS21} with $\gamma=\min\{2,p\}$ to find $\Gamma=\Gamma(n,p)>0$ such that
\begin{equation}
\|\partial_i w\|_{L^\infty(B_{1/4})}\lesssim \biggl(\fint_{B_{1/2}}(1+|\zeta+\nabla w|)^{2n |p-2|}\dx\biggr)^\Gamma\biggl(\fint_{B_{1/2}} |\partial_i w|^{\min\{p,2\}}\dx\biggr)^\frac1{\min\{p,2\}}.
\end{equation}
The claimed estimate on $\|\nabla w\|_{L^\infty(B_{1/4})}$ in \eqref{est:L:reg:ahom:a2:linfty} in case $p\geq 2$ follows by summing over $i=1,\dots,n$, applying Theorem~\ref{thm:Lq} in the form
\begin{equation*}
\fint_{B_{1/2}}(1+|\zeta+\nabla w|)^{2n |p-2|}\dx\lesssim \biggl(\fint_{B_1}(1+|\zeta|+|\nabla w|)^p\dx\biggr)^\frac{2n|p-2|}p\lesssim M^\frac{2n|p-2|}p
\end{equation*}
and redefining $\Gamma$ (still depending only on $n$ and $p$). In the case $p\leq 2$, we additionally apply H\"older's inequality to estimate
\begin{align*}
\|\nabla w\|_{\underline L^p(B_{1/2})}\leq \|V_p(\nabla w)\|_{\underline L^2(B_{1/2})} \|1+|\nabla w|\|_{\underline L^p(B_{1/2})}^\frac{2-p} 2\lesssim \|V_p(\nabla w)\|_{\underline L^2(B_{1/2})} M^\frac{2-p}{2p}.
\end{align*}
 The claimed $H^2$- and $C^{0,\alpha}$ estimates of \eqref{est:L:reg:ahom:a2:linfty} and \eqref{est:L:reg:ahom:a2:holder} follow by standard elliptic regularity for the equation \eqref{pf:L:reg:hom:eq1} using the fact that by the $L^\infty$ bound on $\nabla w$  in \eqref{est:L:reg:ahom:a2:linfty} the coefficient $A$ is uniformly elliptic in $B_{1/4}$ with an ellipticity ratio that is bounded by a power of $M$.

\step 2 Conclusion. Step~1 and a simple scaling and covering argument imply the claimed estimates \eqref{reg:ahom:power} and \eqref{reg:ahom:holder} under the additional assumption that $\bfa\in C^1(\R^n,\R^n)$. In order to remove the differentiability assumption, we use a standard approximation argument. Let $\varphi$ be a non-negative, radially symmetric mollifier, i.e. it satisfies
$$
\varphi\geq0,\quad {\rm supp}\; \varphi\subset B_1,\quad \int_{\R^n}\varphi(x)\,dx=1,\quad \varphi(\cdot)=\widetilde \varphi(|\cdot|)\quad \mbox{for some $\widetilde\varphi\in C^\infty(\R)$}.
$$
For $\delta\in(0,1]$, set $\bfa_\delta:=\bfa\ast \varphi_\delta$, where $\varphi_\delta:=\delta^{-n}\varphi(\cdot /\delta)$. Clearly, $\bfa_\delta$ is smooth and satisfies \eqref{ass:1:monotone:strong} with the same constants and the continuity property \eqref{ass:acont} with the same constants in the case $p\in(1,2]$ and $1+|\xi_1|+|\xi_2|$ replaced by $3+|\xi_1|+|\xi_2|$ in the case $p>2$. Let $w_\delta\in W^{1,p}(B)$ be the unique solution to 
$$
\nabla \cdot \bfa_\delta(\zeta+\nabla w_\delta)=0\quad\mbox{in $B$}\quad w_\delta-w\in W_0^{1,p}(B).
$$
The energy estimate \eqref{L:energyestimate:eq1:basic:lip} implies
\begin{equation*}
|V_{p}(\zeta)|^2+\fint_{B}|V_{p}(\nabla w_\delta)|^2\dx\lesssim M \qquad\mbox{for all $\delta\in(0,1]$}.
\end{equation*}
Hence, $w_\delta$ satisfies \eqref{reg:ahom:power} and \eqref{reg:ahom:holder} for all $\delta\in(0,1]$. By standard monotonicity arguments (using $|\bfa_\delta(\xi)-\bfa(\xi)|\lesssim \delta^{\min\{1,p-1\}}$), we find $w_\delta\to w$ in $W^{1,p}(B)$. The desired estimate simply follow by standard weak lower semicontinuity results (and possibly replacing the H\"older exponent $\alpha$ by $\alpha/2$).
\end{proof}

\begin{remark}\label{lem:HighRegularity:rem}
The above argument can easily be extended to the case when \eqref{lem:HighRegularity:eq} is replaced by $\nabla \cdot \bfa (\nabla u)=f$ with $f\in L^q(B)$ with $q>n$. Indeed, the higher integrability also applies in that case, see Corollary~\ref{cor:Lqx} and for the replacement $\nabla \cdot A\nabla \partial_i w=\partial_i f$ of \eqref{pf:L:reg:hom:eq1}, we can apply similar local boundedness results valid for inhomogeneous equations, see \cite[Corollary 5.4]{T73}. This justifies regularity assumptions for the homogenized solution in \cite[Theorem~2.2]{Clozeau2023}.
\end{remark}

We next provide the two-scale homogenization estimate.

\begin{proposition}\label{Preg:1} Consider the situation in Theorem~\ref{L:nonlinear:largescalereg}.
Fix $M\in[1,\infty)$ and $\e\in(0,1]$.
 There exist $c=c(n,\lambda,p)\in[1,\infty)$, $\beta=\beta(p,n,\lambda,M)>0$, $\alpha=\alpha(M,n,p,\Lambda)>0$ and $C_M=C_M(\Lambda,M,n,p)>0$ such that the following is true: Let $u\in W^{1,p}(B_R)$ and $\zeta\in\R^n$ be such that
\begin{equation}\label{eq:P1Lip}
\nabla\cdot \bfa(\tfrac{x}\e,\nabla u)=0\qquad \text{ in } B_{R}\qquad\mbox{with}\qquad  |V_{p}(\zeta)|^2+\fint_{B_{R}}|V_{p}(\nabla u)|^2\dx\leq M.
\end{equation}
Then there exists $\xi\in\R^n$ satisfying
\begin{equation}\label{P1:est:xi-zetaLip}
|V_{p}(\xi-\zeta)|\leq C_M\biggl(\fint_{B_{R}}|V_{p}(\nabla u-\zeta-\nabla \phi_\zeta)|^2\dx\biggr)^\frac12
\end{equation}
such that for all $\theta\in(0,1)$ it holds
\begin{align}\label{P1:est:claimLip}
\fint_{B_{\theta R}}|V_{p}(\nabla u-(\xi+\nabla\phi_{\xi}))|^2\dx\leq&C_{M}\biggl(\theta^{2\alpha}+\theta^{-n}(\e/R)^\beta\biggr)\fint_{B_{R}}|V_{p}(\nabla u-\zeta-\nabla\phi_\zeta)|^2\dx\notag\\
&+C_M \theta^{-n}(\e/R)^\beta |V_p(\zeta)|^2.
\end{align}
Here $\phi_\zeta$ and $\phi_\xi$ are the correctors defined in Lemma \ref{L:phi}. 
\end{proposition}

\begin{proof}
We prove the statement for $R=2$. The general statement follows by a simple scaling argument. Throughout the proof, we write $\lesssim$ if $\leq$ holds up to a positive multiplicative constant depending on $(\Lambda,n,p)$ and denote by $\Gamma>0$ a constant depending only on $p$ and $n$ which might change in each occurrence.

\step 1 Harmonic approximation. Let $\overline u$ be the unique function satisfying
\begin{equation}\label{eq:uhom}
\nabla \cdot \overline \bfa(\nabla \overline u+\zeta)=0 \quad\text{ in } B_1,\qquad \overline u- u_\zeta\in W_0^{1,p}(B_1)\quad\mbox{where }u_\zeta(x):=u(x)-\zeta\cdot x-\phi_\zeta(x)-b,
\end{equation}
and $b\in\R$ is such that $(u_\zeta)_{B_2}=0$. Next, we gather the needed regularity properties for $\overline u$:
\begin{gather}\label{est:energy0:ubar}
\|V_{p}(\zeta+\nabla \overline u)\|_{\underline L^2(B_1)}\lesssim  \|V_{p}(\zeta+\nabla u_\zeta)\|_{\underline L^2(B_2)}\\
\forall \rho\in (0,\tfrac12)\qquad\|\nabla^2 \overline u\|_{L^{2}(B_{1-\rho})}+\|\nabla \overline u\|_{L^{\infty}(B_{1-\rho})}\lesssim (\rho^{-1}M)^\Gamma \|V_{p}(\nabla u_\zeta)\|_{\underline L^2(B_2)}\label{est:lip:linfty}\\
\|V_{p}(\zeta+\nabla \overline u)\|_{\underline L^{2m}(B_1)}\lesssim  \|V_{p}(\nabla u_\zeta)\|_{\underline L^2(B_2)}+|V_p(\zeta)|\label{est:nonlinearmeyerglob:lip1}
\end{gather}
and there exists $C_M=C_M(\Lambda,M,n,p)\in[1,\infty)$ and $\alpha=\alpha(\Lambda,M,n,p)>0$ such that
\begin{equation}\label{est:lip:holder:ubar}
|V_{p}([\nabla \overline u]_{C^{0,\alpha}(B_{1/4})})|\lesssim C_M \|V_{p}(\nabla u_\zeta)\|_{L^2(B_2)}.
\end{equation}

\eqref{est:energy0:ubar} follows directly from Proposition \ref{prop:regularity:basiclip} (see \eqref{L:energyestimate:eq1:basic:lip}). \eqref{est:nonlinearmeyerglob:lip} in combination with \eqref{eq:triangleV} gives
\begin{align}\label{eq:meyers1}
\|V_p(\zeta+\nabla \overline u)\|_{\underline L^{2m}(B_1)}\lesssim \|V_p(\zeta+\nabla u_\zeta)\|_{\underline L^{2m}(B_1)}\lesssim |V_p(\zeta)|+\|V_p(\nabla u)\|_{\underline L^{2m}(B_1)}+\|V_p(\nabla \phi_\zeta)\|_{\underline L^{2m}(B_1)}.
\end{align}
\eqref{est:nonlinearmeyerglob:lip1} now follows from \eqref{est:nonlinearmeyerloc:lip} applied to $u$ and $\phi_\zeta$ in combination with \eqref{eq:correctorControlled} and \eqref{eq:triangleV}. \eqref{est:lip:linfty} and \eqref{est:lip:holder:ubar} are a consequence of Prop \ref{lem:HighRegularity} in combination with the observation that
\begin{align}\label{eq:L2claim}
\|V_{p}(\nabla \overline u)\|_{\underline L^2(B_1)}\lesssim M^\frac{|p-2|}{\min\{p,2\}}\|V_{p}(\nabla u_\zeta)\|_{\underline L^2(B_2)}.
\end{align}
Indeed, using \eqref{eq:L2claim}, \eqref{est:lip:linfty} follows from \eqref{reg:ahom:power} by a simple covering argument. \eqref{est:lip:holder:ubar}  is a consequence of \eqref{reg:ahom:holder} and \eqref{eq:L2claim}, which give
\begin{align*}
|V_p([\nabla \overline u])_{C^{0,\alpha}(B_{1/4})}|^2\leq C_M |V_p(\|V_{p}(\nabla \overline u)\|_{L^2(B_{1/2})}^2)|^2\leq C_M |V_p(\|V_p(\nabla u_\zeta)\|_{L^2(B_1)}^2)|^2\leq C_M \|V_p(\nabla u_\zeta)\|_{L^2(B_2)}^2
\end{align*}
Thus, it remains to prove \eqref{eq:L2claim}. Note that using \eqref{ass:1:monotone:strong} and \eqref{ass:acont} for $p\geq 2$,
\begin{align*}
\fint_{B_1} |V_p(\nabla \overline u)|^2\dx\lesssim& \fint_{B_1} \langle \overline\bfa(\zeta+\nabla \overline u)-\overline\bfa(\zeta),\nabla \overline u\rangle\dx
\stackrel{\eqref{eq:uhom}}= \fint_{B_1} \langle \overline\bfa(\zeta+\nabla \overline u)-\overline\bfa(\zeta),\nabla u_\zeta\rangle\dx\\
\lesssim& \fint_{B_1} (1+|\nabla \overline u|+|\zeta|)^{p-2}|\nabla \overline u||\nabla u_\zeta|\dx\\
\lesssim&  \fint_{B_1} (1+|\nabla \overline u|+M^\frac12)^{p-2}|\nabla \overline u||\nabla u_\zeta|\dx.
\end{align*}
To obtain the last line, we noted that $|\zeta|\leq |V_p(\zeta)|\leq M^\frac 1 2$.
\eqref{eq:L2claim} now follows with the help of \eqref{eq:YoungV}. The case $p\leq 2$ is obtained by a similar calculation.

\step 2 Definition of two-scale expansion and energy estimate.

Fix $\rho>0$, $\ell \in\mathbb N$ with $\ell\geq2$ and define
\begin{equation}\label{def:mathcalQ:lip}
\mathcal Q:=\{Q_{\ell \e}(z)\,|:z\in \e\ell\Z^d,\,Q_{\ell \e}(z)\cap B_{1-2\rho}\neq\emptyset\}
\end{equation}
and set
\begin{equation}\label{def:overlineu2sLip}
\overline u^{2s}(x):=\overline u(x)+\zeta \cdot x+\sum_{Q\in\mathcal Q}\eta_Q\phi_{\zeta+\xi_Q,\e}
\end{equation}
where $\xi_Q$ and $\phi_{\xi,\e}$ are defined by \eqref{def:xik:lp} and $\eta_Q$ satisfies \eqref{def:etaQ} for every $Q\in\mathcal Q$. As in Step~2 of the proof of Lemma~\ref{Preg:1:lp}, we obtain
\begin{equation}\label{eq:errorLip}
-\nabla \cdot \bfa( \tfrac{\cdot}\e,\nabla \overline u^{2s})=\nabla \cdot \sum_{i=1}^3R^{(i)},
\end{equation}
where
\begin{align*}
R^{(1)}:=&\overline \bfa(\nabla \overline u+\zeta)-\sum_{Q\in\mathcal Q}\eta_Q\overline \bfa(\xi_Q+\zeta)\\
R^{(2)}:=&-\sum_{Q\in\mathcal Q} \sigma_{\xi_{Q}+\zeta,\e}\nabla \eta_{Q}\\
R^{(3)}:=&-\biggl(\bfa(\tfrac x \e,\nabla\overline u^{2s})-\sum_{Q\in\mathcal Q}\eta_Q\bfa(x,\xi_Q+\zeta+\nabla \phi_{\xi_{Q}+\zeta,\e})\biggr).
\end{align*}
Hence, the comparison estimate \eqref{L:energyestimate:eq2:basic:lip} implies that for every $\tau\in(0,1]$
\begin{align}\label{P:1step:claim:s2}
\int_{B_1}|V_{p}(\nabla u-\nabla \overline u^{2s})|^2\dx\lesssim& \tau\|V_{{p}}(\nabla \overline u^{2s})\|_{L^2(B_1)}^2+ (1+\tau^{-\frac{2-p}{p-1}})  \sum_{j=1}^3\|V_{{p'}}(R^{(j)})\|_{L^2(B_1)}^2.
\end{align}

\step 3 We claim that there exists $\beta=\beta(\Lambda,n,p)>0$ and a choice of $\rho\in(0,\frac14]$ and $\ell\in\mathbb N$ satisfying $\ell^{-1}\leq \e^{\beta}$ such that
\begin{align}\label{P:1step:claim:s2:1}
\int_{B_1}|V_{p}(\nabla u-\nabla \overline u^{2s})|^2\dx\lesssim&\e^\beta\left(|V_p(\zeta)|^2+  M^\Gamma \|\nabla u_\zeta\|_{L^2(B_2)}^2\right).
\end{align}
In order to prove \eqref{P:1step:claim:s2:1}, we estimate the right-hand side of \eqref{P:1step:claim:s2} term by term.

\substep{3.1} Estimate $\nabla \overline u^{2s}$. We claim that
\begin{equation}\label{est:nablau2s:lip}
\|V_{p}(\nabla \overline u^{2s})\|_{L^2(B_1)}^2\lesssim |V_p(\zeta)|^2+\|V_{p}(\nabla u_{\zeta})\|_{L^2(B_2)}^2.
\end{equation}
Firstly, we note that a similar computation as in \eqref{est:sumxiQtou:lp} yields
\begin{align}\label{est:V1pxiQ}
\int_{B_1}\sum_{Q\in \mathcal Q}\mathds 1_Q |V_{p}(\zeta+\xi_Q)|^2\dx\lesssim&\sum_{Q\in \mathcal Q}|Q| |V_{p}(\zeta+\xi_Q)|^2
\lesssim\int_{B_1}|V_{p}(\zeta+\nabla \overline u)|^2\dx.
\end{align}
Next, we compute
\begin{equation}\label{eq:nablau2s:lip}
\nabla \overline u^{2s}=\nabla \overline u+\zeta+\sum_{Q\in\mathcal Q}(\phi_{\xi_Q+\zeta}\nabla \eta_Q +\eta_Q \nabla \phi_{\xi_Q+\zeta}).
\end{equation}
Using, $\nabla\eta_Q\lesssim \e^{-1}(\mathds 1_Q-\mathds 1_{(Q)_{-\e}})$, \eqref{measure:1qetaq}, Lemma~\ref{lem:easyperiodic}, \eqref{eq:correctorControlled} and \eqref{est:V1pxiQ}, we obtain
\begin{align}\label{est:phinablaeta:lip}
\int_{B_1}\sum_{Q\in\mathcal Q}|V_{p}(\phi_{\xi_Q+\zeta,\e}\nabla \eta_Q)|^2\dx\lesssim&\sum_{Q\in\mathcal Q}\int_{Q\setminus(Q)_{-\e}}|V_{p}(\phi_{\xi_Q+\zeta}(\tfrac{x}\e))|^2\dx\notag\\
\lesssim&\sum_{Q\in\mathcal Q}\ell^{-1}|Q||V_{p}(\xi_Q+\zeta)|^2\lesssim\ell^{-1}\int_{B_1}|V_{p}(\zeta+\nabla \overline u)|^2\dx
\end{align}
and similarly (using $\eta_Q\leq\mathds 1_Q$)
\begin{align}\label{est:etanablaphin:lip}
\int_{B_1}\sum_{Q\in\mathcal Q}|V_{p}( \eta_Q\nabla \phi_{\xi_Q+\zeta,\e})|^2\dx\lesssim&\sum_{Q\in\mathcal Q}|Q||V_{p}(\xi_Q+\zeta)|^2\lesssim\int_{B_1}|V_{p}(\zeta+\nabla \overline u)|^2\dx.
\end{align}
For later reference, we also note that a similar computation as for \eqref{est:phinablaeta:lip} gives
\begin{align}\label{est:etanablaphin:lip:1}
\int_{B_1}\sum_{Q\in\mathcal Q}(\mathds 1_{Q}-\mathds 1_{(Q)_{-\e}})|V_{p}(\nabla \phi_{\xi_Q+\zeta,\e})|^2\dx\lesssim&\ell^{-1}\int_{B_1}|V_{p}(\zeta+\nabla \overline u)|^2\dx.
\end{align}

The claimed estimate \eqref{est:nablau2s:lip} follows from \eqref{eq:nablau2s:lip}, \eqref{est:phinablaeta:lip}, \eqref{est:etanablaphin:lip} and the energy estimate \eqref{est:energy0:ubar} in combination with the triangle inequality \eqref{eq:triangleV}.

\substep{3.2} The term $R^{(1)}$: We claim that there exists $\nu=\nu(n,\Lambda,p)\in(0,1]$ such that for any $\kappa\in(0,1]$,
\begin{align}\label{claim:est:R1Lip}
\int_{B_1}|V_{p^\prime}(R^{(1)})|^2\dx\lesssim&(\rho^\nu+\kappa+\ell^{-1})|V_{p}(\zeta)|^2\notag\\
&+\biggl((1+\kappa^\frac{2-p} p)(M\rho^{-1})^\Gamma(\ell\e)^{2}+\ell^{-1}+\kappa+\rho^\nu\biggr)\int_{B_2}|V_{p}(\nabla u_\zeta)|^2\dx.
\end{align}
Indeed,  the decomposition
\begin{align}\label{eq:decompR1}
|R^{(1)}|\leq |\overline\bfa(\nabla \overline u+\zeta)\biggl(1-\sum_{Q\in\mathcal Q}\mathds 1_Q\biggr)|+\sum_{Q\in \mathcal Q}\mathds 1_Q|\overline \bfa(\xi_Q+\zeta)-\overline\bfa (\nabla \overline u+\zeta)|+\sum_{Q\in\mathcal Q}(\mathds 1_{Q}-\eta_Q)|\overline\bfa (\xi_Q+\zeta)|
\end{align}
and the continuity of $\overline \bfa$ (see \eqref{ass:acont}) yields 
\begin{align*}
|R^{(1)}|\lesssim& 
 (\mathds 1_{B_{1}}-\mathds 1_{B_{1-2\rho}})(1+|\zeta+\nabla \overline u|)^{p-2}|\zeta+\nabla \overline u|\\
 &+\mathds 1_{\{p\geq2\}}\sum_{Q\in \mathcal Q} \mathds 1_Q(1+|\zeta+\nabla \overline u|+|\nabla \overline u-\xi_Q|)^{p-2}|\nabla \overline u-\xi_Q|\\
&+\mathds 1_{\{p\in(1,2)\}}\sum_{Q\in \mathcal Q}\mathds 1_Q(1+| \nabla \overline u-\xi_Q|)^{p-2}|\nabla\overline u-\xi_Q|+\sum_{Q\in \mathcal Q}(\mathds 1_Q-\eta_Q)(1+|\zeta+\xi_Q|)^{p-2}|\zeta+\xi_Q|.
\end{align*}
Inequality \eqref{est:V1pAB} in combination with \eqref{eq:triangleV} and \eqref{eq:VScaling} implies for all $\kappa\in(0,1]$
\begin{align}\label{est:V1pR1}
&\int_{B_1}|V_{p^\prime}(R^{(1)})|^2\dx\notag\\
\lesssim& \int_{B_1} \mathds 1_{B_1\setminus B_{1-2\rho}} |V_{p}(\zeta+\nabla \overline u)|^2+\sum_{Q\in \mathcal Q}(\mathds 1_Q-\eta_Q) |V_{p}(\zeta+\xi_Q)|^2\dx\notag\\
&+(1+\kappa^\frac{2-p}p)\int_{B_1}\sum_{Q\in \mathcal Q}\mathds 1_Q|V_{p}(\nabla\overline u-\xi_Q)|^2\dx+\kappa \mathds 1_{p>2}\int_{B_{1}}|\zeta+\nabla \overline u|^p\dx.
\end{align}
Next, we estimate the right-hand side in \eqref{est:V1pR1} term by term. By the higher integrability estimate \eqref{est:nonlinearmeyerglob:lip1}, we find $m>1$ such that
\begin{eqnarray}\label{est:boundarylayer:3rho}
\int_{B_1\setminus B_{(1-2\rho)}}|V_{p}(\zeta+\nabla \overline u)|^2\dx&\leq& |B_1\setminus B_{1-2\rho}|^{1-\frac1{ m}}\biggl(\int_{B_1}|V_{p}(\zeta+\nabla \overline u)|^{2m}\dx\biggr)^\frac{1}{ m}\notag\\
&\stackrel{\eqref{est:nonlinearmeyerglob:lip1}}\lesssim& \rho^{1-\frac1m}\left(|V_p(\zeta)|^2+\int_{B_2}|V_{p}(\nabla u_\zeta)|^2\dx\right).
\end{eqnarray}

The pointwise inequalities  $|V_{p}(\nabla \overline u-\xi_Q)|^2\lesssim   (1+\|\nabla \overline u\|_{L^\infty(B_{1-\rho})}^{p-2})|\nabla \overline u-\xi_Q|^2$ (for $p\geq2$ ) and the similar ${|V_{p}(\nabla \overline u-\xi_Q)|^2\lesssim |\nabla \overline u-\xi_Q|^2}$ (for $p\in(1,2)$) in combination with \eqref{est:lip:linfty} and Poincar\'e inequality imply
\begin{align}\label{ast:xikovu}
\int_{B_1}\sum_{Q\in\mathcal Q}\mathds 1_Q |V_{p}(\nabla \overline u-\xi_Q)|^2\dx\lesssim& (\rho^{-1}M)^\Gamma (\ell \e)^{2} \| \nabla^2 \overline u\|_{L^{2}(B_{1-\rho})}^2\notag\\
\lesssim& (\rho^{-1}M)^\Gamma(\ell \e)^2 \|V_{p}(\nabla u_\zeta)\|_{L^2(B_2)}^2.
\end{align}
Finally, using $|\mathds 1_Q-\eta_Q|\leq |\mathds 1_Q-\mathds 1_{(Q)_{-\e}}|$, \eqref{measure:1qetaq} and \eqref{est:V1pxiQ}, we obtain
\begin{align}\label{est:V1pzetadiff}
\int_{B_1}\sum_{Q\in \mathcal Q}(\mathds 1_Q-\eta_Q) |V_{p}(\zeta+\xi_Q)|^2\dx
\lesssim& \ell^{-1}\int_{B_1}|V_{p}(\zeta+\nabla \overline u)|^2\dx.
\end{align}
Combining \eqref{est:V1pR1}--\eqref{est:V1pzetadiff} with the energy estimate \eqref{est:energy0:ubar} and the triangle inequality \eqref{eq:triangleV}, we obtain \eqref{claim:est:R1Lip}.

\substep{3.3} The term $R^{(2)}$: Using the estimate \eqref{eq:correctorControlled} for $\sigma$ and similar computations as in \eqref{est:phinablaeta:lip}, we obtain
\begin{align*}
\int_{B_1}\sum_{Q\in \mathcal Q} |V_{p^\prime}(\sigma_{\xi_Q+\zeta,\e}\nabla \eta_Q)|^2\dx\lesssim& \ell^{-1}\int_{B_1}|V_{p}(\zeta+\nabla \overline u)|^2\dx
\end{align*}
and thus we obtain with help of \eqref{est:energy0:ubar} and the triangle inequality \eqref{eq:triangleV},
\begin{align}\label{claim:est:R2Lip}
\int_{B_1}|V_{{p'}}(|R^{(2)}|)|^2\dx\lesssim \ell^{-1}\left(|V_p(\zeta)|^2+\int_{B_2}|V_{p}(\nabla u_\zeta)|^2\dx\right).
\end{align}

\substep{3.4} The term $R^{(3)}$: We claim that there exists $\nu=\nu(n,\Lambda,p)\in(0,1]$ such that for any $\kappa\in(0,1]$,
\begin{align}\label{claim:est:R3Lip}
\int_{B_1}|V_{p^\prime}(R^{(3)})|^2\dx\lesssim&(\rho^\nu+\kappa+(1+\kappa^{\frac{2-p}p})\ell^{-1})|V_p(\zeta)|^2\notag\\
&+\biggl((1+\kappa^\frac{2-p} p)((M\rho^{-1})^\Gamma(\ell\e)^{2}+\ell^{-1})+\kappa+\rho^\nu\biggr)\int_{B_2}|V_{p}(\nabla u_\zeta)|^2\dx.
\end{align}
We decompose
\begin{align*}
|R^{(3)}|\leq& |\bfa(\tfrac{x}\e,\nabla \overline u+\zeta)||1-\sum_{Q\in\mathcal Q}\mathds 1_Q|+ \sum_{Q\in \mathcal Q}\mathds 1_Q |\bfa(\tfrac{x}\e,\nabla \overline u^{2s})-\bfa(\tfrac{x}\e,\xi_Q+\zeta+\nabla \phi_{\xi_Q+\zeta})|\\
&+\sum_{Q\in \mathcal Q} (\mathds 1_Q-\eta_Q)|\bfa(\tfrac{x}\e,\xi_Q+\zeta+\nabla \phi_{\xi_Q+\zeta})|.
\end{align*}
The continuity of $\bfa$ and
$$
\nabla \overline u^{2s}=\nabla \overline u+\zeta\quad\mbox{on $B_1\setminus \cup_{Q\in\mathcal Q}Q\subset B_1\setminus B_{1-2\rho}$}
$$
give
\begin{align*}
|R^{(3)}|\lesssim& (1-\mathds 1_{B_{1-2\rho}})(1+|\nabla \overline u+\zeta|)^{p-2}|\nabla \overline u+\zeta|\\
&+\mathds 1_{\{p\geq2\}}\sum_{Q\in\mathcal Q}\mathds 1_Q (1+|\nabla \overline u^{2s}|+|\xi_Q+\zeta+\nabla \phi_{\xi_Q+\zeta}|)^{p-2}|\nabla \overline u^{2s}-(\xi_Q+\zeta+\nabla \phi_{\xi_Q+\zeta})|\dx\\
&+\mathds 1_{\{p\in(1,2)\}}\sum_{Q\in\mathcal Q}\mathds 1_Q (1+|\nabla \overline u^{2s}-\xi_Q+\zeta+\nabla \phi_{\xi_Q+\zeta}|)^{p-2}|\nabla \overline u^{2s}-(\xi_Q+\zeta+\nabla \phi_{\xi_Q+\zeta})|\dx\\
&+ \sum_{Q\in \mathcal Q}(\mathds 1_Q-\eta_Q)(1+|\xi_Q+\zeta+\nabla \phi_{\xi_Q+\zeta}|)^{p-2}|\xi_Q+\zeta+\nabla \phi_{\xi_Q+\zeta}|\dx.
\end{align*}
The above estimate in combination with \eqref{eq:triangleV} and \eqref{est:V1pAB} implies for all $\kappa\in(0,1]$
\begin{align*}
\int_{B_1} |V_{p^\prime}(R^{(3)})|^2\dx\lesssim& \int_{B_1\setminus B_{1-2\rho}}|V_{p}(\zeta+\nabla \overline u)|^2\dx\\
&+\max\{1,\kappa^\frac{2-p} p\}\int_{B_1}\sum_{Q\in \mathcal Q} \mathds 1_Q |V_{p}(\nabla \overline u^{2s}-(\xi_Q+\zeta+\nabla \phi_{\xi_Q+\zeta}))|^2\dx\\
&+\kappa \int_{B_1}\sum_{Q\in\mathcal Q}\mathds 1_Q|V_{p}(\xi_Q+\zeta+\nabla \phi_{\xi_Q+\zeta})|^2\dx\\
&+\int_{B_1}\sum_{Q\in \mathcal Q} (\mathds 1_Q-\eta_Q)|V_{p}(\xi_Q+\zeta+\nabla \phi_{\xi_Q+\zeta})|^2\dx.
\end{align*}
The claimed estimate  \eqref{claim:est:R3Lip} follows from \eqref{est:boundarylayer:3rho} in combination with the  following estimates
\begin{align}\label{est:bigclaimR3:0}
\int_{B_1}\sum_{Q\in\mathcal Q}(\mathds 1_Q-\eta_Q)|V_{p}(\xi_Q+\zeta+\nabla \phi_{\xi_Q+\zeta})|^2\dx\lesssim \ell^{-1}\left(|V_p(\zeta)|^2+\int_{B_2}|V_{p}(\nabla u_\zeta)|^2\dx\right)
\end{align}
and
\begin{align}\label{est:bigclaimR3}
&\int_{B_1}\sum_{Q\in\mathcal Q} \mathds 1_Q|V_{p}(\nabla \overline u^{2s}-(\xi_Q+\zeta+\nabla \phi_{\xi_{Q}+\zeta}))|^2\dx\notag\\
\lesssim& \ell^{-1}|V_p(\zeta)|^2+\biggl((M\rho^{-1})^{\Gamma}\ell^2\e^2+\ell^{-1}\biggr)\int_{B_2}|V_{p}(\nabla u_\zeta)|^2\dx
\end{align}
Estimate \eqref{est:bigclaimR3:0} follows from $\mathds 1_{Q}\geq \eta_Q\geq\mathds 1_{Q_{-\e}}$, \eqref{eq:triangleV}, \eqref{est:V1pzetadiff}, \eqref{est:etanablaphin:lip:1} and the energy estimate \eqref{est:energy0:ubar}. Finally, we prove \eqref{est:bigclaimR3}. For this we use the identity
\begin{align*}
&\nabla \overline u^{2s}-(\xi_Q+\zeta+\nabla \phi_{\xi_{Q}+\zeta})\notag\\
=&(\nabla \overline u-\xi_Q)+(\eta_Q-\mathds 1_Q)\nabla \phi_{\xi_Q+\zeta}+\phi_{\xi_{Q}+\zeta}\nabla\eta_{Q}\qquad\mbox{on $Q\in\mathcal Q$}.
\end{align*}
Estimate \eqref{est:bigclaimR3} follows from \eqref{eq:triangleV} in combination with \eqref{ast:xikovu}, \eqref{est:etanablaphin:lip:1}, \eqref{est:phinablaeta:lip} and \eqref{est:energy0:ubar}. 

\substep{3.5} Proof of \eqref{P:1step:claim:s2:1}. We first consider $p\in[2,\infty)$. Sending $\tau\to0$ in \eqref{P:1step:claim:s2}, we obtain
$$
\int_{B_1}|V_{p}(\nabla u-\nabla \overline u^{2s})|^2\dx\lesssim \sum_{j=1}^3\|V_{{p'}}(R^{(j)})\|_{L^2(B_1)}^2
$$
and in combination with \eqref{claim:est:R1Lip}, \eqref{claim:est:R2Lip}, and \eqref{claim:est:R3Lip} with $\kappa=\ell^{-\frac{p}{2(p-1)}}\in(0,1]$, we have
$$
\int_{B_1}|V_{p}(\nabla u-\nabla u^{2s})|^2\dx\lesssim (\rho^\nu+\ell^{-\frac{p}{2(p-1)}})|V_p(\zeta)|^2+ \biggl(\ell^\frac{5p-6}{2(p-1)}(M\rho^{-1})^\Gamma \e^{2}+\ell^{-\frac{p}{2(p-1)}}+\rho^\nu\biggr)\int_{B_2}|V_{p}(\nabla u_\zeta)|^2\dx.
$$
Choosing $\rho\in(0,\frac14]$ the largest number such that $\rho^\nu\leq \ell^{-\frac{p}{2(p-1)}}$ and finally $\ell\geq10$ the smallest integer satisfying $\ell^{-\frac{p}{2(p-1)}}\leq \e^\beta$ with $\beta>0$ such that $\beta=2-\beta(\frac{5p-6}p+\frac\Gamma\nu)$, we obtain \eqref{P:1step:claim:s2:1} for $p\in[2,\infty)$. Next, we consider $p\in(1,2)$. Choosing $\rho\in(0,\frac14]$ the largest number such that $\rho^\nu\leq \ell^{-1}$ and sending $\kappa\to0$ in \eqref{claim:est:R1Lip} and \eqref{claim:est:R3Lip}, we obtain 
\begin{align*}
\int_{B_1}|V_{p}(R^{(1)})|^2+|V_{p}(R^{(3)})|^2\dx\lesssim& \ell^{-1}|V_p(\zeta)|^2+ \biggl(\ell^{2+\frac\Gamma\nu}M^\Gamma \e^{2}+\ell^{-1}\biggr)\int_{B_2}|V_{p}(\nabla u_\zeta)|^2\dx.
\end{align*}
Combining the above estimate with \eqref{P:1step:claim:s2}, \eqref{est:nablau2s:lip} and \eqref{claim:est:R2Lip}, we obtain for $\tau\in(0,1]$
\begin{align*}
\int_{B_1}|V_{p}(\nabla u-\nabla \overline u^{2s})|^2\dx\lesssim& \left(\tau+\tau^\frac{p-2}{p-1}\ell^{-1}\right)|V_p(\zeta)|^2\\
&+(\tau+\tau^{\frac{p-2}{p-1}}(\ell^{2+\frac{\Gamma}\nu}M^\Gamma \e^2+\ell^{-1}))\|V_{{p}}(\nabla \overline u_\zeta)\|_{L^2(B_2)}^2.
\end{align*}
Clearly, choosing first $\tau$ depending on $\ell$ and then $\ell$ depending on $\e$, we find $\beta>0$ such that \eqref{P:1step:claim:s2:1} holds true.

\step 4 We claim that there exists $C_M\in[1,\infty)$ and $\alpha=\alpha_M\in(0,1]$ such that for $\theta\in(0,\frac{1}8)$, we have  
\begin{align}\label{est:P1:s4}
\fint_{B_\theta}|V_{p}(\nabla u-(\xi_0+\zeta+\nabla \phi_{\xi_0+\zeta}))|^2\dx\lesssim&\theta^{-n}\fint_{B_1}|V_{p}(\nabla u-\nabla \overline u^{2s})|^2\dx\notag\\
&+C_M (\theta^{2\alpha}+\ell^{-1}+(\e \ell)^{2\alpha})\fint_{B_{2}}|V_{p}(\nabla u_\zeta)|^2\dx+\ell^{-1}|V_p(\zeta)|^2,
\end{align}
where $\xi_0:=\xi_{Q_{\ell \e}(0)}=(\nabla \overline u)_{Q_{\ell \e}(0)}$. By triangle inequality, we have
\begin{align}\label{est:P1:s4:triangle}
|\nabla u-(\xi_0+\zeta+\nabla \phi_{\xi_0+\zeta})|\leq&|\nabla u-\nabla u^{2s}|+|\nabla \overline u-\xi_0|+|\sum_{Q\in\mathcal Q} \mathds 1_Q(\nabla \phi_{\xi_Q+\zeta}-\nabla \phi_{\xi_0+\zeta})|\notag\\
&+\sum_{Q\in\mathcal Q} (\mathds 1_Q-\eta_Q)|\nabla \phi_{\xi_Q+\zeta}|+|\sum_{Q\in\mathcal Q} \phi_{\xi_Q+\zeta}\nabla \eta_Q|.
\end{align}
We first suppose that $\theta>0$ is sufficiently large that $Q_{\ell \e/2}(0)\subset B_\theta$. Set $Z_\theta:=\{Q\in\mathcal Q\,|\, Q\cap B_\theta \neq\emptyset\}$. In what follows, we frequently use \eqref{est:lip:holder:ubar} in the form
\begin{align*}
\|V_{p}(\nabla \overline u-\xi_0)\|_{L^\infty(B_\theta)}^2+\max_{Q\in Z_\theta}|V_{p}(\xi_Q-\xi_0)|^2\lesssim& |V_{p}(\theta^\alpha [\nabla \overline u]_{C^{0,\alpha}(B_\frac14)})|^2\\
\lesssim&C_M\theta^{\min\{p,2\}\alpha}\int_{B_2}|V_{p}(\nabla u_\zeta)|^2\dx.
\end{align*}
By the triangle inequality \eqref{eq:triangleV} in combination with \eqref{est:phinablaeta:lip} to estimate the fifth term, \eqref{est:etanablaphin:lip:1} to estimate the fourth, a similar computation for the difference of correctors using \eqref{eq:correctorDifference} (with $\tau=\theta^\alpha\in(0,1)$ in the case $p\in(1,2)$) to estimate the third term and $\#Z_\theta|Q|\lesssim \theta^n$ (here we use the assumption $Q_{\ell \e/2}(0)\subset B_\theta$), we have 
\begin{align*}
&\fint_{B_\theta}|V_{p}(\nabla u-(\xi_0+\zeta+\nabla \phi_{\xi_0+\zeta}))|^2\\
\lesssim&\fint_{B_\theta}|V_{p}(\nabla u-\nabla \overline u^{2s})|^2\dx+\|V_{p}(\nabla \overline u-\xi_0)\|_{L^\infty(B_\theta)}^2+\mathds 1_{\{p\geq2\}}\theta^{-n}\sum_{Q\in Z_\theta}|Q|(1+|\xi_Q|+|\xi_0|)^{p-2}|\xi_Q-\xi_0|^2\\
&+\mathds 1_{\{p\in(1,2)\}}\theta^{-n}\sum_{Q\in Z_\theta}|Q|(|V_{p}(\xi_Q-\xi_0)|^2\theta^{\frac{(p-2)\alpha}p}+\theta^\alpha(|V_{p}(\xi_Q)|^2+|V_{p}(\xi_0)|^2))\\
&+\theta^{-n}\sum_{Q\in Z_\theta}\ell^{-1}|Q| |V_p(\xi_Q+\zeta)|^2\\
\lesssim& \fint_{B_\theta} |V_p(\nabla u-\nabla u^{2s})|^2\dx + C_M \theta^{\min\{p,2\}\alpha}\left(1+\theta^\frac{\alpha(p-2)} p\right)\int_{B_2} |V_p(\nabla u_\zeta)|^2\dx\\
&+(\theta^\alpha+\ell^{-1}) \|V_p(\nabla \overline u)\|_{L^\infty(B_{1/4})}^2+\ell^{-1}|V_p(\zeta)|^2\\
\lesssim&\fint_{B_\theta}|V_{p}(\nabla u-\nabla u^{2s})|^2\dx+C_M \left(\theta^{\tilde \alpha}+\ell^{-1}\right)\int_{B_2}|V_{p}(\nabla u_\zeta)|^2+\ell^{-1}|V_p(\zeta)|^2,
\end{align*}
where we use in the second estimate $(1+|\xi_Q|+|\xi_0|)^{p-2}\lesssim (1+\|\nabla\overline u\|_{L^\infty(B_{1/4})})^{p-2}\lesssim C_M$ for $p\geq2$ and $Q\in Z_\theta$ and in the last estimate we set $\tilde\alpha = \min\{\alpha,\left(p+\frac{p-2} p\right)\alpha\}>0$. Hence, \eqref{est:P1:s4} follows by replacing $\alpha$ with $\tilde\alpha$. It remains to consider the case $B_\theta\subset Q_{\ell \e/2}(0)$. In this case, \eqref{est:P1:s4:triangle} in combination with \eqref{def:mathcalQ:lip} and \eqref{def:etaQ} imply
$$
|\nabla u-(\xi_0+\zeta+\nabla \phi_{\xi_0+\zeta})|\leq|\nabla u-\nabla u^{2s}|+|\nabla \overline u-\xi_0|\qquad\mbox{in $B_\theta$}.
$$ 
The claimed estimate \eqref{est:P1:s4} easily follows by using $B_\theta\subset Q_{\ell \e/2}(0)$ and $\xi_0=(\nabla \overline u)_{Q_{\ell \e}(0)}$ in the form
$$
\|V_{p}(\nabla \overline u-\xi_0)\|_{L^\infty(B_\theta)}^2\lesssim |V_{p}((\e \ell)^\alpha [\nabla \overline u]_{C^{0,\alpha}(B_\frac14)})|^2
\lesssim C_M(\e \ell)^{\min\{p,2\}\alpha}\int_{B_2}|V_{p}(\nabla u_\zeta)|^2\dx.
$$
\step 5 Conclusion. Consider the choices for $\rho$ and $\ell$ as in Step~3. Combining the estimates \eqref{P:1step:claim:s2:1} and \eqref{est:P1:s4} , we have
\begin{align*}
&\fint_{B_\theta}|V_{p}(\nabla u-(\xi_0+\zeta+\nabla \phi_{\xi_0+\zeta}))|^2\dx\lesssim C_M\left(\theta^{2\alpha}+\theta^{-n}\e^\beta +\e^\beta\right)\fint_{B_{2}}|V_{p}(\nabla u_\zeta)|^2\dx+\theta^{-n}\e^\beta |V_p(\zeta)|^2.
\end{align*}
for a suitable $\beta=\beta(\Lambda,M,n,p)>0$ as well as
$$
|V_{p}(\xi-\zeta)|^2=|V_{p}(\xi_0)|^2=|V_{p}(\|\nabla \overline u\|_{L^\infty(B_\frac14)})|^2\stackrel{\eqref{est:lip:linfty}}\lesssim  M^\Gamma \fint_{B_2}|V_{p}(\nabla  u_\zeta)|^2\dx
$$
and thus \eqref{P1:est:xi-zetaLip} follows.
\end{proof}

\subsection{Proof of Theorem~\ref{L:nonlinear:largescalereg}}\label{sec:prooflargescalelip}
We are now ready to provide the proof of Theorem~\ref{L:nonlinear:largescalereg}. Combining the regularity result Proposition~\ref{lem:HighRegularity} with Proposition~\ref{Preg:1}, we obtain a suitable decay estimate for a 'modified excess' which is usual in large-scale regularity results in homogenization.

\begin{proof}[Proof of Theorem~\ref{L:nonlinear:largescalereg}]
Throughout the proof we write $\lesssim$ if $\leq$ holds up to a positive multiplicative constant depending only on $\Lambda,n$ and $p$. 

We make the following preliminary observation: Let $c(p)\in[1,\infty)$ be the constant \eqref{eq:VEquiv}. For $z_1,z_2\in \R^n$ with $|z_1|\leq M'$, we have
\begin{align}\label{eq:triangleVprecise}
|V_p(z_2)|\stackrel{\eqref{eq:VEquiv}}\leq |V_p(z_1)|+c(p)|z_1-z_2|(1+|z_1|+|z_1-z_2|)^\frac{p-2} 2\leq |V_p(z_1)|+C_1(M')|V_p(z_1-z_2)|,
\end{align}
where $C_1(M') = 2c(p)\left(1+\mathds 1_{p\geq 2} (M')^\frac{p-2} 2\right)$.

\step 1 One-step improvement.

There exists  $c_1=c_1(\Lambda,n,p)\in[1,\infty)$ such that the following is true: For given $\overline M\in[1,\infty)$, let $\alpha=\alpha(\Lambda,\overline M,n,p)>0$ and $\beta=\beta(\Lambda,\overline M,n,p)\in(0,1]$ be as in Proposition~\ref{Preg:1}. For every $\tau\in(0,1]$ there exist $\theta_0=\theta_0(\Lambda,\overline M,n,p)\in(0,\frac14]$, $\e_0=\e_0(\Lambda,\overline M,n,p,\tau)\in(0,1]$,  and $C_{\overline M}=C_{\overline M}(\Lambda,\overline M,n,p)\in[1,\infty)$ with the following properties: Suppose $u\in W^{1,p}(B_R)$ and $\zeta\in\R^n$ satisfy
$$
\nabla \cdot \bfa(\tfrac{x}\e,\nabla u)=0\qquad\mbox{in $B_R$},
$$
and
\begin{equation}\label{T1:pf:eq0}
\fint_{B_{R}}|V_{p}(\nabla u-(\zeta+\nabla \phi_\zeta))|^2\dx\leq M'\leq \overline M\qquad\mbox{and}\qquad |V_{p}(\zeta)|^2\leq M'\leq {\overline M}.
\end{equation}
Then, there exists $\xi\in \R^n$ such that for every $\e\in(0,\e_0 R]$, it holds
\begin{align}\label{T1:pf:1}
\fint_{B_{\theta_0 R}}|V_{p}(\nabla u-(\xi+\nabla\phi_{\xi}))|^2\dx\leq& \theta_0^\alpha\|V_{p}(\nabla u-(\zeta+\nabla \phi_\zeta))\|_{\underline L^2(B_R)}^2+\tau (\e/R)^{\beta/2}|V_{p}(\zeta)|^2
\end{align}
and
\begin{align}\label{T1:pf:2}
|V_p(\xi-\zeta)|\leq& C_{M'}\|V_p(\nabla u-(\zeta+\nabla \phi_\zeta)\|_{\underline L^2(B_R)}.
\end{align}
%
%
 Moreover, we can choose $\theta_0$ such that it holds,
\begin{equation}\label{T1:pf:2b}
C_1(\overline M)C_{\overline M} \theta_0^{\alpha/2}\leq \frac12\quad\mbox{and}\quad \theta_0^\frac\beta8\leq\frac12.
\end{equation}
Triangle inequality \eqref{eq:triangleV} and the corrector estimates \eqref{eq:correctorControlled} imply that there exists $c=c(\Lambda,n,p)\in[1,\infty)$ with
$$
\fint_{B_{R}}|V_{p}(\nabla u)|^2\dx\leq c\fint_{B_{R}}|V_{p}(\nabla u-(\zeta+\nabla \phi_\zeta))|^2\dx+c|V_{p}(\zeta)|^2\leq 2cM'\leq 2c\overline M.
$$
By Proposition~\ref{Preg:1} (with $M=3cM'\leq 3c\overline M$) we obtain that there exists $\xi\in\R^n$ satisfying \eqref{T1:pf:2},
such that for all $\theta\in(0,1)$  it holds
%
\begin{align}\label{P1:est:claim:00}
\fint_{B_{\theta R}}|V_{p}(\nabla u-(\xi+\nabla\phi_{\xi}))|^2\dx\leq&C_{\overline M}\biggl(\theta^{2\alpha}+\theta^{-n}(\e/R)^\beta\biggr)\fint_{B_{R}}|V_{p}(\nabla u-\zeta-\nabla\phi_\zeta)|^2\dx\notag\\
&+C_{\overline M}\theta^{-n}(\e/R)^\beta |V_{p}(\zeta)|^2.
\end{align}
Let $\theta_0=\theta_0(\overline M,\alpha)\in(0,1]$ be the largest constant satisfying \eqref{T1:pf:2b}
and define $\e_0\in(0,\frac18]$ to be the largest constant satisfying
\begin{equation}\label{pf:T1:overeps}
2C_{\overline M}\theta^{-n}\e_0^\beta\leq \theta_0^\alpha,\quad C_1(\overline M)C_{\overline M}\theta_0^{-n} \e_0^\frac \beta 8\leq \tau 
\end{equation}
Inserting \eqref{T1:pf:2b} and \eqref{pf:T1:overeps} into \eqref{P1:est:claim:00}, we obtain \eqref{T1:pf:1}

\step 2 Iteration. Let $M\in[1,\infty)$ be given and choose $\overline M = 4M$. Let $0<\theta_0,\alpha,\beta$ and $\e_0\in(0,1]$ be as in Step~1 with $\tau=\tau(M)\in(0,1]$ given by
\begin{equation}\label{ass:iteration:tau}
\tau=\frac{1-\theta_0^\alpha}{4}\in(0,1].
\end{equation}
Set $R_k:=\theta_0^kR$ and $\zeta_0=0$. We claim that for all $k\in\mathbb N$ satisfying $\e\leq \e_0 \theta_0^kR$, we find $\zeta_k\in\R^n$ such that it holds
\begin{align}\label{T1:pf:s2:itera}\
&E_{k}\leq (\theta_0^{\alpha k}+(\e/R_{k-1})^{\beta/2})E_0\\ 
&|V_p(\zeta_k)|\leq \sum_{j=1}^{k-1}(2^{- j}+(\e/R_{j-1})^\frac{\beta}{8}) E_{0}^\frac12,\label{T1:pf:s2:iterb}
\end{align}
where 
\begin{equation}
E_{k}:=\fint_{B_{R_k}}|V_{p}(\nabla u-(\zeta_k+\nabla \phi_{\zeta_k}))|^2\dx.
\end{equation}
We proof the claim by induction. The case $k=1$ is a direct consequence of Step~1 with $M'= M$, $\zeta_{1}:=\xi$ and $\zeta=\zeta_0=0$. Suppose \eqref{T1:pf:s2:itera} and \eqref{T1:pf:s2:iterb} hold for some $1\leq k$. Then, we have for $\e\leq \e_0\theta^kR$
\begin{align}
E_{k}\leq& (\theta_0^{\alpha k}+(\e/R_{k-1})^{\beta/2})E_{0}\leq (\theta_0^\alpha+(\e_0\theta_0)^\frac{\beta}2)E_0\leq  (\theta_0^\alpha+(\e_0 \theta_0)^\frac{\beta}2)M\stackrel{\eqref{T1:pf:2b}}\leq \overline M\label{T1:pf:s2:iterb:1}\\
|V_p(\zeta_k)|\leq&  2E_0^\frac12\label{T1:pf:s2:iterb:2}\end{align}
where we use $\e\leq \e_0\theta_0^kR$, $R_j=\theta^jR$ and \eqref{T1:pf:2b} in the form of
$$
\sum_{j=1}^{k-1}(\e/R_{j-1})^\frac{\beta}{8}\leq \e_0^\frac{\beta}8\sum_{j=1}^{k-1}\theta_0^{(k-j+1)\frac\beta8}\leq\e_0^\frac\beta8\sum_{j=1}^\infty2^{-j}=\e_0^\frac\beta8\leq1.
$$
Suppose now, that $\e\in(0,\e_0 \theta_0^{k+1}R]$. Estimate \eqref{T1:pf:s2:iterb:2} and the choice of $\overline M$ imply that $|V_{p}(\zeta_k)|^2\leq \overline M$. Hence, we deduce from Step~1 with $M'= \overline M$, $\zeta_{k+1}:=\xi$ with $\zeta=\zeta_k$ that 
\begin{align}\label{T1:pf:s2:step1induction}
E_{k+1}\leq& \theta_0^\alpha (E_k+\tau (\e/R_k)^{\beta/2} |V_{p}(\zeta_k)|^2),\quad
|V_p(\zeta_{k+1}-\zeta_k)|\leq C_{\overline M} E_k^\frac12.
\end{align}
Appealing to the induction hypothesis, we find using the choice of $\tau$,
\begin{eqnarray*}
E_{k+1}&\leq& \theta_0^\alpha\biggl(\theta_0^k E_0+(\e/R_{k-1})^{\beta/2}E_0\biggr)+ \tau (\e/R_k)^{\beta/2} 4E_0 \\
&\stackrel{\eqref{ass:iteration:tau}}\leq&\theta_0^{\alpha(k+1)}E_0+(\e/R_k)^{\beta/2}E_0.
\end{eqnarray*}
Thus, we obtain the claim \eqref{T1:pf:s2:itera} with $k+1$. In order to show \eqref{T1:pf:s2:iterb} for $k+1$, we observe,
\begin{eqnarray*}
C_1(\overline M)|V_p(\zeta_{k+1}-\zeta_k)|&\stackrel{\eqref{T1:pf:s2:step1induction}}\leq& C_1(\overline M)C_{\overline M} E_k^\frac 1 2
\stackrel{\eqref{T1:pf:s2:iterb:1}}\leq C_1(\overline M)C_{\overline M}(\theta_0^\frac{\alpha k} 2+(\e/R_{k-1})^\frac \beta 4)E_0^\frac 1 2\\
&\leq &\left((C_1(\overline M)C_{\overline M}\theta_0^\frac \alpha 2)^k+
C_1(\overline M)C_{\overline M} (\e_0 \theta_0)^\frac \beta 8(\e/R_{k-1})^\frac \beta 8\right) E_0^\frac 1 2\\
&\stackrel{\eqref{T1:pf:2b},\eqref{pf:T1:overeps}}\leq &\left(2^{-k}+(\e/R_{k-1})^\frac \beta 8\right)E_0^\frac 1 2
\end{eqnarray*}
Note by elementary calculations that for $p\geq 2$, $|\zeta_k|\leq |V_p(\zeta_k)|\leq \overline M$. Hence, triangle inequality \eqref{eq:triangleVprecise} and \eqref{T1:pf:s2:iterb} (for $k$) imply \eqref{T1:pf:s2:iterb} with $k$ replaced by $k+1$. Clearly this finishes the induction argument.

\step 3 Conclusion. For given $M\in[1,\infty)$, let $\overline M=4M$ and $\alpha,\beta$ be as in Step~2. Let $u\in W^{1,p}(B_R)$ be such that \eqref{eq:T1} and \eqref{eq:T1:b} are satisfied. Fix $\ell\in\mathbb N$ as the largest integer such that  $\e\leq \e_0\theta^\ell R$. Let $r\in(\theta^{\ell+1}R,R]$ and let $k\in\mathbb N$ be the largest integer satisfying $r\leq \theta^kR$. Clearly, we have $k\geq \ell$ and thus $\e\leq \e_0\theta^kR$. In view of Step~2, there exists $\zeta_k\in \R^n$ satisfying \eqref{T1:pf:s2:iterb:2} and
\begin{align*}
\fint_{B_r}|V_{p}(\nabla u-(\zeta_k+\nabla \phi_{\zeta_k}))|^2\dx\leq&\theta_0^{-n}E_k\leq \theta_0^{-n}(\theta_0^{\alpha k}+(\e/R_{k-1})^{\beta/2})E_0\leq 2\theta_0^{-n}\fint_{B_{R}}|V_{p}(\nabla u)|^2\dx.
\end{align*}
Hence, triangle inequality and \eqref{T1:pf:s2:iterb:2}  imply
\begin{align*}
\fint_{B_r}|V_{p}(\nabla u)|^2\dx\lesssim \theta_0^{-n}\fint_{B_{R}}|V_{p}(\nabla u)|^2\dx+|V_{p}(\zeta_k)|^2\leq C\fint_{B_{R}}|V_{p}(\nabla u)|^2\dx,
\end{align*}
where $C=C(\Lambda,\overline M,n,p)\in[1,\infty)$. This finishes the proof.
\end{proof}

\begin{remark}\label{rem:rhslipschitz}
A straightforward variation of the above proof of Theorem~\ref{L:nonlinear:largescalereg} applies to the case when there is a small forcing term, that is \eqref{eq:T1} is replaced by $\nabla \bfa(\frac{x}\e,\nabla u)=f$, where $\|f\|_{L^q(B)}\leq \delta$ for some $\delta=\delta(\Lambda,M,n,p)>0$ and $q>n$. The smallness assumption is needed to counterbalance the dependency on $M$ in various estimates (in particular \eqref{P1:est:claimLip}). We do not expect that the smallness is necessary and it might by possible to use recent advances in Schauder-theory for non-uniformly elliptic equations \cite{deFilippis2023,deFilippis2024a} to overcome this problem.
\end{remark}

\textbf{Data availability statement:}
This paper has no associated data.

\appendix
\section{Proof of the energy estimates Proposition \ref{prop:regularity} and Proposition \ref{prop:regularity:basiclip}}\label{sec:appendixEnergyEstimates}
We provide the proofs of Proposition \ref{prop:regularity} and Proposition \ref{prop:regularity:basiclip} here.
\begin{proof}[Proof of Proposition \ref{prop:regularity}]
\step 1 Proof of the energy inequality. For $p\in[2,\infty)$ this follows directly from
\begin{align*}
\Lambda^{-1}\|\nabla w\|_{\underline L^p(B)}^p\leq& \fint_B \langle\overline\bfa(x,\nabla w)-\overline\bfa (x,0)),\nabla w\rangle\dx\\
=&\fint_B\langle\overline\bfa(x,\nabla w),\nabla u\rangle+\langle\overline \bfa(x,0),\nabla w\rangle\dx\leq \Lambda\fint_B(\mu+|\nabla w|)^{p-1}|\nabla u|+\mu^{p-1}|\nabla w|\dx
\end{align*}
and a suitable application of Young's inequality. In the case $p\in(1,2)$, we obtain with the same computation as in the case $p\in[2,\infty)$ together with the pointwise inequality
\begin{align}\label{ineq:pointwisesubq}
(\mu+|\nabla w|)^p\leq 2((\mu+|\nabla w|)^{p-2}|\nabla w|^2+(\mu+|\nabla w|)^{p-2}\mu^2)\leq 2(\mu+|\nabla w|)^{p-2}|\nabla w|^2+2\mu^p
\end{align}
the claim \eqref{L:energyestimate:eq1}.

\step 2 Proof of the comparison estimate. In case $p\in[2,\infty)$, we have
\begin{align*}
\Lambda^{-1}\|\nabla (u-w)\|_{\underline L^p(B)}^p\leq& \fint_{B}\langle\bfa(x,\nabla u)-\bfa(x,\nabla w)), \nabla (u-w)\rangle\dx\\
=&\fint_B \langle F,\nabla (u-w)\rangle+\langle\overline \bfa (x,\nabla w)-\bfa(x,\nabla w), \nabla (u-w)\rangle\dx\\
\leq& (\|F\|_{\underline L^{p'}(B)}+\delta^{p-1}\|\mu+|\nabla w|\|_{\underline L^p(B)}^{p-1})\|\nabla (u-w)\|_{\underline L^p(B)}
\end{align*}
which implies
\begin{equation*}
\|\nabla (u-w)\|_{\underline L^p(B)}\leq (\Lambda (\|F\|_{\underline L^{p'}(B)}+\delta^{p-1}\|\mu+|\nabla w|\|_{\underline L^p(B)}^{p-1}))^\frac1{p-1}
\end{equation*}
and thus \eqref{L:energyestimate:eq2} follows with help of \eqref{L:energyestimate:eq1}. Next, we consider $p\in(1,2)$. The same computations as above yield
\begin{align*}
\Lambda^{-1}\fint_B (\mu+|\nabla u|+|\nabla w|)^{p-2}|\nabla (u-w)|^2\dx\leq (\|F\|_{\underline L^{p'}(B)}+\delta^{p-1}\|\mu+|\nabla w|\|_{\underline L^p(B)}^{p-1})\|\nabla (u-w)\|_{\underline L^p(B)}.
\end{align*}
The above estimate in combination with H\"older inequality yields
\begin{align*}
\|\nabla (u-w)\|_{\underline L^p(B)}\leq& \|\mu+|\nabla u|+|\nabla w|\|_{\underline L^p(B)}^\frac{2-p}2\|(\mu+|\nabla u|+|\nabla w|)^\frac{p-2}2\nabla (u-w)\|_{\underline L^2(B)}\\
\leq&\|\mu+|\nabla u|+|\nabla w|\|_{\underline L^p(B)}^\frac{2-p}2(\Lambda (\|F\|_{\underline L^{p'}(B)}+\delta^{p-1}\|\mu+|\nabla w|\|_{\underline L^p(B)}^{p-1})^\frac 1 2\|\nabla (u-w)\|_{\underline L^p(B)})^\frac12
\end{align*}
Hence, we obtain with help of \eqref{L:energyestimate:eq1}
\begin{align*}
\|\nabla (u-w)\|_{\underline L^p(B)}\leq& C\|\mu+|\nabla u|\|_{\underline L^p(B)}^{2-p}(\|F\|_{\underline L^{p'}(B)}+\delta^{p-1}\|\mu+|\nabla u|\|_{\underline L^p(B)}^{p-1})
\end{align*}
and the claim \eqref{L:energyestimate:eq2} follows by Young's inequality in the form $a^{2-p}b\leq (2-p)\tau a+\tau^\frac{p-2}{p-1}b^\frac1{p-1}$.

\step 3 Caccioppoli inequality: By scaling and translation arguments we may assume that $B=B_1$. Take $\eta\in C^1_0(B_1)$ with $\eta =1$ in $B_{1/2}$ and $\lvert \nabla \eta\rvert \leq 8$. We first consider the case $p\geq2$. Testing \eqref{eq:basicRegularityEquation} with $\eta^p (u-b)$, we find
\begin{align*}
&\Lambda^{-1}2^{-p}\|\eta(\mu+|\nabla u|)\|_{L^p(B)}^p\\
\leq& \int_B \langle\bfa(x,\nabla u)-\bfa (x,0),\eta^p\nabla u\rangle+\eta^p\mu^p\dx\\
=&\int_B \langle F,\nabla (\eta^p(u-b))\rangle-p\langle\bfa(x,\nabla u),\eta^{p-1}\nabla \eta (u-b)\rangle-\langle\overline \bfa(x,0),\eta^p\nabla u\rangle+\eta^p\mu^p\dx\\
\leq& \|F\|_{L^{p'}(B)}\left(\|\eta\nabla u\|_{L^p(B)}+8p\|u-b\|_{L^p(B)}\right)+8p\Lambda\|\mu+|\nabla u|\|_{L^p(B)}^{p-1}\|u-b\|_{L^p(B)}\\
&+\int_B\eta^p(\mu^{p-1}|\nabla u|+\mu^p)\dx.
\end{align*}
The claimed estimate follows by suitable applications of Young's inequality. In the case $p\in(1,2)$, we use the same computation as above, together with the pointwise inequality
\begin{align*}
(\mu+|\nabla u|)^p\leq 2((\mu+|\nabla u|)^{p-2}|\nabla u|^2+(\mu+|\nabla u|)^{p-2}\mu^2)\leq 2\Lambda (\bfa(x,\nabla u)-\bfa (x,0))\cdot\nabla u+2\mu^p
\end{align*}
to deduce the claim.

\step 4 Proof of higher differentiability. By scaling and translation we may assume that $B=B_1$. Estimate \eqref{P1:pf:step1:claim0} established in the proof of Theorem~\ref{thm:Lq} in Section~\ref{sec:thm:Lqx}, for $\sigma=1$ reads
\begin{align*}
\|\nabla \tau_h w\|_{L^p(B_{r+(\rho-r)/4})}
\lesssim& \|\mu+|\nabla w|\|_{L^{p}(B_\rho)}^{\frac{\sigma-1+p}p}\begin{cases}(\frac{|h|}{\rho-r})^{\frac{1}{p-1}}&\mbox{if $p\in[2,\infty)$}\\(\frac{|h|}{\rho-r})&\mbox{if $p\in(1,2)$}\end{cases}
\end{align*}
Clearly, the above estimate in combination with \eqref{eq:defbesov} and the energy estimate \eqref{L:energyestimate:eq1} imply the claim.

\step 5 Proof of Meyer's estimates. The Caccioppoli inequality in combination with Poincar\'e-Sobolev inequality and Gehring's lemma \cite[Theorem 6.6.]{Giusti2003}, imply in a standard fashion, c.f. \cite{Giusti2003}, that there exist $m=m(\Lambda,n,p)>0$ and $c=c(\Lambda,n,p)$ such that for every $B_{r}(x)\subset B$ it holds
\begin{align*}
\|\mu+|\nabla u|\|_{\underline L^{mp}(B_{r/2}(x))}\leq c \|\mu+|\nabla u|\|_{\underline L^p(B_r(x))}.
\end{align*}
The claimed estimate \eqref{est:nonlinearmeyerloc:lp} follows by covering $\rho B=B_{\rho r}(x_0)$ with balls with radius $(1-\rho)r/4$ and using the above estimate. 

Combining the fact that a boundary Caccioppoli inequality can be derived for $w$ by a standard modification of the argument in Step 3 with boundary Gehring estimates, see \cite[Section 6.5]{Giusti2003} for both, we find \eqref{est:nonlinearmeyerglob:lp} with possibly different $m>1$ and $c>0$ but still depending only on $\Lambda,n$ and $p$.
\end{proof}

\begin{proof}[Proof of Proposition \ref{prop:regularity:basiclip}]
\step 1 Proof of the energy inequality. This follows directly from
\begin{align*}
\Lambda^{-1}\|V_{p}(\nabla w)\|_{\underline L^2(B)}^2\leq& \fint_B\langle\overline\bfa(x,\nabla w),\nabla w\rangle\dx=\fint_B\langle \overline\bfa(x,\nabla w),\nabla u\rangle\dx\leq \Lambda\fint_B(\mu+|\nabla w|)^{p-2}|\nabla w||\nabla u|\dx
\end{align*}
and \eqref{eq:YoungV}.

\step 2 Proof of the comparison estimate. For all $p\in(1,\infty)$, we have
\begin{align*}
\Lambda^{-1}\|W_{p}(\nabla u,\nabla w)\|_{\underline L^2(B)}^2\leq& \fint_{B}(\bfa(x,\nabla u)-\bfa(x,\nabla w))\cdot \nabla (u-w)\dx\\
=&\fint_B (F-(F)_B)\cdot\nabla (u-w)\dx.
\end{align*}
For $p\in[2,\infty)$, the above estimate together with \eqref{eq:YoungV:0} and \eqref{def:Wpz} imply \eqref{L:energyestimate:eq2:basic:lip}. For $p\in(1,2)$, we combine the above estimate with \eqref{eq:YoungV:0}, \eqref{est:WpVsub} and \eqref{L:energyestimate:eq1:basic:lip} and obtain for $\tau,\kappa \in(0,1]$
\begin{align*}
\|V_{p}(\nabla u-\nabla w)\|_{L^2(B)}\lesssim& \|W_{p}(\nabla u,\nabla w)\|_{\underline L^2(B)}\tau^{-\frac{2-p}{p}}+\tau\|V_{p}(\nabla u)\|_{L^2(B)}\\
\lesssim& \||F-a|\cdot |\nabla (u-w)|\|_{L^1(B)}^\frac12\tau^{-\frac{2-p}{p}}+\tau\|V_{p}(\nabla u)\|_{L^2(B)}\\
\lesssim&\kappa\tau^{-\frac{2-p}{p}} \|V_{p}(\nabla (u-w))\|_{L^2(B)}+\kappa^{-\frac1{(p-1)}}\tau^{-\frac{2-p}{p}}\|V_{p'}(F-a)\|_{L^2(B)}\\
&+\tau\|V_{p}(\nabla u)\|_{L^2(B)}.
\end{align*}
Choosing $\kappa=\theta \tau^{\frac{2-p}p}$ with $\theta\lesssim1$ sufficiently small, we can absorb the $\|V_{p}(\nabla u-\nabla w)\|_{ L^2(B)}$ term and the claimed estimate \eqref{L:energyestimate:eq2:basic:lip} follows.

\step 3 Proof of Caccioppoli inequality.
Take $\eta\in C_c^1(B)$ with $\eta=1$ in $B_{r/2}(x_0)$ and $|\nabla \eta|\leq \frac{8}r$. Testing \eqref{eq:basicRegularityEquation:basic:lip} with $\eta^p (u-b)$ and using \eqref{eq:YoungV},
\begin{eqnarray*}
\Lambda^{-1}\int \eta^p |V_{p}(\nabla u)|^2&\leq& \int \eta^p\langle \bfa(\nabla u),\nabla u\rangle\dx= -p \int \langle \bfa(\nabla u),\nabla \eta\rangle (u-b)\eta^{p-1}\dx\\
&\leq& 8p\Lambda \int (1+|\nabla u|)^{p-2}|\nabla u| |r^{-1}(u-b)|\eta^{p-1}\dx\\
&\stackrel{\eqref{eq:YoungV}}\leq& \frac1{2\Lambda} \int \eta^p|V_{p}(\nabla u)|^2\dx+C(\Lambda,p)\int_{B} |V_{p}(r^{-1}(u-b))|^2\dx.
\end{eqnarray*}
Re-arranging gives \eqref{est:cacc}.

\step 4  Proof of Meyer's estimates. 
The Caccioppoli inequality in combination with the Sobolev inequality for $V$-functions \eqref{eq:PoincareSobolev} and Gehring's lemma, see \cite{Giusti2003}, imply that there exist $m=m(\Lambda,n,p)>0$ and $c=c(\Lambda,n,p)$ such that for every $B_{r}(x)\subset B$ it holds
\begin{align*}
\|V_{p}(\nabla u)|\|_{\underline L^{2m}(B_{r/2}(x))}\leq c \|V_{p}(\nabla u)|\|_{\underline L^2(B_r(x))}.
\end{align*}
The claimed estimate \eqref{est:nonlinearmeyerloc:lip} follows by covering $\rho B=B_{\rho r}(x_0)$ with balls with radius $(1-\rho)r/4$ and using the above estimate. Combining the fact that a boundary Caccioppoli inequality can be derived for $w$ by a standard modification of the argument in Step 3 with boundary Gehring estimates, we find \eqref{est:nonlinearmeyerglob:lip} with possibly different $m>1$ and $c>0$ but still depending only on $\Lambda,n$ and $p$.
\end{proof}

\section{Standard controlled growth conditions}\label{sec:appendixB}

Finally, we state a result regarding local Lipschitz estimates under standard controlled growth assumptions that we used in the proof of Theorem~\ref{thm:Linfty:uniformeps}.

\begin{proposition}\label{eq:LipschitzStandard}
Suppose $\bfa:\R^n\times\R^n\to\R^n$ satisfies Assumption~\ref{ass:standard} for some $1<p<\infty$, $\Lambda\in[1,\infty)$ and $\mu=1$. Moreover, assume that \eqref{ass:x} holds with $\mu=1$ and $\omega(t)\leq \Lambda \max\{1,t^\alpha\}$ for some $\alpha\in (0,1)$. Let $\Omega\subset\R^n$ be open and $u\in W^{1,p}(\Omega)$ be a weak solution of
$$
\nabla \cdot\bfa(x,\nabla u)=0\qquad\mbox{in $\Omega$}.
$$
Then, there is $\gamma=\gamma(\alpha,n,\Lambda,p)\in (0,1)$ such that $u\in C^{1,\gamma}_{\rm loc}(\Omega)$. Moreover, there exists a constant $c=c(\alpha,n,\Lambda,p,R)\in[1,\infty)$ such that for all balls $B=B_R(x_0)\Subset\Omega$ it holds
\begin{equation}\label{eq:Lipschitzstandard}
\|V_p(\nabla u)\|_{L^\infty(\frac12 B)}^2+R^{-2\gamma}\|V_p(\nabla u)\|_{C^{0,\gamma}(\frac12 B)}^2\leq c  \fint_{B} |V_p(\nabla u)|^2\dx.
\end{equation}
\end{proposition}

Proposition~\ref{eq:LipschitzStandard} is certainly known, but we did not find an exact reference in the literature. Similar results can be found in classical textbooks, see e.g.\ \cite{Giusti2003}, or in more recent works in which the H\"older continuity of $x\mapsto \bfa(x,z)$ is relaxed, see e.g. \cite{DM11}. However in those works, estimate \eqref{eq:Lipschitzstandard} is given with an additional $+1$ on the right hand side which is a consequence of a slightly different continuity assumption compared to \eqref{ass:x} with $(\mu+|z|)^{p-2}|z|$ replaced by $(\mu+|z|)^{p-1}$.

\begin{proof}

Assuming that $B_1\Subset\Omega$, we show that there is $\gamma=\gamma(\alpha,\Lambda,n,p)>0$ and ${c=c(\alpha,\Lambda,n,p)\in[1,\infty)}$ such that for all $0<r<1$ it holds
\begin{equation}\label{claim:proof:standardregularity}
\fint_{B_r}|V_p(\nabla u)|^2\dx+r^{-2\gamma}\fint_{B_r}|V_p(\nabla u)-(V_p(\nabla u))_{B_r}|^2\dx\leq c \fint_{B_1}|V_p(\nabla u)|^2\dx.
\end{equation}
From this, the claimed regularity estimates easily follows by translation, scaling and covering arguments together with the Campanato characterization of H\"older continuity, see e.g. \cite[Section 5.1.1.]{Giaquinta}.

\step 1 Comparison estimate. Let $0<R\leq1$ and let $v\in W^{1,p}(B_R)$ be such that
\begin{align*}
v\in u+W^{1,p}_0(B_R)\qquad\mbox{and}\qquad \nabla \cdot \bfa(0,\nabla v) = 0 \quad \text{ in } B_R.
\end{align*}
Then, there exists $c=c(\Lambda,n,p)\in[1,\infty)$ such that
\begin{align}\label{eq:comparisonEstimate}
\int_{B_R} |V_p(\nabla u)-V_p(\nabla v)|^2\dx\leq c R^{\alpha p_2} \int_{B_R} |V_p(\nabla v)|^2\dx,
\end{align}
where we set $p_2:=\min\{p',2\}$. We first note that the same argument yielding \eqref{eq:YoungV:0} implies that there exists $c=c(p)>0$ such that for any $\tau>0$, $\mu>0$ and $z,w\in \R^n$,
\begin{align}\label{eq:YoungV:1}
(\mu+|z|)^{p-2}|z||w|\leq \tau |V_{\mu,p}(z)|^2 + c \max\{\tau^{-1} , \tau^{-(p-1)}\}|V_{\mu,p}(w)|^2,
\end{align}
where $V_{\mu,p}(z):=(\mu+|z|)^{\frac{p-2}2}z$. We have
\begin{align*}
&\Lambda^{-1}\int_{B_R} (1+|\nabla u|+|\nabla v|)^{p-2}|\nabla u-\nabla v|^2\dx\\
\leq& \int_{B_R} \langle \bfa(x,\nabla u)-\bfa(x,\nabla v),\nabla(u-v)\rangle\dx\\
=& \int_{B_R} \langle \bfa(0,\nabla v)-\bfa(x,\nabla v),\nabla(u-v)\rangle\dx\\
\leq& \Lambda R^\alpha\int_{B_R} (1+|\nabla v|)^{p-2}|\nabla v| |\nabla(u-v)|\dx\\
\leq& 2^{2-p}\Lambda R^\alpha \int_{B_R} (1+|\nabla v|+|\nabla v|)^{p-2}|\nabla v| |\nabla(u-v)|\dx.
\end{align*}
Using \eqref{eq:YoungV:1} with $\mu = 1+|\nabla v|$, $z = |\nabla v|$ and $w=|\nabla u-\nabla v|$, we find for any $\e=\tau^{-1}>0$,
\begin{align*}
(1+|\nabla v|+|\nabla v|)^{p-2}|\nabla v| |\nabla(u-v)|&\leq c\e |V_{\mu,p}(\nabla (u-v))|^2+\max\{\e^{-1},\e^{-\frac1{p-1}}\}|V_{\mu,p}(\nabla v)|^2,
\end{align*}
for some $c=c(p)>0$. Moreover, by the definition of $\mu$, we find $c=c(p)\in[1,\infty)$ such that 
$$
|V_{\mu,p}(\nabla v)|^2\leq c|V_{p}(\nabla v)|^2\quad\mbox{and}\quad  |V_{\mu,p}(\nabla (u-v))|^2\leq c(1+|\nabla u|+|\nabla v|)^{p-2}|\nabla (u-v)|^2.
$$
Combining the last three displayed formulas with the choice $\e=\frac{1}{C\Lambda^2R^\alpha}$ for a suitable constant ${C=C(p)\in[1,\infty)}$, we obtain \eqref{eq:comparisonEstimate}.

\step 2 Regularity for $v$. Let $0<R\leq1$ and $v$ be as in Step~1. Then there exists $\beta=\beta(\Lambda,n,p)>0$ and $c=c(\Lambda,n,p)\in[1,\infty)$ such that for all $0<r\leq R$ it holds
\begin{align}\label{eq:autonomousLipschitz}
\fint_{B_r}|V_p(\nabla v)|^2\dx\leq c \fint_{B_R} |V_p(\nabla u)|^2\dx,
\end{align}
and
\begin{align}\label{eq:autonomousHolder2}
\fint_{B_r} |V_{p}(\nabla v)-(V_{p}(\nabla v))_{B_r}|^2\dx\leq c \biggl(\frac{r}{R}\biggr)^{2\beta}\fint_{B_R} |V_p(\nabla v)-(V_p(\nabla v))_{B_R}|^2\dx.
\end{align}
The key is to obtain \eqref{eq:autonomousHolder2} from which \eqref{eq:autonomousLipschitz} will be deduced. \eqref{eq:autonomousHolder2} is certainly not new, but we are not aware of a precise reference. The estimate \eqref{eq:autonomousHolder2} is proven in \cite{DM10} under stronger assumptions, see also the references therein for further similar results. A proof of \eqref{eq:autonomousHolder2} can be obtained following \cite[Section 8.7]{Giusti2003}, choosing (in the notation of \cite{Giusti2003}) $w_0= (1+|\nabla v|^2)^\frac p 2-1$ and $w_k = \left((1+|\nabla v|)^{p-1}-1\right)\partial_k v$.

\step 3 Morrey-regularity for $u$. For every $\e>0$, there exists $c=c(\alpha,\e,\Lambda,n,p)\in[1,\infty)$ such that the following is true: For all $0<r<R\leq1$ it holds
\begin{align}\label{eq:MorreyReg}
\int_{B_r} |V_p(\nabla u)|^2\dx\leq c \left(\frac r R\right)^{n-\e}\int_{B_R} |V_p(\nabla u)|^2\dx.
\end{align}
Consider $0<r<R\leq R_0\leq1$ and let $v\in u+W_0^{1,p}(B_R)$ be as in Step~1. Using \eqref{eq:comparisonEstimate} and \eqref{eq:autonomousLipschitz} it holds
\begin{eqnarray*}
\int_{B_r} |V_p(\nabla u)|^2\dx &\leq& 2\int_{B_r} |V_p(\nabla u)-V_p(\nabla v)|^2\dx + 2\int_{B_r} |V_p(\nabla v)|^2\dx\\
&\stackrel{\eqref{eq:autonomousLipschitz},\eqref{eq:comparisonEstimate}}\leq& c\left(R_0^{\alpha p_2}+\left(\frac r R\right)^n\right)\int_{B_R} |V_p(\nabla u)|^2\dx.
\end{eqnarray*}
Appealing to a standard iteration argument, see e.g. \cite[Lemma 5.13]{Giaquinta}, we obtain that \eqref{eq:MorreyReg} holds for ${0<r<R\leq R_0}$ for some $R_0=R_0(\alpha,\e,\Lambda,n,p)\in(0,1]$. Clearly, from this the general claim follows.

\step 4 Conclusion. Set $\gamma=\min\{\frac12 p_2\alpha,\beta\}>0$. We claim that there exists $c=c(\alpha,\Lambda,n,p)\in[1,\infty)$ such that for all $0<r\leq1$ it holds
\begin{align}\label{est:step4:holder}
\fint_{B_r} |V_p(\nabla u)-(V_p(\nabla u))_{B_r}|^2\dx \leq cr^{\gamma}\int_{B_1} |V_p(\nabla u)|^2\dx.
\end{align}

Let $0<r<R\leq1$ and let $v\in u+W_0^{1,p}(B_R)$ be as in Step~1. By triangle inequality, we have
\begin{align*}
\int_{B_r} |V_p(\nabla u)-(V_p(\nabla u))_{B_r}|^2\dx \leq 2\int_{B_r} |V_p(\nabla v)-(V_p(\nabla v))_{B_r}|^2\dx + 2\int_{B_R} |V_p(\nabla u)-V_p(\nabla v)|^2\dx.
\end{align*}
For the first term, we use \eqref{eq:autonomousHolder2} and triangle inequality to obtain
\begin{align*}
&\int_{B_r} |V_p(\nabla v)-(V_p(\nabla v))_{B_r}|^2\dx \\
\leq& c\biggl(\frac{r}{R}\biggr)^{n+2\beta}\int_{B_R} |V_p(\nabla v)-(V_p(\nabla v))_{B_R}|^2\dx\\
\leq&2 c\biggl(\frac{r}{R}\biggr)^{n+2\beta}\int_{B_R} |V_p(\nabla u)-(V_p(\nabla u))_{B_R}|^2\dx+2c\int_{B_R}|V_p(\nabla u)-V_p(\nabla v)|^2\dx.
\end{align*}
The term involving $|V_p(\nabla u)-V_p(\nabla v)|$ can be estimated by
\begin{align*}
\int_{B_R} |V_p(\nabla u)-V_p(\nabla v)|^2\dx\stackrel{\eqref{eq:comparisonEstimate}}\leq c R^{\alpha p_2}\int_{B_R} |V_p(\nabla u)|^2\dx\leq c R^{n+\gamma}\int_{B_1} |V_p(\nabla u)|^2\dx,
\end{align*}
where we use \eqref{eq:MorreyReg} with $\e=\gamma$ and $R\leq1$ in the second estimate, and $c=c(\alpha,\Lambda,n,p)\in[1,\infty)$. Combining the previous estimates we find
\begin{align*}
\int_{B_r} |V_p(\nabla u)-(V_p(\nabla u))_{B_r}|^2\dx \leq& c\biggl(\frac{r}{R}\biggr)^{n+2\gamma}\int_{B_R} |V_p(\nabla u)-(V_p(\nabla u))_{B_R}|^2\dx \\
&\quad+ c R^{n+\gamma}\fint_{B_1} |V_p(\nabla u)|^2\dx
\end{align*}
and the claim \eqref{est:step4:holder} follows by iteration, see \cite[Lemma 5.13]{Giaquinta}. The remaining estimate claimed in \eqref{claim:proof:standardregularity} follows from \eqref{est:step4:holder} and a standard iteration argument using the fact that for all $0<\rho\leq r<2\rho\leq1$ it holds
\begin{align*}
|(V_p(\nabla u))_{B_{r}}-(V_p(\nabla u))_{B_{\rho}}|
\leq& 2^n\fint_{B_{r}}|V_p(\nabla u)-(V_p(\nabla u))_{B_{r}}|\dx+\fint_{B_{\rho}}|V_p(\nabla u)-(V_p(\nabla u))_{B_{\rho}}|\dx\\
\leq& c(2^n+1)r^{\gamma/2}\biggl(\fint_{B_1}|V_p(\nabla u)|^2\dx\biggr)^\frac12
\end{align*}
\end{proof}

\bibliographystyle{plain}
\bibliography{references}

\end{document}